\documentclass[11pt]{article}
\usepackage[english]{babel}
\usepackage[latin1]{inputenc}
\usepackage{exscale}
\usepackage[centertags]{amsmath}
\usepackage{amsfonts,amstext,amssymb,amsthm}
\usepackage{newlfont}
\usepackage{enumerate}
\usepackage{graphicx}
\usepackage{these}
\usepackage{color}
\usepackage[usenames,dvipsnames]{xcolor}
\definecolor{BlueFonse}{rgb}{0,0,1}
\definecolor{BlueFonse1}{cmyk}{1,0,0,0.7}
\usepackage[colorlinks=true,linkcolor=blue,citecolor=ForestGreen]{hyperref}
\usepackage{hyperref}
\usepackage{url}
\usepackage{esint}
\usepackage{tikz}
\usetikzlibrary{intersections}
\usepackage{pgfplots}
\usetikzlibrary{calc,3d,shapes, pgfplots.external, intersections}
\usepackage{pgfplots}
\usetikzlibrary{decorations.pathmorphing,snakes}
\usetikzlibrary{patterns}
\usepackage{float}
\DeclareMathOperator{\diam}{diam}
\title{Sharp isoperimetric inequalities for small volumes in complete noncompact Riemannian manifolds of bounded geometry involving the scalar curvature}
\author{Stefano Nardulli, Luis Eduardo Osorio Acevedo\footnote{Partially supported by CNPq}}
\begin{document}
      \maketitle
\begin{center}
\noindent {\sc abstract}. We provide an isoperimetric comparison theorem for small volumes in an $n$-dimensional Riemannian manifold $(M^n,g)$ with strong bounded geometry, as in Definition \ref{Def:StrongBoundedGeometry}, involving the scalar curvature function. Namely in strong bounded geometry, if the supremum of scalar curvature function $S_g<n(n-1)k_0$ for some $k_0\in\R$, then for small volumes the isoperimetric profile of $(M^n,g)$ is less then or equal to the isoperimetric profile of $\mathbb{M}^n_{k_0}$ the complete simply connected space form of constant sectional curvature $k_0$. This work generalizes Theorem $2$ of \cite{DruetPAMS} in which the same result was proved in the case where $(M^n, g)$ is assumed to be just compact. As a consequence of our result we give an asymptotic expansion in Puiseux's series up to the second nontrivial term of the isoperimetric profile function for small volumes. Finally, as a corollary of our isoperimetric comparison result, it is shown, in the special case of manifolds with strong bounded geometry, and $S_g<n(n-1)k_0$ that for small volumes the Aubin-Cartan-Hadamard's Conjecture in any dimension $n$ is true. 
\newline\newline
\noindent{\it Key Words:} Isoperimetric comparison, isoperimetric inequalities, small volumes, bounded geometry, Aubin-Cartan-Hadamard's conjecture, finite perimeter sets, metric geometry, calculus of variations, geometric measure theory, Sobolev's inequalities on manifolds, scalar curvature, partial differential equations on manifolds, smoothing Riemannian manifolds, Ricci flow.
\bigskip

\centerline{\bf AMS subject classification: }
49Q20, 58E99, 53A10, 53C23, 53C21, 49Q05, 49Q10, 49Q15, 53C42, 83C99.
\end{center}
      \tableofcontents       
      \newpage
\section{Introduction}\label{1}
For the rigorous technical meaning of the concepts involved in this informal introductory section we refer the reader to Sections \ref{Sec:MainResults} and \ref{Sec:Theorem1}. Let us start recalling what is the isoperimetric problem. The isoperimetric problem in a Riemannian manifold $(M^n, g)$ of dimension $n$, consists in the study of the existence and geometric characterization of isoperimetric regions (see Definition \ref{Def:IsPWeak}), i.e., domains which minimizes area under a fixed volume constraint. Together with this problem comes naturally the study of the \textit{isoperimetric profile} function $I_{M,g}$  assigning to every $v\in[0,V_g(M)[$ the infimum of the areas of open bounded sets with smooth boundary and fixed volume $v$ (which is well known to be equal to $I_M$ of Definition \ref{Def:IsPWeak}, as could be seen comparing with Theorem $1$ of \cite{FloresNardulli015}). A natural question arises when dealing with isoperimetric profiles.
\begin{Question} It is possible to compare the isoperimetric profile of differents Riemannian manifolds just comparing the way in which they are curved? 
\end{Question}
During the lasts decades various attempts were done to answer to this question and several partial results were obtained. We outline here the more relevant (at our knowledge) to the present work. It is well known that the isoperimetric profile $I_{\R^n}(v)=c_nv^{\frac{n-1}n}$, for some explicit positive constant $c_n$ depending only on the dimension $n$ of the manifold. We call $c_n$ the \textit{Euclidean isoperimetric constant}. Furthermore it is a classical result that the isoperimetric regions (see Definition \ref{Def:IsPWeak}) in Euclidean space are the geodesic balls. 
%Let us introduce the reader to this longstanding conjecture. 
%In $\R^n$ with the Euclidean metric $g_0$ it is true that for every finite perimeter set $\Omega$ we have that $\P_{g_0}(\Omega)\ge c_nV_{g_0}(\Omega)^{\frac{n-1}n}$. 
so the constant $c_n$ is known explicitly. %and was determined independently by Thierry Aubin in \cite{AubinJDG76} and Giorgio Talenti in \cite{Talenti}.  
The Aubin-Cartan-Hadamard's conjecture states that for any Cartan-Hadamard manifold $(M^n,g)$ of dimension $n$ (i.e., a simply connected manifold with sectional curvature less than or equal to $0$), every finite perimeter set $\Omega\subseteq M$ satisfies the inequality $\P_g(\Omega)\ge c_nV_g(\Omega)^{\frac{n-1}n}$. This conjecture was proved to be true if $n=2$ in \cite{Weil26}, if $n=3$ first in \cite{Kleiner92} and later in \cite{RitoreCartanHadamard05} and \cite{Schulze08}, if $n=4$ in \cite{Croke}, finally as a corollary of Theorem $2$ of \cite{DruetPAMS} a weaker form of the conjecture follows for any dimension $n$ assuming $M$ to be compact (of course when $M$ is compact we do not have that the manifold is simply connected and with nonpositive curvature) and $\Omega$ having small volume. 
It is worth noting that in \cite{Croke} it is shown that $\P_g(\Omega)\ge\tilde{c}_nV_g(\Omega)^{\frac{n-1}n}$ is true in any dimension but with a positive constant $\tilde{c}_n$ that is not sharp in general. However the proof of C. Croke gives in the special case $n=4$ that $\tilde{c}_4=c_4$. In dimension $n\ge 5$ the Aubin-Cartan-Hadamard conjecture is still open. At our knowledge the only previous partial results in any dimension $n$ with the sharp constant, but restricted to the small volume regime, are Theorem $4.3$ of \cite{MJ} which require additional assumptions on the Gauss-Bonnet-Chern integrand in even dimension, and Theorem $2$ of \cite{DruetPAMS} in case of compact manifolds and Corollary $2$ of \cite{FloresNardulli015} in case of noncompact manifolds with $C^2$-locally asymptotically bounded geometry at infinity (compare Definition \ref{Def:StrongBoundedGeometrySmoothInfinity}) that is the noncompact version of Theorem $4.4$ of \cite{MJ}, but this requires a bound on the sectional curvature. Our Corollary \ref{CorRes:AubinCartanHadamard} extends these partial results in any dimension to domains $\Omega$ of small volume inside a Cartan-Hadamard manifold $(M,g)$ having strong bounded geometry. The difference between our result and Theorem $2$ of \cite{DruetPAMS} is that we relax the assumption on the manifold $M$ of being compact and replace it by just requiring that $M$ have $C^2$-locally asymptotically bounded geometry at infinity (see Definition \ref{Def:StrongBoundedGeometrySmoothInfinity}) or requiring just strong bounded geometry but loosing a little bit on the strict inequality sign in \eqref{Eq:ResTheorem3.1}. For a more exhaustive treatment about the state of the art of the Aubin-Cartan-Hadamard's conjecture we suggest the reading of the very good surveys of Olivier Druet \cite{DruetSurvey} available online, Section $3.2$ of Manuel Ritor\'e in \cite{RitoreSinestrari}, and the very recent and interesting paper \cite{KloeckerKuperberg}. Let us state here Theorem $1$ of \cite{DruetPAMS}. 
\begin{Thm}[Theorem $1$ of \cite{DruetPAMS}]\label{Thm:DruetPAMS1} Let $(M^n,g)$ be a complete Riemannian manifold, with $n\ge 2$, $x\in M$ such that there exists $k_0\in\R$ satisfying $Sc_g(x)<n(n-1)k_0$. Then there exists $r_x>0$ such that for every finite perimeter set $\Omega$ contained in the geodesic ball of center $x$ and radius $r_x$, 
\begin{equation}\label{Eq:Thm1DruetPAMs}
\P_g(\Omega)>\P_{g_{k_0}}(B),
\end{equation}
where $B$ is a ball enclosing a volume $v$ in the model simply connected space form $(\mathbb{M}^n_{k_0}, g_{k_0})$ of constant sectional curvature $k_0$. 
\end{Thm}
It should be seen from the proof of the preceding result and Theorems \ref{Res:Theorem1}-\ref{Res:Theorem3.1} of this paper that a lower bound of the optimal $r_x$ is continuous with respect to $Sc_g(x)$ and $C^0$ convergence of metrics, so if $M$ is compact there exists $r:=\inf\{r_x:\:x\in M\}>0$ such that the conclusion of the Theorem \ref{Thm:DruetPAMS1} holds for any $\Omega$ contained in a ball of radius $r$. Unfortunately the radius $r_x$ could go to zero when $x$ tends to infinity in an arbitrary noncompact complete Riemannan manifold. Hence some extra assumptions on the geometry at infinity of $M$ are needed to allow us to find such a positive uniform lower bound $r$. Actually using the last equation at page 2353 of \cite{DruetPAMS} and reasoning by contradiction it appears evident from the proof that to have $C^2$-locally asymptotically strong bounded geometry smooth at infinity, gives such a lower bound. A necessary condition to have $r>0$ is that the volume of balls of a fixed radius for example $r/2$ does not vanish when the centers go to infinity. This is a non collapsing condition that for example follows assuming $Ricci_g\ge(n-1)kg$ for some $k\in\R$ and positivity of the injectivity radius. Thus it seems natural to make these assumptions in our Theorem \ref{Res:Theorem1}. Actually, in other parts of the proof we will need to strengthen a little more our assumptions on the geometry of $M$ and we are lead to assume that $M$ have strong bounded geometry in the sense of our Definition  \ref{Def:StrongBoundedGeometry}. To obtain our main result about small volumes, in first we prove a global isoperimetric comparison for small diameters in Theorem \ref{Res:Theorem1} when $M$ has $C^2$-locally asymptotically bounded geometry at infinity, then we refine it in Theorem \ref{Res:Theorem1.1} requiring only that $M$ have strong bounded geometry, using smoothing of the metric via Ricci flow. Our Theorem \ref{Res:Theorem1} replaces the local isoperimetric comparison given by Theorem $1$ of \cite{DruetPAMS} with a global one under the assumptions of $C^2$-locally asymptotically bounded geometry at infinity (see Definition \ref{Def:StrongBoundedGeometry}).   
Given granted the proof of Theorem \ref{Res:Theorem1} we then prove Theorem \ref{Res:Theorem3} by geometric measure theory and Gromov-Hausdorff convergence of manifolds. The proof of Theorem \ref{Res:Theorem1} goes along the same lines of Theorem $1$ of \cite{DruetPAMS}. So the main ingredients used in its proof are results about local optimal Sobolev inequalities in $W^{1,p}$ via PDE techniques when $p>1$ which are easier to obtain than when $p=1$. After the limit problem when $p\to 1^+$ is studied. These local optimal Sobolev inequalities in $W^{1,p}$ are combined with an asymptotic analysis of solutions of quasi-elliptic equations involving the $p$-Laplacian when the parameter $p\to 1^+$. The importance of the scalar curvature when studying sharp Sobolev inequalities on Riemannian manifolds was first observed by Olivier Druet in \cite{DruetJFA}, later by Hebey in \cite{Hebey02} and appears evident when deducing Theorem $1$ \cite{DruetPAMS} from Proposition $1$ of \cite{DruetPAMS}.
The modifications required to achieve our goals are nontrivial, so to make the paper self-contained we wrote the entire proof of Theorem \ref{Res:Theorem1} in Section \ref{Sec:Theorem1}. Theorem \ref{Res:Theorem3} is a consequence of Theorem  \ref{Res:Theorem1} using techniques of geometric measure theory, say the theory of sets of finite perimeter, comparison geometry, and Gromov-Hausdorff convergence of manifolds. The proof follows the scheme traced by the proof of Theorem $2$ of \cite{DruetPAMS}, however the required changes in the proof are highly nontrivial and original. The two main difficulties that are encountered when one tries to apply the proof of Theorem $2$ of \cite{DruetPAMS} (working only for compact manifolds) to our more general context consist in the fact that 
existence of isoperimetric regions for every volume in a noncompact Riemannian manifold is no longer guaranteed and that one needs to prove that isoperimetric regions of small volumes  are also of small diameter. 
%on a result of \cite{MJ} refined later in \cite{NarAnn} and \cite{NarCalcVar} that asserts that for small volumes isoperimetric regions are also of small diameter. As it is widely known it is no longer true that isoperimetric regions exist in just complete Riemannian manifolds. 
For an account of results on the problem of existence of isoperimetric regions (see Definition \ref{Def:IsPWeak}) in complete Riemannian manifolds the reader is referred to \cite{NarAsian}, \cite{MondinoNardulli} and the references therein. Our approach to solve this difficulty is to use the theory of generalized existence and generalized compactness developed by the first author in \cite{NarAsian}, \cite{FloresNardulliGenComp}, and replace genuine isoperimetric regions in $M$ by generalized isoperimetric regions lying in some pointed limit manifold. This is possible because by a compactness theorem of the theory of convergence of manifolds namely Theorem $76$ of \cite{Petersen} the hypotheses of Theorem $1$ of \cite{FloresNardulliGenComp} are automatically fulfilled in the context of $C^2$-locally asymptotic strong bounded geometry smooth at infinity that we consider here. To finish the proof of Theorem \ref{Res:Theorem3} we need to prove that in $C^2$-locally asymptotically strong bounded geometry (compare Definition \ref{Def:StrongBoundedGeometry}) for an isoperimetric region having small volumes implies small diameter. With this aim in mind we replace the proof of \cite{MJ} based on Nash's isometric embeddings by another intrinsic one. We carry out this task proving a little more general result in Lemma \ref{Lemma:smallvolumesimpliessmalldiametersmild}, which asserts that if just $Ricci$ is bounded below and the injectivity radius is positive, isoperimetric regions of small volumes are of small diameter. In this proof we don't need to use any monotonicity formula; this fact constitutes a novelty with respect to the classical extrinsic proof of \cite{MJ}. Our proof is completely intrinsic and uses a cut and paste argument inspired by Proposition $2.5$ of \cite{NarCalcVar} (which works only for manifolds of strong bounded geometry) adapted to the case of weak bounded geometry (see Definition \ref{Def:BoundedGeometry}) joint with others non trivial intrinsic arguments aimed to encompass some technical difficulties of geometric measure theory, which arise when passing from the \textit{Euclidean space} $\R^n$ to an arbitrary Riemannian manifold $(M^n,g)$ without using Nash's isometric embedding theorem. The arguments of the proof permit also to give an effective estimate of the constants of  Lemma \ref{Lemma:smallvolumesimpliessmalldiametersmild} as functions of the bounds of the geometry of $(M^n,g)$. The main result of this paper is Theorem \ref{Res:Theorem3.1}. To prove it we need to prove beforehand Theorem \ref{Res:Theorem3} and then apply Theorem \ref{Res:Theorem3} to a suitable smoothing $(M,g_t)$ with initial data $(M,g)$ along the Ricci flow. The main reason is that for every $m\in\N=\{0,1,2,...\}$ the smoothed manifolds $(M,g_t)$ have $C^{m,\alpha}$ bounded geometry in the sense of Definition \ref{Def:BoundedGeometryPetersen} and so at infinity are more regular than the original $(M,g)$. This permits the use of the arguments of the proof of Theorem \ref{Res:Theorem1} to $(M,g_t)$ which does not work if applied directly to $(M,g)$. Hence passing to the limit when $t\to0^+$, and using the results of \cite{Shi}, \cite{Kapovitch}, we can transport the isoperimetric comparison from $(M,g_t)$ to $(M,g)$ without any further difficulty. As corollaries of our main Theorem  \ref{Res:Theorem3.1} we get immediately Corollary \ref{CorRes:AubinCartanHadamard1} that is a special case of the Aubin-Cartan-Hadamard's conjecture and the expansion of the isoperimetric profile in Puiseux's series given by Corollary \ref{CorRes:PuiseuxSeries}. As a final remark we have that all the constants involved in our statements of Section \ref{Sec:MainResults} are effectively computed in terms of the minimal bounds on the geometry that we are assuming. 
\section{Main results}\label{Sec:MainResults}
In the sequel we always assume that all the Riemannian manifolds $M^n$ considered are smooth with smooth Riemannian metric $g$. We denote by $V_g$ the canonical Riemannian measure induced on $M$ by $g$, and by $A_g$ the $(n-1)$-Hausdorff measure associated to the canonical Riemannian length space metric $d$ of $M$, that we also denote by $\mathcal{H}^{n-1}_g$. When it is already clear from the context, explicit mention of the metric $g$ will be suppressed. We will denote by $Ric_g$ the Ricci tensor of $(M,g)$, by $Sec_g$ the sectional curvature of $(M,g)$, $Sc_g$ the scalar curvature function, $S_g:=\sup_{x\in M}\left\{Sc_g(x)\right\}$ and by $\mathbb{M}_k^n$ the simply connected space form endowed with the standard metric of constant sectional curvature $k\in\R$ that we denote by $g_k$, by $inj_{(M,g)}$ the injectivity radius of $M$, for any $D\subseteq M$, $diam_g(D)$ the diameter of $D$ in the metric space $(M,g)$, $dv_g$ the Riemannian measure with respect to the metric $g$. In what follows we will consider as a key object the set of all finite perimeter sets (see Definition \ref{Def:LocallyFinitePerimeterSets}) of $M$ that we will denote by $\tilde{\tau}$. So a little technical discussion is in order here. By classical results of geometric measure theory (see Proposition $12.19$ and Formula $(15.3)$ of \cite{Maggi}) we know that if $E$ is a set of locally finite perimeter in $M$, then $spt(\nabla\chi_E)=\{x\in M:\:0<V_g(E\cap B(x,r)),\forall r>0\}\subseteq\partial E$, furthermore there exists an equivalent Borel set $F$ (i.e., $V_g(E\Delta F)=0$) such that $spt(\nabla\chi_F)=\partial F=\overline{\partial^*F}$, where $\partial^*F$ is the reduced boundary of $F$. It is not too hard to show that if $E$ has $C^1$ boundary, then $\partial^*E=\partial E$, where $\partial E$ is the topological boundary of $E$. De Giorgi's structure theorem (compare Theorem $15.9$ of \cite{Maggi}) guarantees that for every set $E$ of locally finite perimeter, $A_g(\partial^* E)=\mathcal{H}_g^{n-1}(\partial^*E)=\mathcal{P}_g(E)$. Hence without loss of generality we will adopt the assumption that all the locally finite perimeter sets considered in this text satisfy $\overline{\partial^*E}=\partial E$. It is worth to mention that the results in the book \cite{Maggi} are stated and proved in $\R^n$ but they are valid mutatis mutandis also in an arbitrary complete Riemannian manifold, the required details could be easily provided using the work about $BV$-functions on a Riemannian manifold accomplished in \cite{MPPP}.  
\begin{Def}\label{Def:BoundedGeometry}
A complete Riemannian manifold $(M, g)$, is said to have \textbf{weak bounded geometry}, if there exists a constant $k\in\mathbb{R}$, such that $Ric_M\geq k(n-1)$ (i.e., $Ric_M\geq k(n-1)g$ in the sense of quadratic forms) and $V(B_{(M,g)}(p,1))\geq v_0>0$, for some positive constant $v_0$, where $B_{(M,g)}(p,r)$ is the geodesic ball of $M$ centered at $p$ and of radius $r> 0$.
\end{Def}
\begin{Rem} In this paper we differ from the nomenclature used by the first author in his preceding works. What we call here weak bounded geometry is what is called, in all previous articles, just bounded geometry. 
\end{Rem}
\begin{Def}\label{Def:MildBoundedGeometry}
A complete Riemannian manifold $(M, g)$, is said to have \textbf{mild bounded geometry}, if there exists a constant $k\in\mathbb{R}$, such that $Ric_M\geq k(n-1)$ (i.e., $Ric_M\geq k(n-1)g$ in the sense of quadratic forms) and $inj_M>0$, where $inj_M$ is the injectivity radius of $M$. 
\end{Def}
\begin{Rem}
It is known that mild bounded geometry implies weak bounded geometry, but the converse is not true. For more details about this point the reader is referred to Remark $2.5$ of \cite{MondinoNardulli} and to the references therein. 
\end{Rem}
\begin{Def}\label{Def:StrongBoundedGeometry}
A complete Riemannian manifold $(M, g)$, is said to have \textbf{strong bounded geometry}, if there exists a positive constant $K>0$, such that $|Sec_M|\le K$ and $inj_M\ge i_0>0$ for some positive constant $i_0$. Sometimes we will use the condition $\Lambda_1\le Sec_M\le\Lambda_2$, for some given constants $\Lambda_1,\Lambda_2\in\R$ instead of $|Sec_M|\le K$ to express that $M$ have a two sided bound on the sectional curvature.
\end{Def}
\begin{Rem}
It turns out that it is easy to check that strong bounded geometry implies mild bounded geometry, with the converse being not true in general. 
\end{Rem}
%\begin{Rem} Observe that an analysis of the proof of Theorem \ref{Res:Theorem2} shows that $\tilde{v}_0$ could go to $0$ when $\varepsilon\to0$.
%\end{Rem}
\begin{Def} For any $m\in\mathbb{N}$, $\alpha\in [0, 1]$, a sequence of pointed smooth complete Riemannian manifolds is said to \textbf{converge in
the pointed $C^{m,\alpha}$, respectively $C^{m}$ topology to a smooth manifold $M$} $($denoted by $(M_i, p_i, g_i)\rightarrow (M,p,g)$$)$, if for every $R > 0$ we can find a domain $\Omega_R$ with $B(p,R)\subseteq\Omega_R\subseteq M$, a natural number $\nu_R\in\mathbb{N}$, and $C^{m+1}$ embeddings $F_{i,R}:\Omega_R\rightarrow M_i$, for large $i\geq\nu_R$ such that $B(p_i,R)\subseteq F_{i,R} (\Omega_R)$ and $F_{i,R}^*(g_i)\rightarrow g$ on $\Omega_R$ in the $C^{m,\alpha}$, respectively $C^m$ topology. 
\end{Def}
\begin{Def}\label{Def:BoundedGeometryInfinity} For any $m\in\N$ and $\alpha\in[0,1]$, we say that a smooth Riemannian manifold $(M^n, g)$ has $C^{m,\alpha}$-\textbf{locally asymptotic bounded geometry}, if it is with weak bounded geometry and if for every diverging sequence of points $(p_j)$, there exists a subsequence $(p_{{j}_{l}})$ and a pointed smooth manifold $(M_{\infty}, g_{\infty}, p_{\infty})$ with $g_{\infty}$ a smooth Riemannian metric such that the $C^{m,\alpha}$ norm is finite and the sequence of pointed manifolds $(M, p_{{j}_{l}}, g)\rightarrow (M_{\infty}, g_{\infty}, p_{\infty})$, in  $C^{m,\alpha}$-topology. When $\alpha=0$ we write $C^m$ instead of $C^{m,0}$. 
\end{Def}
\begin{Def}[Page 308 of \cite{Petersen}]\label{Def:BoundedGeometryPetersen}
A subset $A$ of a Riemannian n-manifold $M$  has \textbf{bounded $C^{m,\alpha}$ norm on the scale of $r$}, $||A||_{C^{m,\alpha},r}\leq Q$, if every point p of M lies in an open set $U$ with a chart $\psi$ from the Euclidean $r$-ball into $U$ such that
\begin{enumerate}[(i):]
      \item For all $p\in A$ there exists $U$ such that $B(p,\frac{1}{10}e^{-Q}r)\subseteq U$.
      \item $|D\psi|\leq e^Q$ on $B(0,r)$ and $|D\psi^{-1}|\leq e^Q$ on $U$.
      \item $r^{|j|+\alpha}||D^jg||_{\alpha}\leq Q$ for all multi indices j with $0\leq |j|\leq m$, where $g$ is the matrix of functions of metric coefficients in the $\psi$ coordinates regarded as a matrix on $B(0,r)$. 
\end{enumerate} 
We write that $(M,g,p)\in\mathcal{M}^{m,\alpha}(n, Q, r)$, if $||M||_{C^{m,\alpha},r}\leq Q$. 
\end{Def}
%\begin{Rem} If $M$ have $C^2$-locally asymptotically bounded geometry and positive injectivity radius then $M$ have strong bounded geometry. But the converse does not hold. 
%\end{Rem}
\begin{Rem} The condition of being smooth at infinity is used just in the last equation \eqref{Eq:ResTheorem3Last} of the proof of Theorem \ref{Res:Theorem3} when we apply Theorem \ref{Res:Theorem1} to a possibly limit manifold $(M_{\infty}, g_{\infty})$ that even in strong bounded geometry is a $C^{3,\beta}$ differentiable manifold but with a metric that is just $C^{1,\beta}$ and no more regular. There are examples of this phenomenon as explained in Example $1.8$ of \cite{Peters}. Actually the limit metric is $W^{2,p}$ for any $p>1$, as showed in \cite{Nikolaev}. This last regularity result is not enough strong to allow the use the arguments of the proof of Theorem \ref{Res:Theorem1} in $(M_{\infty},g_\infty)$.
\end{Rem}
This last remark justifies the following definitions.
\begin{Def}\label{Def:StrongBoundedGeometrySmoothInfinity}
We say that a smooth Riemannian manifold $(M^n, g)$ is \textbf{smooth at infinity}, if for every diverging sequence of points $(p_j)$, there exists a subsequence $(p_{{j}_{l}})$ and a pointed smooth manifold $(M_{\infty}, g_{\infty}, p_{\infty})$ with $g_{\infty}$ of class $C^{\infty}$. We say that a smooth Riemannian manifold $(M^n, g)$ has \textbf{strong bounded geometry smooth at infinity}, if it is of strong bounded geometry and is smooth at infinity. We say that $(M^n, g)$ has $C^{m,\alpha}$-locally asymptotically strong bounded geometry smooth at infinity, if it is of strong bounded geometry, smooth at infinity, and has  $C^{m,\alpha}$-\textbf{locally asymptotic bounded geometry}.   
\end{Def}
\begin{Rem} Observe that by Theorems $76$ and $72$ of \cite{Petersen} or Theorem $4.4$ of \cite{Peters} it is easily seen that to have strong bounded geometry smooth at infinity implies to have $C^{1,\beta}$-locally asymptotic bounded geometry, for any $\beta$. 
\end{Rem} 
We have now all the definitions needed to state our results.
\begin{Res}[Small diameters in $C^2$-locally asymptotically strong bounded geometry smooth at infinity]\label{Res:Theorem1} Let $(M^n, g)$ be a complete Riemannian manifold with $n\geq 2$ and with $C^2$-locally asymptotically strong bounded geometry smooth at infinity. Let us assume that there exists a real constant $k_0\in\R$ such that $S_g<n(n-1)k_0$. Then there exists $d=d(n, k, k_0, inj_M, S_g)>0$, which depends only on $n, k, k_0, inj_M, S_g$ such that for every $\Omega\subseteq M^n$ finite perimeter set with diameter $\diam_g(\Omega)\leq d$ holds 
\begin{equation}\label{Eq:Theorem1}
\P_g(\Omega)>\P_{g_{k_0}}(B),
\end{equation}
where $B\subseteq\mathbb{M}_{k_0}^n$ is a geodesic ball having $V_{g_{k_0}}(B)=V_g(\Omega)$. Moreover we have the following lower bound on the greatest $d$ for which \eqref{Eq:Theorem1} holds, namely $d=d(n, k, k_0, inj_M, S_g)$ could be chosen to be equal to 
\begin{equation}\label{Eq:Theorem1DiameterEstimates}
C(n,k)^{-\frac{1}{n}}\left\{\frac{(n+2)K(n,1)^2}{nC_0(n,k_0)}\left[n(n-1)k_0-S_g\right]\right\}^{\frac{1}{4}},
\end{equation}
see equation \eqref{Eq:RadiusEstimates} for the exact meaning of the constants involved here. Notice that the dependence of $d$ on $inj_M>0$ appears just to ensure that such a $d$ exists and is positive.
\end{Res}
\begin{Rem} Strict inequality is necessary because as pointed out in $\cite{DruetPAMS}$ Theorem \ref{Res:Theorem1} is false if we have just $Ric_g\le(n-1)k_0$ and not $S_g<n(n-1)k_0$, as pointed out in $\cite{DruetPAMS}$. The comparison result is false also on $S^2\times S^2$, as noticed in $\cite{MJ}$, compare again $\cite{DruetPAMS}$ page $2352$.
\end{Rem}
In the next theorem we refine the results contained in Theorem \ref{Res:Theorem1}. The price to pay to have this stronger result is that the proof of Theorem \ref{Res:Theorem3} is much more involved.
\begin{Res}[Sharp small volumes in $C^2$-locally asymptotically strong bounded geometry smooth at infinity]\label{Res:Theorem3} Let $(M^n, g)$ be a complete Riemannian manifold, $n\geq 2$, with $C^2$-locally asymptotically strong bounded geometry smooth at infinity. Let us assume that there exists a real constant $k_0\in\R$ such that $S_g<n(n-1)k_0$. Then there exists a positive constant $\tilde{v}_0=\tilde{v}_0(n,k,k_0, inj_{M,g}, S_g)>0$ such that for every $\Omega\subseteq M$ finite perimeter set with $V_g(\Omega)\le\tilde{v}_0$ it holds 
\begin{equation}\label{Eq:ResTheorem3}
\P_g(\Omega)>\P_{g_{k_0}}(B),
\end{equation}
where $B\subseteq\mathbb{M}_{k_0}^n$ is a geodesic ball having $V_{g_{k_0}}(B)=V_g(\Omega)$. Moreover $\tilde{v}_0$ can be chosen as an arbitrary number 
\begin{equation}
0<\tilde{v}_0\le\left\{\frac{2(n+2)K(n,1)^2}{nC_0}\left[n(n-1)k_0-S_g\right]\right\}^{\frac{n}{4}}, 
\end{equation}
where $C_0=C_0(n,k_0)>0$. %and the dependence of $\tilde{v}_0$ on $inj_M>0$ appears just to ensure that such a $\tilde{v}_0$ exists and is positive.
\end{Res}
\begin{Res}[Small diameters in strong bounded geometry]\label{Res:Theorem1.1} Let $(M^n, g)$ be a complete Riemannian manifold with $n\geq 2$ and with strong bounded geometry. Let us assume that there exists a real constant $k_0\in\R$ such that $S_g<n(n-1)k_0$. Then there exists $d=d(n,k, k_0, inj_M, S_g)>0$ such that for every $\Omega\subseteq M^n$ finite perimeter set with diameter $\diam_g(\Omega)\leq d$ holds 
\begin{equation}\label{Eq:Theorem1.1}
\P_g(\Omega)>\P_{g_{k_0}}(B),
\end{equation}
where $B\subseteq\mathbb{M}_{k_0}^n$ is a geodesic ball having $V_{g_{k_0}}(B)=V_g(\Omega)$. Moreover we have the following lower bound on the greatest $d$ for which \eqref{Eq:Theorem1} holds, namely $d=d(n, k, k_0, inj_M, S_g)$ could be chosen to be equal to 
\begin{equation}\label{Eq:Theorem1DiameterEstimatesStrongBoundedGeometry}
C(n,k)^{-\frac{1}{n}}\left\{\frac{(n+2)K(n,1)^2}{nC_0(n,k_0)}\left[n(n-1)k_0-S_g\right]\right\}^{\frac{1}{4}},
\end{equation}
see equation \eqref{Eq:RadiusEstimates} for the exact meaning of the constants involved here. Notice that the dependence of $d$ on $inj_M>0$ appears just to ensure that such a $d$ exists and is positive.
\end{Res} 
A first consequence of Theorem \ref{Res:Theorem1.1} is the following result whose proof is much more simpler than that of our main Theorem \ref{Res:Theorem3.1}.
\begin{Res}[Small volumes \`a la B\'erard-Meyer]\label{Res:Theorem2} Let $(M^n, g)$ be a complete Riemannian manifold, $n\geq 2$, and with strong bounded geometry. Let us assume that there exists a real constant $k_0\in\R$ such that $S_{sup}<n(n-1)k_0$. Then for every $\varepsilon>0$ there exists a positive constant $\tilde{v}_0=\tilde{v}_0(M,\varepsilon)>0%n, k, v_0, k_0, \varepsilon)>0
$
such that for every $\Omega\subseteq M$ finite perimeter set with $V_g(\Omega)\leq\tilde{v}_0$ holds 
\begin{equation}\label{Eq:ResTheorem2}
\P_g(\Omega)>(1-\varepsilon)\P_{g_{k_0}}(B),
\end{equation}
where $B\subseteq\mathbb{M}_{k_0}^n$ is a geodesic ball having $V_{g_{k_0}}(B)=V_g(\Omega)$.
\end{Res}
\begin{Rem} This gives a refinement of the classical result of B\'erard-Meyer in \cite{BerardMeyer}. Of course Theorem \ref{Res:Theorem2} follows immediately from the stronger Theorem \ref{Res:Theorem3.1} below.
\end{Rem}
\begin{Res}[Sharp small volumes in strong bounded geometry]\label{Res:Theorem3.1} Let $(M^n, g)$ be a complete Riemannian manifold, $n\geq 2$, with strong bounded geometry. Let us assume that there exists a real constant $k_0\in\R$ such that $S_g<n(n-1)k_0$. Then there exists a positive constant $\tilde{v}_0=\tilde{v}_0(n,k,k_0, inj_{M,g}, S_g)>0$ such that for every $\Omega\subseteq M$ finite perimeter set with $V_g(\Omega)\le\tilde{v}_0$ it holds 
\begin{equation}\label{Eq:ResTheorem3.1}
\P_g(\Omega)\ge\P_{g_{k_0}}(B),
\end{equation}
where $B\subseteq\mathbb{M}_{k_0}^n$ is a geodesic ball having $V_{g_{k_0}}(B)=V_g(\Omega)$, and $\tilde{v}_0$ can be chosen as an arbitrary real number satisfying
\begin{equation}
0<\tilde{v}_0\le\left\{\frac{(n+2)K(n,1)^2}{nC_0}\left[n(n-1)k_0-S_g\right]\right\}^{\frac{n}{4}}, 
\end{equation}
where $C_0=C_0(n,k_0)>0$. Moreover up to take a smaller $\tilde v_1=\tilde v_1(n,k,k_0, inj_M, S_g, I_{M,g})\le\tilde v_0$ we have that for any finite perimeter set $\Omega$ with $V_g(\Omega)\le\tilde v_1$ we have 
\begin{equation}\label{Eq:ResTheorem3.1Strict}
\P_g(\Omega)>\P_{g_{k_0}}(B).
\end{equation}  
\end{Res}
\begin{Rem} Observe that the upper bound on $\tilde v_0$ of Theorem \ref{Res:Theorem3.1} is strictly less than the corresponding upper bound on $\tilde v_0$ of Theorem \ref{Res:Theorem3} as it could be noted looking at their respective expressions they differ by the multiplication of a factor $2^{\frac{n}{4}}$.
\end{Rem}
A particular case of the more general situation considered in Theorem \ref{Res:Theorem3} gives a positive answer to a special case of Aubin-Cartan-Hadamard's conjecture for small volumes as stated in the following corollaries.
\begin{CorRes}[Aubin's Conjecture in $C^2$-locally asymptotically strong bounded geometry smooth at infinity for small volumes]\label{CorRes:AubinCartanHadamard} Let $(M^n, g)$ be a Cartan-Hadamard manifold, $n\geq 2$ with $C^2$-locally asymptotically strong bounded geometry smooth at infinity, and $S_g<0$. Then there exists a positive constant $\tilde{v}_0=\tilde{v}_0(n,k,k_0, inj_{M,g}, S_g)>0$ such that for every $\Omega\subseteq M$ finite perimeter set with $V_g(\Omega)\le\tilde{v}_0$ it holds 
\begin{equation}\label{Eq:CorRes:AubinCartanHadamardC2}
\P_g(\Omega)>\P_{g_{k_0}}(B),
\end{equation}
where $B\subseteq\mathbb{M}_{k_0}^n$ is a geodesic ball having $V_{g_{k_0}}(B)=V_g(\Omega)$. Moreover $\tilde{v}_0$ can be chosen as an arbitrary number 
\begin{equation}
0<\tilde{v}_0\le\left\{\frac{2(n+2)K(n,1)^2}{nC_0}\left[-S_g\right]\right\}^{\frac{n}{4}}, 
\end{equation}
where $C_0=C_0(n,k_0)>0$.
\end{CorRes}
As a corollary of Theorem \ref{Res:Theorem3.1} we have the following statement.
\begin{CorRes}[Aubin's Conjecture in strong bounded geometry for small volumes]\label{CorRes:AubinCartanHadamard1} Let $(M^n, g)$ be a complete Riemannian manifold, $n\geq 2$, with strong bounded geometry, and $S_g<0$. Then there exists a positive constant $\tilde{v}_0=\tilde{v}_0(n,k,k_0, inj_{M,g}, S_g)>0$ such that for every $\Omega\subseteq M$ finite perimeter set with $V_g(\Omega)\le\tilde{v}_0$ it holds 
\begin{equation}\label{Eq:CorRes:AubinCartanHadamardC3}
\P_g(\Omega)\ge\P_{g_{k_0}}(B),
\end{equation}
where $B\subseteq\mathbb{M}_{k_0}^n$ is a geodesic ball having $V_{g_{k_0}}(B)=V_g(\Omega)$. Moreover $\tilde{v}_0$ can be chosen as an arbitrary number 
\begin{equation}
0<\tilde{v}_0\le\left\{\frac{2(n+2)K(n,1)^2}{nC_0}\left[-S_g\right]\right\}^{\frac{n}{4}}, 
\end{equation}
where $C_0=C_0(n,k_0)>0$. Moreover the all the conclusions of Theorem \ref{Res:Theorem3.1} holds.
\end{CorRes}
As a last consequence of Theorem \ref{Res:Theorem3.1} we get Corollary \ref{CorRes:PuiseuxSeries} which gives an asymptotic expansion of the isoperimetric profile in Puiseux's series up to the second non trivial order generalizing previous results of \cite{NarCalcVar}. Before to state the corollary we recall here the definition of the isoperimetric profile.
\begin{Def}\label{Def:IsPWeak} Let $(M,g)$ be an arbitrary Riemannian manifold. For every $v\in]0, V(M)[$ we define $I_{M,g}(v):=\inf\{\P_g(\Omega)\}$, where the infimum is taken over the family of finite perimeter subsets $\Omega\subseteq M$ having fixed volume $V(\Omega)=v$ that will be denoted in the sequel $\tilde{\tau}_v$. If there exists a finite perimeter set $\Omega$ satisfying $V(\Omega)=v$, $I_M(V(\Omega))=A(\partial\Omega)= \mathcal{P}(\Omega)$ such an $\Omega$ will be called an \textbf{isoperimetric region}, and we say that $I_M(v)$ is \textbf{achieved}. 
\end{Def}
\begin{CorRes}[Asymptotic expansion of the isoperimetric profile]\label{CorRes:PuiseuxSeries} If $(M,g)$ have strong bounded geometry, then 
\begin{equation*} 
I_{M,g}(v)=c_nv^{\frac{(n-1)}n}\left(1-\gamma_nS_gv^{\frac2n}\right)+O\left(v^{\frac4n}\right),
\end{equation*} 
when $v\to0^+$, where $S_g:=\sup_{x\in M}\left\{Sc_g(x)\right\}$ and $\gamma_n=\frac 1{2n(n+2)\omega_n^{\frac 2n}}$ is a positive dimensional constant. Here $\omega_n$ is the volume of a geodesic ball of radius $1$ in $\R^n$. 
\end{CorRes}
\begin{Rem} The preceding corollary roughly speaking means that up to the second nontrivial term the asymptotic expansion of $I_M$ coincides with $I_{\mathbb{M}^n_k, g_k}$, where $n(n-1)k=S_g$. 
\end{Rem}
\begin{Rem} Via the same smoothing results of \cite{Shi}, \cite{Kapovitch} one can prove the preceding asymptotic expansion just using the theory of pseudo-bubbles developped by the first author in \cite{NarCalcVar}. However the theory of pseudo bubbles does not give the sharp isoperimetric comparison of Theorem \ref{Res:Theorem3.1}.
\end{Rem}
\subsection{Aknowledgments} This paper is part of the second author's Ph.D. thesis written under the supervision of the first author at the Federal University of Rio de Janeiro, UFRJ-Universidade Federal do Rio de Janeiro. We wish to thank Pierre Pansu for attracting our attention on this problem and for many helpful insights and suggestions. We are also indebted to Olivier Druet for useful discussions and for teaching us a lot about the content of his nice and inspiring work \cite{DruetPAMS}. Finally we want to thank Frank Morgan for its useful remarks and comments.% and for useful discussions and suggestions. 
\section{Sobolev inequalities and the proof of Theorem \ref{Res:Theorem1.1} in strong bounded geometry}\label{Sec:Theorem1}
In this section we closely follow the proof of Theorem $1$ of \cite{DruetPAMS}. We just make the needed changes to get the proof of our Theorem \ref{Res:Theorem1}. First we set some notations and make the definitions that will be required in the sequel. By $\xi$ we denote the standard Euclidean metric of $\R^n$. For every $1\le p<n$, $K(n,p)>0$ is the best constant in the Sobolev inequalities on $(\R^n,\xi)$ defined as 
\begin{equation}
K(n,p)^{-p}:=\inf_{u\not\equiv 0,u\in C_c (\R^n)}\left\{\dfrac{\int_{\R^n}|\nabla  u|_\xi^p dv_\xi}{\left(\int_{\R^n} |u|^{p^*} dv_\xi\right)^{\frac{p}{p^*}}}\right\},
\end{equation}  
where  $p^*:=\frac{np}{n-p}$ is the critical Sobolev's exponent. The explicit value of $K(n,p)$ is computed in \cite{AubinJDG76}, \cite{Talenti} namely
\begin{eqnarray*}
K(n,1)&:=&\frac{1}{n}\left(\frac{n}{\omega_{n-1}}\right)^{1/n}=c_n^{-1},\\
K(n,p)&:=&\frac{1}{n}\left(\frac{n(p-1)}{n-p}\right)^{1-1/p}\left(\frac{\Gamma(n+1)}{\Gamma(n/p)\Gamma(n+1-n/p)\omega_{n-1}}\right)^{1/n}.
\end{eqnarray*}
However the only property of $K(n,p)$ that we will use is that  
$$\lim_{p\to 1^+}K(n,p)=K(n,1)=\frac{1}{n}\left(\frac{n}{\omega_{n-1}}\right)^{1/n}.$$
We will use frequently the $L^p$ and the $W^{1,p}$ norm on $M$ defined by
$$||u||_{p,g}:=\left(\int_M|u|^pdv_g\right)^\frac1p,$$
$$||u||_{1,p,g}:=||u||_{p,g}+||\nabla_g u||_{p,g},$$
For any function $u$ belonging respectively to $L^p(M)$ and $W^{1,p}(M)$.
%\textcolor{red}{Definir os espaços de Sobolev sobre variedades}.
When $1\le p<n$ we will need to work inside $W^{1,p}(\R^n)$ that will denote the standard Sobolev space defined as the completion of $C_c^{\infty}(\R^n)$ with respect to the norm 
\begin{equation}\label{GagliardoNirenberg}
||u||_{1,p,\xi}:=\left(\int_{\R^n}|\nabla_\xi u|_\xi^pdv_\xi\right)^\frac1p.
\end{equation} 
\begin{Def} Let $(M,g)$ be a Riemannian manifold of dimension $n$, $U\subseteq M$ an open subset, $\mathfrak{X}_c(U)$ the set of smooth vector fields with compact support on $U$. Given a function $u\in L^{1}(M)$, define the variation of $u$ by
      \begin{equation}
                 |Du|(M):=sup\left\{\int_{M}u div_g(X)dV_g: X\in\mathfrak{X}_c(M), ||X||_{\infty}\leq 1\right\},
      \end{equation}  
      where $||X||_{\infty}:=\sup\left\{|X_p|_{{g}_p}: p\in M\right\}$ and $|X_p|_{{g}_p}$ is the norm of the vector $X_p$ in the metric $g_p$ on $T_pM$.
      We say that a function $u\in L^{1}(M)$, has \textbf{bounded variation}, if $|Du|(M)<\infty$ and we define the set of all functions of bounded variations on $M$ by $BV(M):=\{u\in L^1(M):\:|Du|(M)<+\infty\}$. A function $u \in L^1_{loc}(M)$ has \textbf{locally bounded variation in $M$}, if
for each open set $U \subset\subset M$, 
$$|Du|(U):=\sup\left\{ \int_U u div_g(X)dV_g : X\in\mathfrak X_c(U), \|X\|_\infty\le1  \right\}<\infty,$$
and we define the set of all functions of locally bounded variations on $M$ by $BV_{loc}(M) := \{u \in L^1_{loc}(M) : |Du|(U) < +\infty, U\subset\subset M\}$.
      \end{Def}
\begin{Def}\label{Def:LocallyFinitePerimeterSets} Let $(M,g)$ be a Riemannian manifold of dimension $n$, $U\subseteq M$ be an open subset, $\mathfrak{X}_c(U)$ the set of smooth vector fields with compact support in $U$. Given $E\subset M$ measurable with respect to the Riemannian measure, the \textbf{perimeter of $E$ in $U$}, $ \mathcal{P}(E, U)\in [0,+\infty]$, is
      \begin{equation}
                 \mathcal{P}(E, U):=sup\left\{\int_{U}\chi_E div_g(X)dV_g: X\in\mathfrak{X}_c(U), ||X||_{\infty}\leq 1\right\},
      \end{equation}  
where $||X||_{\infty}:=\sup\left\{|X_p|_{{g}_p}: p\in M\right\}$ and $|X_p|_{{g}_p}$ is the norm of the vector $X_p$ in the metric $g_p$ on $T_pM$. If $\mathcal{P}(E, U)<+\infty$ for every open set $U\subset\subset M$, we call $E$ a \textbf{locally finite perimeter set}. Let us set $\mathcal{P}(E):=\mathcal{P}(E, M)$. Finally, if $\mathcal{P}(E)<+\infty$ we say that \textbf{$E$ is a set of finite perimeter}.    
\end{Def}
%\begin{Lemme}
% Let $n \ge 1$, $D$ be a bounded domain of $\R^n$, $h$ be a $C^\infty$ Riemannian metric in a neighborhood of $\overline D$, $1 < p < \infty$, $a\in L^\infty(\Omega)$. Suppose that $u \in W^{1,p}(D) \cap  C^1(\overline D)$ and $w\in C^2(D)\cap C^0(\overline D)$ satisfy
 %$$
 %\Delta_{p,h}u + a(x)|u|^{p-2}u\le0\quad \text{in $D$},
% $$
%and 
%$$
%\Delta_{p,h}w +a(x)|w|^{p-2}w\ge 0\quad \text{in $D$}.
%$$
%In addition we assume
%$$
%w(x)>0, \quad \forall x\in \overline D,\quad |\nabla_hw(x)|\ne 0, \quad \forall x\in D
%$$
%and 
%\begin{center}
% $u$ is $C^{1,1}$ in any open set where $\nabla_hu \ne 0$.
%\end{center}
%Then if $u\le w$ on $\partial D$, we have $u\le w$ in $\overline D$.
%\end{Lemme}
Before to prove Theorem \ref{Res:Theorem1} we prove Proposition \ref{Proposition:Sufficient} which is sufficient to prove Theorem \ref{Res:Theorem1}. We postpone the proof of this last fact to the end of this section.  In the proof of Proposition \ref{Proposition:Sufficient} we make frequent use of de Moser's iterative scheme, so we give an ad-hoc version of it in the following lemma which is suitable for our applications.
\begin{Lemme}[Ad-hoc De Giorgi-Nash-Moser]\label{Lemma:DeGiorgiNashMoser}
Let $(M^n,g)$ be a complete Riemannian manifold, $n \ge 2$, $1<p<n$, and $v \in W_g^{1,p}(M)$ with $0\le v \le 1$, satisfying\begin{equation}\label{eq:11remterm}
 %\int_{B_g(x,s)}\langle \nabla_{g} v, \nabla_{g} \phi \rangle dv_{g} \le \int_{B_g(x,s)}   \Lambda v^{p^*-1}\phi dv_{g}
\Delta_{p,g} v\le  \Lambda v^{p^*-1},
\end{equation}
where $\Lambda$ only depends on $n$. Then for any $x$ in $M$, 
%any $K>0$, any $p>1$, and any $q > p^*$, 
for any $\delta > 0$ %if $v$ is some nonnegative function of $W^{1,p}(M)$ satisfying \eqref{eq:11remterm} then
\begin{equation}
 \sup_{y\in M\cap B_g(x,\delta/2)}\{v(y)\}\le C\left( \int_{M\cap B_g(x,\delta)}v^{p^*}dv_{g} \right)^{\frac1{p^*}},
\end{equation}
 where $C> 0$ does not depend on $p$.
\end{Lemme}
\begin{Rem}
 Substituting the condition $0\le v\le 1$, by $\| v\|_{L^q(B(x,2\delta))}<K$ for a suitable value of $K$, and $q> p^*$, get the same result, the proof is based on the Moser iterative scheme applied to \eqref{eq:11remterm}. See for example Lemmas $3.1$ and $3.2$ of $\cite{AubinLi}$.
\end{Rem}

\begin{proof}
 Consider the inequality $\Delta_{p,g} v\le \Lambda v^{p^*-1}$ in $M$, and $v\le1$,
for some positive constant $\Lambda$ independent of $p$.
Consider a non-negative $\eta\in C_c^\infty(B_g(x,\delta))$ such that
for $0<r<s\leq\delta$ satisfies 
\begin{itemize}
\item[i.] $0\leq \eta \leq 1$,
\item[ii.] $\eta\equiv 1$, in $B_g(x,r)$,
\item[iii.] $\eta\equiv 0$, in $B_g(x,\delta)\setminus B_g(x,s)$,
\item[iv.] $|\nabla_{g} \eta|_g\leq \frac{C_0}{s-r}$, 
where $C_0$ depends only on the geometry  of $(M,g)$ or on the bounds of the geometry in case $M$ satisfy some condition of bounded geometry, for example instrong bounded geometry $C_0$ depends on $n,\Lambda_1,\Lambda_2, inj_M$.
\end{itemize}
Multiplying Equation \eqref{eq:11remterm} by $\eta^pv^{k+1}$, for $0<k\le p^*-p$, and integrating by parts over $B_g(x,s)$ leads to
\begin{equation}\label{eq:23900}
\int_{B_g(x,s)} |\nabla_{g} v|_g^{p-2}\langle\nabla_{g} v,\nabla_{g} (\eta^pv^{k+1}) \rangle dv_{g} \le \int_{B_g(x,s)}\Lambda v^{p^*+k}\eta^pdv_{g}.
\end{equation}
Let $w=v^{\frac{k+p}{p}}$, then $|\nabla_{g} w|_g^p=\left(	\frac{k+p}{p}\right)^p v^k |\nabla_{g} v|_g^p$.

We observe that
{\small
\begin{align*}
%\int_{B_g(x,s)}&\eta^pv^{k+1}\Delta_pvdv_{g} =
\int_{B_g(x,s)}&|\nabla_{g} v|_g^{p-2}\langle \nabla_{g} v,\nabla_{g} (\eta^pv^{k+1})\rangle dv_{g}\notag\\
=&(k+1)\int_{B_g(x,s)}|\nabla_{g} v|_g^{p} v^k\eta^p dv_{g}+\int_{B_g(x,s)}|\nabla_{g} v|_g^{p-2}v^{k+1}\langle \nabla_{g} v,\nabla_{g} (\eta^p)\rangle dv_{g}\\
=&(k+1)\left(\dfrac{k+p}{p}\right)^{-p}\int_{B_g(x,s)}|\eta\nabla_{g} w|_g^p  dv_{g}
+\int_{B_g(x,s)}|\nabla_{g} v|_g^{p-2}v^{k+1}\langle \nabla_{g} v,\nabla_{g} (\eta^p)\rangle dv_{g}\\
%+\int_{B_g(x,s)}|\nabla_{g} v|^{p-2}v^{k+1}\langle \nabla_{g} v,p\eta^{p-1}\nabla_{g} \eta\rangle dv_{g},\\
\ge&(k+1)\left(\dfrac{k+p}{p}\right)^{-p}\int_{B_g(x,s)}|\eta\nabla_{g} w|_g^pdv_{g}-\int_{B_g(x,s)}|\nabla_{g} v|_g^{p-1}v^{k+1}|\nabla_{g} \eta^p|_g dv_{g}\\
=&(k+1)\left(\dfrac{k+p}{p}\right)^{-p}\int_{B_g(x,s)}|\eta\nabla_{g} w|_g^pdv_{g}-p\int_{B_g(x,s)}|\nabla_{g} v|_g^{p-1}v^{k+1}\eta^{p-1}|\nabla_{g} \eta|_g dv_{g}\\
=&(k+1)\left(\dfrac{k+p}{p}\right)^{-p}\int_{B_g(x,s)}|\eta\nabla_{g} w|_g^pdv_{g}-p\int_{B_g(x,s)}\left(|\nabla_{g} \eta|_gv^{\frac{k+p}{p}}\right)\left(|\nabla_{g} v|_g^{p-1}\eta^{p-1}v^{\frac {k(p-1)}p}\right) dv_{g}\\
=&(k+1)\left(\dfrac{k+p}{p}\right)^{-p}\int_{B_g(x,s)}|\eta\nabla_{g} w|_g^pdv_{g}-p\left(\dfrac{k+p}{p}\right)^{1-p}\int_{B_g(x,s)}(w|\nabla_{g} \eta|_g)(\eta |\nabla_{g} w|_g) ^{p-1}dv_{g}
\end{align*}
}
where we applied the  Cauchy-Schwarz's inequality. Later we have to use in the second integral on the right the Young's inequality in the following form  
$$ab\le \dfrac{(\theta^{-1} a)^p}p+ \dfrac{(p-1)(\theta b)^{\frac p{p-1}}}{p},$$
with $\theta\in]0,+\infty[$.
Set $a=w|\nabla_{g} \eta|_{g}$, $b=(\eta |\nabla_{g} w|_{g})^{p-1}$, and choose 
$\theta>0$ such that
\begin{eqnarray}\label{036}
\left(\dfrac{k+p}{p}\right)^{1-p}(p-1)\theta^{\frac{p}{p-1}}=\frac{1}{2}(k+1)\left(\dfrac{k+p}{p}\right)^{-p},
\end{eqnarray}
we get 
\begin{tiny} 
\begin{eqnarray}\label{eq:23901}
\int_{B_g(x,s)}|\nabla_{g} v|^{p-2}\langle \nabla_{g} v,\nabla_{g} (\eta^pv^{k+1})\rangle dv_{g} & \ge & \dfrac{(k+1)}2\left(\dfrac{k+p}{p}\right)^{-p}\int_{B_g(x,s)}|\eta\nabla_{g} w|^pdv_{g}\\ \nonumber
& - & 2^{p-1}\left(\frac{p-1}{k+1}\right)^{p-1}\int_{B_g(x,s)} |w\nabla_{g} \eta|^p dv_{g}.
\end{eqnarray}
\end{tiny}
Combining \eqref{eq:23900} and \eqref{eq:23901} leads to 
\begin{align}\label{eq:23902}
\int_{B_g(x,s)} |\eta \nabla_{g} w|_g^p dv_{g}\le& C_1\int_{B_g(x,s)}\eta^pv^{k+p^*}dv_{g}
 +C_2\int_{B_g(x,s)}|w \nabla_{g} \eta|_g^pdv_{g}\notag\\
 =&C_1\int_{B_g(x,s)} (\eta w)^p v^{p^*-p}dv_{g}
 +C_2\int_{B_g(x,s)}|w \nabla_{g} \eta|_g^pdv_{g},
\end{align}
where $C_1(p)=\frac2{k+1}\left(\frac{k+p}{p}\right)^{p}\Lambda$, and $C_2(p)=2^{p}(p-1)^{(p-1)}$.

Independently we have the following computations
\begin{align}\label{eq:mink456}
\int_{B_g(x,s)}|\nabla_{g} (\eta w)|_g^pdv_{g}\leq 2^{p-1}\left(\int_{B_g(x,s)} |\eta\nabla_{g} w|_g^pdv_{g}+\int_{B_g(x,s)}|w\nabla_{g} \eta|_g^pdv_{g}\right),
\end{align}
%pois $|X+Y|^p\leq2^{p-1}(|X|^p+|Y|^p)$, usando novamente o item (iv.) na desigualdade precedente temos que
and by the Sobolev embedding we get
\begin{eqnarray*}
\left(\int_{B_g(x,s)}\right.|\left. \eta w|^{p^*} dv_{g}\right)^{\frac{p}{p^*}} 
& \le &  C(n,p) \int_{B_g(x,s)}|\nabla_{g} (\eta w)|_g^pdv_{g},\\
&\stackrel{\eqref{eq:mink456}}{\le} & C_3\int_{B_g(x,s)}|\eta\nabla_{g} w|_g^pdv_{g}+C_3\int_{B_g(x,s)}|w\nabla_{g} \eta|_g^pdv_{g}\\
&\stackrel{\eqref{eq:23902}}{\le} & C_4\int_{B_g(x,s)}(\eta w)^p v^{p^*-p}dv_{g}+C_5\int_{B_g(x,s)}|w\nabla_{g} \eta|_g^pdv_{g},
\end{eqnarray*}
where $C_3=2^{p-1}C(n,p)$, $C_4=C_1C_3$, and $C_5=C_2C_3+C_3$.

But since $v\le 1$, $0\le \eta\le 1$ over $B_g(x,s)$, and $|\nabla_{g} \eta|_g\le \frac{C_0}{s-r}$ we have that

  {\small
\begin{align*}
\left(\int_{B_g(x,s)} |\eta w|^{p^*} dv_{g}\right)^{\frac{p}{p^*}}
&\le \left[ C_4 +  C_5\left( \frac{C_0}{s-r}\right)^p\right]\int_{B_g(x,s)} | w|^{p} dv_{g}.
\end{align*}
}
On the other hand we have $\left(\frac{k+p}{p}\right)^{p}\le (k+1)^p$ for $k>0$ and $1<p<n$, then
$C_1\le 2 (k+1)^{n-1}\Lambda$, $C_2\le 2^n(n-1)^{n-1}$, $C_3\le2^{n-1}C(n)$, $C_4\le 2^n(k+1)^{n-1}\Lambda C(n)$, $C_5\le 2^{n-1}(2^{n}(n-1)^{n-1}(k+1)^{n-1}+1)C(n)$.
Then 
{\scriptsize 
$$C_4 +  C_5\left( \frac{C_0}{s-r}\right)^p\le 2^{n-1}C(n)\left( 2 (k+1)^{n-1}\Lambda+(2^n(n-1)^{n-1}+1)\left( \frac{C_0}{s-r}\right)^p \right).
$$
}
Thus setting 
$$B_0=2^{n-1}C(n)\left( 2 (k+1)^{n-1}\Lambda+(2^n(n-1)^{n-1}+1)\left( \frac{C_0}{s-r}\right)^p \right),$$ we get
\begin{equation}\label{eq:moserite1}
\left(\int_{M\cap B_g(x,r)} |v|^{\frac{p^*(k+p)}{p}} dv_{g}\right)^{\frac{p}{p^*(k+p)}}
\le  B_0^{\frac{1}{k+p}} \left(\int_{M\cap B_g(x,s)} | v|^{k+p} dv_{g}\right)^{\frac{1}{k+p}}.
\end{equation}

Now we want to use the Moser's iterative scheme. Let us call
\begin{eqnarray*}
F(t,\rho)=\left(\int_{M\cap B_g(x,\rho)}v^tdv_{g}\right)^{1/t},
\end{eqnarray*}
by the inequality \eqref{eq:moserite1} we get that 
\begin{eqnarray}\label{eq:031}
F\left((k+p)\frac{p^*}{p},r\right)\leq B_0^{\frac{1}{k+p}}F((k+p),s).
\end{eqnarray}
Choose $k_0$ such that $(k_0+p)=p^*$, $s_0=\delta$ and define for every  $i\geq 1$
\begin{eqnarray*}
(k_i+p)=\frac{p^*}{p}(k_{i-1}+p)=\left( \frac{p^*}{p}\right)^i(k_0+p), s_i=\frac{\delta}{2}+\frac{\delta}{2^{i+1}}. 
\end{eqnarray*}
Make $k=k_i$, $s=s_i$ and $r=s_{i+1}$. Note that
$s_i-s_{i+1}=\frac\delta{2^{i+2}}$, furthermore  we get
$k_i\rightarrow +\infty$ when $i\rightarrow +\infty$, because 
\begin{equation*}
k_{i+1}-k_i%&=&\left[\left(\frac{p^*}{pq}\right)^{i+1}(k_0+p)-\left(\frac{p^*}{pq}\right)^{i}(k_0+p)\right],
%\\
=(k_0+p)\left(\frac{p^*}{p}\right)^{i}\left[\frac{p^*}{p}-1\right]>0.
\end{equation*}
Now we apply this to \eqref{eq:031}, and we obtain
\begin{eqnarray*}
F\left(\frac{k_i+p}{p}p^*,s_{i+1}\right)=F((k_{i+1}+p),s_{i+1})\leq B_i^{\frac{1}{k_i+p}}F((k_i+p),s_i) .
\end{eqnarray*}
Then making the iteration yields
\begin{eqnarray*}
F((k_{i+1}+p),s_{i+1})\le \prod_{j=0}^{i}B_j^{\frac{1}{k_j+p}}F((k_0+p),s_0) =\prod_{j=0}^{i}B_j^{\frac{1}{k_j+p}}F(p^*,\delta).
\end{eqnarray*}
Taking $i\to \infty$ the expression above becomes
\begin{eqnarray}\label{047}
\|v\|_{L_g^\infty(B_g(x,\delta/2))}\leq \prod_{i=0}^{+\infty}B_i^{\frac{1}{k_i+p}}\left(\int_{M\cap B_g(x,\delta)}v^{p^*}dv_{g}\right)^{\frac{1}{p*}}.
\end{eqnarray}
It remains to prove the convergence of $\prod_{i=0}^{+\infty}B_i^{\frac{1}{k_i+p}}$ to a constant independent of $p$.

Since for $p$ sufficiently close to $1$ we can get that
 $\frac{p^*}{p}=\frac{n}{n-p}\le 2$, we have 
$$(k_i+1)< (k_i+p)=\left( \frac{p^*}{p}\right)^i(k_0+p)\le 2^i(k_0+n),$$
and making the choice $\tilde{\delta}=\min\{1,\delta \}$, we get
\begin{eqnarray*}
B_i&=&2^{n-1}C(n)\left( 2 (k_i+1)^{n-1}\Lambda+(2^n(n-1)^{n-1}+1)\left( \frac{C_0}{s_{i+1}-s_i}\right)^p \right)
\\
&\le&2^{n-1}C(n)\left( 2 (k_0+n)^n2^{in}\Lambda+(2^n(n-1)^{n-1}+1)\left( \frac{2^{(i+2)^p}C_0^p}{\delta^p}\right)\right)\\
&\le&2^{in} 2^{n-1}C(n)\left( 2 (k_0+n)^n\Lambda+(2^n(n-1)^{n-1}+1) 2^{2n} C_0^n \tilde \delta^n\right)\\
&=&2^{in}\tilde C(n).
 \end{eqnarray*}

As it is easy to see from the definition of $k_i$ we have $$\frac{1}{2^i(k_0+n)}\le \frac{1}{k_i+p}\le \frac{1}{k_i+1}.$$ Let us define
$\alpha_i:=\frac{1}{2^i(k_0+n)}$, if $B_i<1$ and $\alpha_i:=\frac{1}{k_i+1}$, if $B_i\ge1$, in any case we have
$$
B_i^\frac{1}{k_i+p}\le B_i^{\alpha_i}.
$$
Then passing to the infinite products 
\begin{eqnarray*}
\prod_{i=0}^{+\infty} B_i^{\frac{1}{k_i+p}}&\leq & \prod_{i=0}^{+\infty}B_i^{\alpha_i}\leq \left( \prod_{i=0}^{+\infty}\tilde{C}^{\alpha_i}\right)\left( \prod_{i=0}^{+\infty} \left(2^{in}\right)^{\alpha_i}\right)
\\ 
&= &\left(\tilde{C}^{\sum_{i=0}^{+\infty}\alpha_i}\right)\left(2^{n\sum_{i=0}^{+\infty}i\alpha_i}\right).
\end{eqnarray*}
Notice that  
$$\sum_{i=0}^{+\infty}\alpha_i \quad \text{and }\quad \sum_{i=0}^\infty i\alpha_i,$$
are convergent series. Then for values of $p$ close to $1$ we have
\begin{eqnarray}\label{eq:04700}
\|v\|_{L_g^\infty(B_g(x,\delta/2))}\leq C\left(\int_{M\cap B_g(x,\delta)}v^{p^*}dv_{g}\right)^{\frac{1}{p*}},
\end{eqnarray}
where $C$ does not depend on $p$.
\end{proof}
\begin{Prop}[Proposition \cite{DruetPAMS} page $2353$]\label{Prop:Druet}
 Let $(M^n, g)$ be a complete Riemannian manifold 
of dimension $n \ge 2$.
 Let  $x_0 \in M$, and define $\alpha_\varepsilon:= \dfrac{n}{n-2} Sc_g(x_0)+\epsilon$. Then 
for any $\varepsilon>0$, there exists $r_\varepsilon >0$ such that for any $u$ in
$C_c^\infty( B_g ( x_0 , r_\epsilon ) )$ we have that
\begin{equation}\label{Eq:PropDruet}
\|u\|^{2}_{\frac{n}{n-1},g}\le K(n,1)^{2}\left(\|\nabla u\|^2_{1,g} +\alpha_\varepsilon\|u\|^2_{1,g}\right).
\end{equation}
\end{Prop}
The preceding proposition justifies the following definition.
\begin{Def}
Let $(M,g)$ be a Riemannian manifold. Let us define $r^*_\varepsilon(M,g, x)\in[0,+\infty]$ as the supremum of all $r>0$ such that  \eqref{Eq:RespropSufficient} is satisfied. Of course $r^*_\varepsilon(M,g, x)=0$, if there is no such positive $r_\varepsilon$. We call $r^*_\varepsilon(M,g,x)$ the \textbf{Druet's radius of $(M,g)$ at $x$}. Let us define $r^*_\varepsilon(M,g)\in[0,+\infty]$ as the infimum of $r^*_\varepsilon(M,g, x)$ taken over all $x\in M$. We call $r^*_\varepsilon(M,g)$ the \textbf{Druet's radius of $(M,g)$}.
\end{Def}
By Proposition \ref{Prop:Druet}, if $(M^n,g)$ is complete then for any $x\in M$ we have $r^*_\varepsilon(M,g, x)>0$. By  Proposition \ref{Proposition:Sufficient} if $(M^n,g)$ has $C^2$-locally asymptotic strong bounded geometry smooth at infinity, then we have $r^*_\varepsilon(M,g)>0$.
We want to study now a little of stability properties of Druet's radius with respect to the convergence of manifolds.
\begin{Lemme}\label{Lemma:stability} Suppose to have a sequence of pointed complete smooth Riemannian manifolds $(M,g_i,p_i)\to (M_{\infty}, g_{\infty}, p_{\infty})$ in $C^0$ topology with $(M_{\infty}, g_{\infty}, p_{\infty})$ smooth and $Sc_{g_i}(p_i)\to Sc_{g_{\infty}}(p_\infty)$. Then  \begin{equation}
\liminf_{i\to+\infty}r^*_\varepsilon(M,g_i, p_i)\le r^*_\varepsilon(M_{\infty},g_{\infty}, p_{\infty}).
\end{equation} 
\end{Lemme}

\begin{ResProp}\label{Proposition:Sufficient}  Let $(M, g)$ be a complete Riemannian manifold of dimension $n\ge2$ with $C^2$-locally asymptotic strong bounded geometry smooth at infinity. For any $\varepsilon>0$ there exists $r_{\varepsilon}=r_{\varepsilon}(M,g)>0$ such that for any point $x_0\in M$, any function $u\in C^{\infty}_c(B_g(x_0,r_{\varepsilon}))$, we have
\begin{equation}\label{Eq:RespropSufficient}
||u||^2_{\frac n{n-1},g}\le K(n,1)^2\left(||\nabla_g u||_{1,g}^2+\alpha_\varepsilon||u||_{1,g}^2\right),
\end{equation} 
where $\alpha_{\varepsilon}=\frac n{n+2}S_g(x_0)+\varepsilon$.
\end{ResProp}
\begin{Rem} We notice that the constant $r_{\varepsilon}=r_{\varepsilon}(M)>0$ is obtained by contradiction and that the proof does not give an explicit effective lower bound on it. 
\end{Rem} 
We state here a Corollary of Lemma \ref{Lemma:stability} that could be used to give a slightly different proof of Theorem \ref{Res:Theorem1.1}.
\begin{Cor} Let $(M, g)$ be a complete Riemannian manifold of dimension $n\ge2$ with $C^2$-locally asymptotic strong bounded geometry smooth at infinity, $p_i\to\infty$ and $(M, g, p_i)\to(M_{\infty}, g_{\infty}, p_{\infty})$. Then for every $\varepsilon>0$ we have $r^*_\varepsilon(M_{\infty}, g_{\infty})\ge r^*_\varepsilon(M,g)>0$.
\end{Cor}
\begin{proof}[Proof of Proposition \ref{Proposition:Sufficient}] For any $x_0\in M$, for any $r>0$, any $p>1$ and any $\varepsilon>0$, set
\begin{tiny}
\begin{equation*}
\lambda_{p,r,g}(x_0):=\inf_{\substack{u\in C^\infty_c (B_g (x_0,r))\\ u\not\equiv 0}}   \dfrac{ \left(\int_{B_g(x_0,r)} |\nabla_g u|^pdv_g\right)^{2/p} +\alpha_\varepsilon\left(\int_{B_g(x_0,r)} |u|^pdv_g\right)^{2/p}}{\left(\int_{B_g(x_0,r)} |u|^{p^*}dv_g\right)^{2/{p^*}}},
\end{equation*}
\end{tiny}
where $B_g(x_0,r)\subseteq M$ is the geodesic ball $(M,g)$ centered at $x_0\in M$ and of radius $r>0$.
We will argue the theorem by contradiction. With this aim in mind suppose that there exists $\varepsilon_0>0$ such that for every $r>0$ there exists a point $x_{0,r}$ depending on $r$ such that it holds
\begin{equation*}
\lambda_{1,r,g}(x_{0,r})<K(n,1)^{-2}.
\end{equation*} 
As it is easy to check from the very definition of $\lambda_{p,r,g}(x_{0,r})$, we have that $\limsup_{p\to1^+}\lambda_{p,r,g}(x_{0,r})\le\lambda_{1,r,g}(x_{0,r})$, which implies that for any $r>0$, there exists $p_r(x_{0,r})>1$ such that
\begin{tiny}
\begin{eqnarray}\label{eq:2.2}
	\lambda_{p_r,r,g}(x_{0,r})<K(n,1)^{-2}\left(\dfrac{n-p_r(x_{0,r})}{p_r(x_{0,r})(n-1)}\right),
	\:\lambda_{p_r(x_{0,r}),r,g}<K(n,p_r(x_{0,r}))^{-2}.
\end{eqnarray}
\end{tiny}
We may assume that $r\searrow 0$ and we may choose $p_r(x_{0,r})$ decreasing when $r$ is decreasing.
Then inverting this sequence we get a sequence $p > 1$ going to $1^+$ a sequence $r_p > 0$ going to $0^+$ as $p$ goes
to $1^+$, and a sequence of points $x_{0,p}:=x_{0,r_p}\in M$ which verify \eqref{eq:2.2}. Notice here that in general the sequence of points $x_{0,p}$ could go to infinity when $p\to1^+$. This is the main difficulty that we encounter in adapting the original proof of Theorem $1$ of \cite{DruetPAMS} in case of noncompact ambient manifolds. Set $\alpha_p:=\dfrac n{n+2} S_g(x_{0,p})+\varepsilon_0$. Now up to a subsequence we can assume that 
\begin{equation}\label{Eq:Prop1alphap}
\lim_{p\to1^+}\alpha_p=\frac n{n+2}l_1+\varepsilon_0,
\end{equation} 
for some $l_1\in [S_{inf,g}, S_g]$, where $S_{inf,g}:=\inf\{Sc_g(x):x\in M\}$ and $S_g:=\sup\{Sc_g(x):x\in M\}$. It is worth to note that $S_{inf,g}$ and $S_g$ are finite real numbers, because $S_g$ is bounded from below by $n(n-1)k$ and from above by $n(n-1)k_0$. The second equation in \eqref{eq:2.2} can be written like $\lambda_{p,r_p}(x_{0,p})<K(n,p)^{-2}$, and by Theorems $1.1$, $1.3$, $1.5$ of \cite{DruetPRSEdinburgh}, we have the existence of a minimizer $u_p$ which satisfies  
\begin{align}\label{eq:druet2.3}
C_p\Delta_{p,g}u_p+\alpha_p\|u_p\|_{p,g}^{2-p}u_p^{p-1}=\lambda_pu_p^{p^*-1}, \quad\text{in } B_g(x_{0,p},r_p),\\ \notag
u_p \in C^{1,\eta} (B_g (x_{0,p},r_p)), \quad\text{for some } \eta> 0,\\ \notag
u_p>0, \quad\text{in } B_g (x_{0,p},r_p),\quad u_p=0, \quad\text{in } \partial B_g (x_{0,p},r_p),
\end{align}

\begin{equation}\label{eq:druet2.4}
 \int_{B_g(x_{0,p},r_p)}u_p^{p^*}dv_g=1,
\end{equation}
\begin{equation}\label{eq:druet2.5}
 \lambda_{p}<K(n,p)^{-2}, \quad \lambda_p<K(n,1)^{-2}\left(\dfrac{n-p}{p(n-1)}\right)^{2},
\end{equation}
\begin{equation}\label{eq:druet2.6}
 C_p:=\left(\int_{B_g(x_{0,p},r_p)}|\nabla_g u_p|_g^pdv_g\right)^{\frac{2-p}p},
\end{equation}
where $\lambda_p:=\lambda_{p,r_p}(x_{0,r_p})$,
%Observe that $\eta$ is independent of $p$ as explained after Equation $(12)$ p. $361$ of  \cite{AubinLi} and the references therein mentioned \textcolor{red}{(.............To write some more words about the proof of this fact.........)}. 
$\Delta_{p,g}$ is the $p$-Laplacian with respect to $g$, defined by $\Delta_{p,g}u:= -div_g (|\nabla_g u|^{p-2}\nabla_g u)$, with $\nabla_gu$ being the gradient of $u$ with respect to the metric $g$. The strategy that will adopt to go head in this proof is concerned with the study of the sequence $(u_p)$ as $p\to1^+$. With this aim in mind, let $x_p$ be a point in $B_g(x_{0,p},r_p)$ where $u_p$ achieves its maximum ($x_p$ tends to infinity, iff $x_{0,p}$ tends to infinity) and we define
$$
u_p(x_p)=\mu_p^{1-\frac np}.
$$
Observing that $u_p(x_p)^{p^*}=\mu_p^{-n}$ we get
\begin{equation}\label{Eq:PropSuffOrderofMagnitude}
1=\int_{B_g(x_{0,p},r_p)}u_p^{p^*}dv_g\le V_g(B_g(x_{0,p},r_p))\mu_p^{-n}\le C_0(n,k)r_p^n\mu_p^{-n},
\end{equation}
where the last inequality is due to Bishop-Gromov. Since $r_p$ goes to $0$, $\mu_p$ goes to $0$ as $p$ goes to $1^+$, moreover $\mu_p=O(r_p)$ and the constant $C_0=C_0(n,k)$ is uniform with respect to $p$, i.e., is uniform with respect to the location of $x_{0,p}$ inside $M$. Analogously, applying H\"older's inequalities, with $q=\frac n{n-p}>1$ yields 
\begin{equation}\label{eq:druet2.7}
\lim_{p\to1^+}\int_{B_g(x_{0,p},r_p)} u_p^pdv_g\le\lim_{p\to1^+}\left\{\int_{B_g(x_{0,p},r_p)}u_p^{p^*}dv_g\right\}^{\frac{1}{q}}V_g(B_g(x_{0,p},r_p))^{\frac1{q'}}=0,
\end{equation}
here $q'$ denotes the conjugate exponent of $q$, i.e., $\frac1q+\frac1{q'}=1$.
\paragraph{Step 1} In this first step we want to show the validity of the two following equations
\begin{equation}\label{eq:druet2.8}
\lim_{p\to1^+}\lambda_p=K(n,1)^{-2},
\end{equation}
and
\begin{equation}\label{eq:druet2.9}
\lim_{p\to1^+}\int_{B_g(x_{0,p},r_p)}|\nabla_g u_p|^p_gdv_g=K(n,1)^{-1}.
\end{equation}
By Theorem $7.1$ of \cite{Hebey}, it follows that for all $\varepsilon > 0$ there exists $B_\varepsilon=B_\varepsilon(n, k, inj_{M,g})>0$ such that for any $p>1$,
\begin{scriptsize}
\begin{eqnarray*}
\left(\int_{B_g(x_{0,p},r_p)}u_p^{p^*}dv_g\right)^{2\frac{n-1}{n}} & \le &
(K(n,1)+\varepsilon)^2\left(\int_{B_g(x_{0,p},r_p)} \left| \nabla_g \left(u_p^{\frac{p(n-1)}{n-p}}\right) \right|_gdv_g \right)^2\\ 
& + & B_\varepsilon\left(\int_{B_g(x_{0,p}0,r_p)}u_p^{\frac{p(n-1)}{n-p}}dv_g\right)^2,
\end{eqnarray*}
\end{scriptsize}
which gives with \eqref{eq:druet2.3}, \eqref{eq:druet2.4} and Hölder's inequalities
\begin{eqnarray*}
1 & \le & \left(K(n,1)+\varepsilon\right)^2\left(\dfrac{p(n-1)}{n-p}\right)^2
(\lambda_p-\alpha_p\|u_p\|_p^2)+B_\varepsilon \|u_p\|_p^2\\
 & \le & \left(K(n,1)+\varepsilon\right)^2\left(\dfrac{p(n-1)}{n-p}\right)^2
(\lambda_p-\alpha_0\|u_p\|_p^2)+B_\varepsilon \|u_p\|_p^2,
\end{eqnarray*}
where $\alpha_0:=n(n-1)k+\varepsilon_0\le\alpha_p$.
This combined with \eqref{eq:druet2.7} give
$$
1\le\left(1+\varepsilon K(n,1)^{-1}\right)^2\liminf_{p\to1^+}\left(\lambda_p K(n,1)^2\right).
$$
Since this inequality is valid for every $\varepsilon>0$, letting $\varepsilon\to 0$ we obtain that
$\liminf_{p\to1^+}\lambda_p\ge K(n,1)^{-2}$.
Using the fact that 
$\lambda_p<K(n,p)^{-2}$, and $\lambda_p<K(n,1)^{-2}\left(\dfrac{n-p}{p(n-1)}\right)^2$, we conclude that
$$
K(n,1)^{-2}\le \liminf_{p\to1^+}\lambda_p\le\liminf_{p\to1^+} K(n,p)^{-2},
$$
and thus \eqref{eq:druet2.8} is proved.
\begin{Rem}
Until this point we have just used the assumption of mild bounded geometry, i.e., Ricci bounded below and positive injectivity radius.
\end{Rem}
Now, it easily seen that \eqref{eq:druet2.9} is an obvious consequence of  \eqref{eq:druet2.3}, \eqref{eq:druet2.4}, \eqref{eq:druet2.7}, and \eqref{eq:druet2.8}.
\begin{Rem}
To prove rigorously \eqref{eq:druet2.9} we have used that the scalar curvature is bounded from both sides that Ricci is bounded below and the injectivity radius is positive.
\end{Rem}
To prove the equality \eqref{eq:druet2.9}, we observe that as a consequence of \eqref{eq:druet2.3} we have
\begin{equation}\label{eq:druet2.3equivalent}
\|\nabla_g u_p\|^2_{p,g}+\alpha_p\|u_p\|_{p,g}^2=\lambda_p\|u_p\|^{p^*}_{p^*,g},
\end{equation}
and that under the assumptions of Proposition \ref{Proposition:Sufficient}, $n(n-1)k\le Sc_g\le n(n-1)k_0$. This yields immediately that $(\alpha_p)$ is a bounded sequence, i.e., $\alpha_p=O(1)$. Taking $p\to1^+$, in \eqref{eq:druet2.3equivalent} and \eqref{eq:druet2.4} we obtain that
$$
\lim_{p\to1^+} \|\nabla_g u_p\|^2_{p,g}+\alpha_p\|u_p\|_{p,g}^2= \lim_{p\to1^+} \lambda_p\|u_p\|^{p^*}_{p^*,g}=\lim_{p\to1^+} \lambda_p,
$$
then using the above result, we can conclude that
$$
\lim_{p\to1^+} \|\nabla_g u_p\|_{p,g}= K(n,1)^{-1}.
$$
\paragraph{Step 2} Let 
$\Omega_p:=\mu_p^{-1}\exp_{x_p}^{-1}(B_g(x_{0,p},r_p))\subset T_{x_p}M\cong\R^n$, the metric $g_p(x):=\exp_{x_p}^*g(\mu_px)$ for $x\in \Omega_p$, and the function given by $v_p(x)=\mu_p^{\frac np-1}u_p(\exp_{x_p}(\mu_px))$ for $x\in\Omega_p$, $v_p(x)=0$ in $x\in \R^n\setminus\Omega_p$.
It is worth to mention that the definition of $v_p$ for $p$ close to $1$ is well posed and does not give any problem, since we suppose that the injectivity radius of $M$ is strictly positive and $r_p\to 0$ as $p\to1^+$. Then making the substitution in \eqref{eq:druet2.3} we obtain that $v_p$ satisfies
\begin{equation}\label{eq:druet2.10}
C_p\Delta_{p,g_p}v_p+\alpha_p\mu_p^2 \|v_p\|_{p,g_p}^{2-p}v_p^{p-1}=\lambda_pv_p^{p^*-1}, \text{in $\Omega_p$}, 
\end{equation}
with $v_p=0$ in $\partial\Omega_p$, and
\begin{equation}\label{eq:druet2.11}
\int_{\Omega_p}v_p^{p^*}dv_g=1.
\end{equation}
Thus $v_p$ satisfy also
\begin{equation}\label{eq:druet2.3equivalent1} 
\|\nabla_{g_p} v_p\|^2_{p,g_p}+\alpha_p\|v_p\|_{p,g_p}^2=\lambda_p\|v_p\|^{p^*}_{p^*,g_p}.
\end{equation}
Unfortunately the sequence $(v_p)$ is not bounded in $W_\xi^{1,1}(\R^n)$ so we need to pass to another auxiliary sequence $(\tilde v_p)$ related in some way to the preceding one and that is bounded in $W_\xi^{1,1}(\R^n)$. We do this because we are interested in a limit function $v_0$ that realizes the minimum of the problem at infinity and that we expect to be the characteristic function of a ball.
To realize this strategy we look for powers of the function $v_p$. As we will see later a suitable choice is the following
\begin{equation}\label{eq:druet2.12}
\tilde v_p(x)=v_p(x)^{\frac{p(n-1)}{n-p}}.
\end{equation}
It is useful to recall here that for every $x\in M$ the exponential map $exp_x$ is a bi-Lipschitz map of an open geodesic ball centered at $x$ having radius $inj_x$ over a ball of $\R^n$ having the same radius, with Lipschitz constant $L_x$ that in general depend on $x$, however by the Rauch's comparison Theorem we know that $L_x$ can be bounded by a constant that depends just on $\Lambda_1,\Lambda_2, inj_M$, which in turn permit to conclude that under our assumption of strong bounded geometry the constants $L_x$ are uniformly bounded with respect to $x$ by a positive constant that depends only on the bounds on the geometry, namely $\Lambda_1, \Lambda_2, inj_M$. Hence using the Cartan's expansion of the metric $g_p$ close to $x_p$ we can show the existence of a positive constant $C=C(\Lambda_1,\Lambda_2, inj_M)>0$, such that for any $x\in \Omega_p$,
\begin{equation}\label{Eq:Comparison}
(1-C\mu^2_p|x|^2)dv_{g_p}\leq dv_{\xi} \leq(1+C\mu^2_p |x|^2)dv_{g_p}.
\end{equation}
From this we conclude that there exists another constant again denoted by $C=C(\Lambda_1,\Lambda_2, inj_M)>1$ such that
\begin{eqnarray}\label{Eq:Prop1BiLip}
dv_{g_p}\ge\left(1-\frac1C\mu_p^2\right)dv_\xi,%\textcolor{red}{Olivier improvement????????}\\
\end{eqnarray}
\begin{equation}\label{Eq:Prop1BiLip1}
|\nabla v_p|^p_{g_p}dv_{g_p}\le(1+C\mu^2_p)|\nabla_\xi v_p|^p_\xi d v_\xi,%\textcolor{red}{Olivier improvement????????}
\end{equation}
where $\xi$ is the Euclidean metric. Equations \eqref{Eq:Comparison}, \eqref{Eq:Prop1BiLip},  \eqref{Eq:Prop1BiLip1} with \eqref{eq:druet2.9}, \eqref{eq:druet2.11} and Hölder's inequalities leads to
\begin{equation}\label{eq:druet2.13}
\lim_{p\to1^+}\dfrac{\int_{\R^n}|\nabla_\xi \tilde v_p|_\xi dv_\xi}{\left(\int_{\R^n}\tilde v_p^{\frac n{n-1}} dv_\xi\right)^{\frac{n-1}{n}}}=K(n,1)^{-1}.%\textcolor{red}{OlivierCheckit?????}
\end{equation}
To show this observe that by \eqref{Eq:Comparison}, \eqref{Eq:Prop1BiLip},  \eqref{Eq:Prop1BiLip1} $\left(\int_{\R^n}\tilde v_p^{\frac n{n-1}} dv_\xi\right)^{\frac{n-1}{n}}\sim\int_{\Omega_p}v_p^{p^*}dv_g=1$, when $p\to 1^+$. To see what happens to the numerator of \eqref{eq:druet2.13} just look at \eqref{Eq:Prop1Holder} below
\begin{eqnarray}\nonumber
\int_{\R^n}|\nabla_\xi\tilde v_p|_\xi dv_\xi & = & \int_{\R^n}\frac{p(n-1)}{n-p}v_p^{\frac{n(p-1)}{n-p}}|\nabla_\xi v_p|_\xi dv_\xi\\ \nonumber
 & \le & \frac{p(n-1)}{n-p}\left\{\int_{\R^n} v_p^{p^*}dv_\xi\right\}^{\frac p{p-1}}\left\{\int_{\R^n}|\nabla_\xi v_p|^p_\xi dv_\xi\right\}^{\frac 1p}\\ \nonumber
 & \le & \frac{p(n-1)}{n-p}(1+C\mu_p^2)||\nabla _{g_p}v_p||_{p,g_p}\\ \label{Eq:Prop1Holder}
 & = & \frac{p(n-1)}{n-p}(1+C\mu_p^2)||\nabla_g u_p||_{p,g}\to K(n,1)^{-1}.
\end{eqnarray} 
The last equality is a consequence of \eqref{eq:druet2.9} and the following calculation
\begin{align*}
 \int_{\Omega_p}|\nabla_{g_p} v_p(x)|_{g_p} ^{r}dv_{g_p}
 &=\int_{\Omega_p} |\nabla_{g_p} (\mu_p^{\frac{n-p}{p}} u_p(\exp_{x_p}(\mu_px)))|_{g_p}^rdv_{g_p}(x)\\
 &=\mu_p^{(\frac{n-p}{p})r}\int_{\mu_p^{-1}\exp_{x_p}^{-1}(B_g(x_{0,p},r_p))} |\nabla_{g_p}  u_p(\exp_{x_p}(\mu_px))|_{g_p}^rdv_{g_p}(x)\\
% &=\mu_p^{(\frac{n-p}{p})r}\int_{B_0(R)} |\nabla  u_p(\exp_{x_p}(\mu_px))|^rdv_{\tilde g_p}(x)\\
% &=\mu_p^{\frac{n-p}{p}r}\int_{B_0(R)} u_p^{r}(\exp_{x_p}(\mu_px))dv_{\tilde g}(x)\\
 &=\mu_p^{\frac{n-p}{p}r}\mu_p^{-n}\mu_p^{r}\int_{\exp_{x_p}^{-1}(B_g(x_{0,p},r_p))} |\nabla_{g_p}  u_p(\exp_{x_p}(x))|_{g_p}^r  dv_{g_p}(\mu_p^{-1}x)\\
 &=\mu_p^{\frac{n(r-p)}{p}}\int_{\exp_{x_p}^{-1}(B_g(x_{0,p},r_p))} |\nabla_{g_p}  u_p(\exp_{x_p}(x))|_{g_p}^r  dv_{ g_p}(\mu_p^{-1}x)\\
 &=\mu_p^{\frac{n(r-p)}{p}}\int_{B_g(x_{0,p},r_p)}  |\nabla_g u_p(x)|_g^{r}dv_{g},
\end{align*}
from which follows
$$
\|\nabla_{g_p} v_p\|_{r,g_p}^r=\mu_p^{\frac{n(r-p)}{p}}\|\nabla_gu_p\|_{r,g}^r.
$$
Remember here that $r_p\to0$ as $p\to1^+$. Notice that by \eqref{GagliardoNirenberg}, $(\tilde v_p)$ is bounded in $W_\xi^{1,1}(\R^n)$. Thus there exists $v_0 \in BV_{loc}(\R^n)$ such that
$$
\lim_{p\to1^+}\tilde v_p = v_0,\quad\text{strictly  in } BV_{loc}(\R^n),
$$
this means that $\tilde v_p\to v_0$ in $L^1_{loc}(\R^n)$ e $||\nabla\tilde v_p||_{1,\xi}(K)\to|Dv_0|(K)$, $\forall K\subset\subset\R^n$. For a proof of this fact see Thm. $3.23$ of \cite{AmbrosioFuscoPallara}. If we apply the concentration-compactness principle of P.L. Lions (\cite{LionsCCIAIHP}, \cite{LionsCCIRMIA}, see also \cite{Struwe} for an exposition in book form) to $|v_p|^{p^*} dv_\xi$ , four situations may occur: compactness, concentration, dichotomy or vanishing. Dichotomy is classically forbidden by \eqref{eq:druet2.13}. To be convinced of this fact the reader could mimic the proof of Theorem $4.9$ of \cite{Struwe}. Concentration without compactness cannot happen since $\sup_{\Omega_p}  v_p = v_p(0) = 1$. As for vanishing, since $v_p$ is bounded in $L^\infty$, by applying Moser's iterative scheme (see for instance Theorem $1$ \cite{Serrin}) to \eqref{eq:druet2.10}, one gets the existence of some $C=C(n,\alpha_pC_p^{-1}\mu_p^2||v_p||_p^{2-p},C_p^{-1}\lambda_p, ||g||_{0,r})> 0$ such that for any $p > 1$,
\begin{equation}\label{Eq:Prop1Moser}
1= \sup_{\Omega_p\cap B_{g_p}(0,1/2)} v_p \le C\left(\int_{\Omega_p\cap B_{g_p}(0,1)}  v_p^{p^*} dv_{g_p}\right)^{\frac1{p^*}},
\end{equation}
where $||g||_{0,r}$ is the norm defined at page $308$ of \cite{Petersen} (see Definition \ref{Def:BoundedGeometryPetersen}). Since a careful analysis of the proof of Theorem $1$ of \cite{Serrin} combined with \eqref{eq:druet2.6}, \eqref{eq:druet2.7}, \eqref{eq:druet2.8}, \eqref{eq:druet2.9}, (which imply, by a change of variables in the integrals, that $\alpha_p\mu_p^2||v_p||_p^{2-p}=\alpha_p\mu_p^p||u_p||^{2-p}\to 0$, thanks to the fact that $\alpha_p\to l_1+\varepsilon_0\in\R$, hence $\alpha_p$ is uniformly bounded), and the $C^0$ convergence of the metric tensor due to Theorems $72$ and $76$ of \cite{Petersen}, when $p\to1^+$, shows that $C$ is uniformly bounded with respect to $p$. 
Thus vanishing cannot happen. Another way to see that our problem have no vanishing is to apply directly Lemma \ref{Lemma:DeGiorgiNashMoser} with $g=g_p$ and $v=v_p$, this is justified because in equation \eqref{eq:druet2.10}
%\begin{equation*}
%\Delta_{p,g_p}v+C_p^{-1}\alpha_p\mu_p^2 \|v\|_p^{2-p}v^{p-1}=C_p^{-1}\lambda_pv^{p^*-1},
%\end{equation*}
%$$
%C_p=\left(\int_{\Omega_p} v^{p*}dv_{g}\right)^{\frac {2-p}{p}},
%$$
we know that $\alpha_p\mu_p^2 \|v_p\|_{p,g_p}^{2-p}\to 0$, and $C_p^{-1}\lambda_p\to K(n,1)^{-1}$. Then for $p$ close to $1$, we can consider that $v$ satisfies the following inequality
$$
\Delta_{p,g}v_p\le \Lambda v_p^{p^*-1},
$$
where $\Lambda$ depends only on $n$. %Then it is possible to apply the De Giorgi-Nash-Moser iterative scheme to $v$. For this particular case we present this form of the following lemma.
Compactness implies that $|v_p|^{p^*} dv_\xi\to|v_0|^{\frac{n}{n-1}}dv_\xi $ that is $\|\tilde v_p\|_{\frac n{n-1}}\to \|v_0\|_{\frac n{n-1}}$. To see this we observe that by the compactness case of the concentration-compactness principle we have that for all $\varepsilon>0$ there exist $R_\varepsilon>0$ and $p_\varepsilon>0$ such that 
$$
1-\varepsilon\le \int_{B_\xi(0,R_\varepsilon)}v_p^{p^*}dv_\xi\le 1+\varepsilon, \quad p\le p_\varepsilon,
$$
passing to the limit when $p\to1^+$ yields $\int_{\R^n}v_0^{\frac n{n-1}}=1$, since 
%and how $\lim_{p\to 1} \tilde v_p=v_0$ and $\tilde v_p^{1^*}=v_p^{p^*}$ we get that $v_0\ge 0$ a.e.  
$$
1-\varepsilon\le \int_{B_\xi(0,R_\varepsilon)}v_0^{\frac n{n-1}}dv_\xi
=\lim_{p\to1^+}\int_{B_\xi(0,R_\varepsilon)} v_p^{p^*}dv_\xi\le 1+\varepsilon,
$$
and
\begin{equation}\label{Eq:PropLnorm}
\int_{\R^n}v_0^{\frac n{n-1}}dv_\xi
=\lim_{\varepsilon\to0^+}\int_{B_\xi(0,R_\varepsilon)} v_0^{\frac n{n-1}}dv_\xi.
\end{equation}
It is clear that $\|\tilde v_p\|_{\frac n{n-1}}$ is bounded by all $p>1$, on the other hand, as is well known $L^{\frac n{n-1}}(\R^n)$ is a reflexive Banach space thus $\tilde v_p\rightharpoonup v_0$ weakly in $L^{\frac n{n-1}}(\R^n)$. A classical result ensures that weak convergence and convergence of norms as in \eqref{Eq:PropLnorm} gives $\tilde v_p\to v_0$ strongly in $L^{\frac n{n-1}}(\R^n)$.

Since we have that $\int_{\R^n}|\nabla \tilde v_p|dv_\xi\to \int_{\R^n}|\nabla v_0|dv_\xi=K(n,1)^{-1}$. 
Then $v_0$ is a minimizer for the $W^{1,1}$ Euclidean Sobolev inequality which verifies
$$
\int_{\R^n}v_0^{\frac n{n-1}}dv =1.
$$ 

Thus there exists $y_0\in\R^n$, $\lambda_0>0$ and $R_0 >0$ such that

\begin{equation}\label{eq:druet2.15}
 v_0 = \lambda_0{\bold 1}_{B_\xi(y_0,R_0)},
\end{equation}
where ${\bold 1}_{B_\xi(y_0,R_0)}$ denotes the characteristic function of the Euclidean ball $B_\xi(y_0,R_0)$, and moreover, since $v_p \le  1$ in $\Omega_p$ we obtain by pointwise convergence a.e. $dv_\xi$ that $0\le \lambda_0\le1$. On the other hand $v_p\le1$ and the strong convergence in $L^{\frac n{n-1}}(\R^n)$  give that for all $q\ge\frac n{n-1}$, $\tilde v_p\to v_0$ strongly in $L^q(\R^n)$. Therefore   
$$\lambda_0^{q}V_\xi ( B_\xi(y_0 , R_0 ))=\lim_{p\to1^+}\int_{\R^n}\tilde v_p^{\frac n{n-1}}dv_\xi=1,\:\forall q\ge\frac n{n-1}.$$
Taking the limit when $q\to+\infty$ we deduce that $\lambda_0$ cannot be strictly less than $1$, thus we get $\lambda_0 =1$. So we have 
\begin{equation}\label{eq:druet2.16}
 V_\xi ( B_\xi(y_0 , R_0 )) = \frac{\omega_{n-1}}n R_0^n = 1 .
 \end{equation}
Up to changing $x_p$ into $\exp_{x_p} (\mu_py_0)$ in the definition of $v_p$, $\Omega_p$  and $g_p$, we may assume that $y_0 = 0$. 
In particular we have
\begin{equation}\label{eq:druet2.17}
\lim_{p\to1}\tilde v_p ={\bold 1}_{B_\xi(0,R_0)},\quad\text{strongly  in } L^{\frac{n}{n-1}}(\R^n),
\end{equation}
and that for any $\varphi\in C_c^\infty (\R^n)$,   
\begin{equation}\label{eq:druet2.18}
 \lim_{p\to1}\int_{\R^n}|\nabla \tilde v_p|_\xi^p\varphi dv_\xi=\int_{\partial B_\xi(0,R_0)}\varphi d\sigma_{\xi},
 \end{equation}
 where $d\sigma_\xi$ is the $(n-1)$-dimensional Riemannian measure of $\partial B_\xi(0,R_0)$ induced by the metric $\xi$ of $\R^n$.
Consider the extremal functions $V_p\in W^{1,p}(\R^n)$ for $K(n,p)^{-p}$ defined below
\begin{equation}\label{eq:druet2.19}
V_p(x)=\left(1+\left(\frac{|x|}{R_0}\right)^{\frac p{p-1}}\right)^{1-\frac np},\quad x\in\R^n,
\end{equation}

a simple application of the concentration-compactness principle, using \eqref{eq:druet2.18}, gives

\begin{equation}\label{eq:druet2.20}
\lim_{p\to1^+} \int_{\R^n} |\nabla_\xi( \tilde v_p -V_p ) |_\xi dv_\xi =0.
\end{equation}
Applying again the Moser's iterative scheme Lemma \ref{Lemma:DeGiorgiNashMoser} to \eqref{eq:druet2.10} with the help of \eqref{eq:druet2.17}, we also get that for any $R > R_0$,
\begin{equation}\label{eq:druet2.21}
\lim_{p\to1^+} \sup_{\Omega_p\setminus B_{g_p}(0,R)}  v_p =0.
\end{equation}
The application of Moser's iterative scheme is possible in strong bounded geometry because of the same arguments leading to \eqref{Eq:Prop1Moser}.
\paragraph{Step 3}In this step we want to obtain from the $L^{\frac n{n-1}}$-estimate \eqref{eq:druet2.17} the pointwise estimates \eqref{eq:druet2.24}, \eqref{eq:druet2.25} which gives estimates on the decay rate to zero of $v_p(z)$ when $|z|\to+\infty$. For more details one can see for instance \cite{DruetGeoDed2002} and \cite{DruetMathAnn1999}. With this aim in mind let us define
\begin{eqnarray}\label{46}
w_p(x)=|x|^{\frac{n}{p}-1}v_p(x),
\end{eqnarray}
and let  $z_p\in\Omega_p$ be a point where $w_p$ attains its maximum, i.e.,
\begin{eqnarray}\label{48}
w_p(z_p)=\| w_p \|_\infty.
\end{eqnarray}
Suppose by contradiction that
\begin{eqnarray}\label{47}
\lim_{p\to  1} \| w_p \|_\infty =\lim_{p\to  1}w_p(z_p)=+\infty.
\end{eqnarray}
Now we set
$$
\nu_p^{1-\frac np}=v_p(z_p),
$$
this implies by \eqref{47} that
\begin{eqnarray}\label{eq:druet2.22}
\lim_{p\to  1}|z_p|v_p(z_p)^{\frac p{n-p}}=
\lim_{p\to  1}\frac{|z_p|}{\nu_p}=+\infty.
\end{eqnarray}
Using the fact that $v_p\le 1$ in $\Omega_p$ and \eqref{46}, we conclude that
\begin{eqnarray}\label{eq:druet2.23}
\lim_{p\to  1}|z_p|=+\infty.
\end{eqnarray}

Consider $\exp_{g_p, z_p}$ the exponential map associated to $g_p$ at $z_p$, let
$\tilde \Omega_p=\nu_p^{-1}\exp_{g_p,z_p}^{-1}(\Omega_p)$, the metric $\tilde g_p(x)=\exp_{g_p,z_p}^*g_p(\nu_px)$ for $x\in \tilde\Omega_p$, and the function given by
$$
\phi_p(x)=\nu_p^{\frac np-1}v_p(\exp_{z_p}(\nu_px))\text{ for $x\in\tilde \Omega_p$, $\phi_p(x)=0$ in $x\in \tilde\Omega_p^c$}.
$$
Then for $x\in B(0,1)$, by \eqref{eq:druet2.22}, and \eqref{eq:druet2.23}, we can prove that $\phi_p$ is uniformly bounded in $B(0,1)$, and verifies \eqref{eq:139Druet}. %with the metric $\tilde g_p$

In fact, for $x\in B(0,1)$, and by the definition of the $\exp_{g_p,z_p}$ map we have
\begin{eqnarray*}
	\nu_p &\geq &d_{g_p}(z_p,\exp_{z_p}(\nu_px)),
	%u_p(z_p)^{\frac{p}{p-n}}&\geq &d_{g_p}(z_p,\exp_{z_p}(u_p(z_p)^{\frac{p}{p-n}}x)),
\end{eqnarray*}
using the triangular inequality get
\begin{eqnarray}\label{024}
|\exp_{z_p}\left(\nu_px\right)|&\geq & |z_p|-d_{g_p}(z_p,\exp_{z_p}(\nu_px)), \notag\\ \notag
&\geq & |z_p|-\nu_p\\\notag
&=& |z_p|-\left(w_p(z_p)|z_p|^{1-\frac{n}{p}}\right)^{\frac{p}{p-n}}\\
&=&\left(1-w_p(z_p)^{\frac{p}{p-n}}\right)|z_p|.
\end{eqnarray}

Since $w_p(z_p)\to  +\infty$ when $p\to  1$ and $\frac{p}{p-n}<0$ for values of $p$ very close to $1$, and since \eqref{48} and \eqref{47}, are valid for $x\in B(0,1)$ we obtain
\begin{eqnarray}\label{54}
|\exp_{z_p}\left(\nu_px\right)|\geq\frac{1}{2}|z_p|.
\end{eqnarray}
Rewriting in terms of $w_p$ we get
%De consequ\^encia para qualquer $x\in B(0,1)$, rescrevendo $(\ref{51})$ pela defini\c c\~ao temos
\begin{eqnarray*}
\phi_p(x)&=&\nu_p^{\frac np-1}v_p( \exp_{z_p}(\nu_px))\\
&=&\nu_p^{\frac np-1} w_p(\exp_{z_p}(\nu_px))|\exp_{z_p}\left(\nu_px\right)|^{1-\frac{n}{p}}.
\end{eqnarray*}
Since  $1-\frac{n}{p}<0$ for values of  $p$ close to 1, we obtain 
\begin{eqnarray*}
\phi_p(x)&\leq & \nu_p^{\frac np-1}w_p(\exp_{z_p}(\nu_px))\left(\frac{1}{2}|z_p|\right)^{1-\frac{n}{p}},
\end{eqnarray*}
and since $z_p$ is the maximum of $w_p$,  we have
$w_p(\exp_{z_p}(\nu_px))\leq w_p(z_p)$, thus
\begin{eqnarray*}
\phi_p(x)&\leq & 2^{\frac{n}{p}-1}|z_p|^{1-\frac{n}{p}} \nu_p^{\frac np-1}w_p(z_p),
\end{eqnarray*}
and we know by definition that   $v_p(z_p)^{-1}w_p(z_p)=|z_p|^{\frac{n}{p}-1}$, we are lead to
\begin{eqnarray*}
\phi_p(x)\leq 2^{\frac{n}{p}-1},
\end{eqnarray*}
that is, $\| \phi_p \|_{L^\infty(B_\xi(0,1))} \leq 2^{\frac{n}{p}-1}$.
Making the needed substitution in \eqref{eq:druet2.10} a straightforward computation gives that $\phi_p$ satisfies
\begin{align}\label{eq:139Druet}
C_p\Delta_{p,\tilde g_p}\phi_p+\alpha_p\mu_p^2\nu_p^{2} \|\phi_p\|_p^{2-p}\phi_p^{p-1}=\lambda_p\phi_p^{p^*-1}, \text{ in $\tilde\Omega_p$,}
\end{align}
and $\phi_p=0$ in $\partial\tilde\Omega_p$.
Since $\phi_p$ is uniformly bounded we can apply Moser's iterative scheme Lemma \ref{Lemma:DeGiorgiNashMoser} to the equation \eqref{eq:139Druet} and to get the existence of some $C > 0$ independent of $p$ such that
\begin{eqnarray}\label{eq:refgal118}
1=\phi_p(0)\le\sup_{\tilde \Omega_p\cap B(0,1/2)} \phi_p \le C\left(\int_{\tilde \Omega_p\cap B(0,1)}\phi_p^{p^*} dv_{\tilde g_p}\right)^{\frac 1{p^*}}.
\end{eqnarray}
For a subsequent use remember that 
\begin{eqnarray}\label{eq:refgal119}
\int_{\tilde \Omega_p\cap B_\xi(0,1)}\phi_p^{p^*} dv_{\tilde g_p}=
\int_{\Omega_p\cap B_{g_p}(z_p,\nu_p)}v_p^{p^*} dv_{g_p}.
\end{eqnarray}
Again an application of the Moser's iterative scheme Lemma \ref{Lemma:DeGiorgiNashMoser} is legitimate by the same arguments leading to \eqref{Eq:Prop1Moser}.
Therefore by \eqref{eq:druet2.17} we get immediately that 
\begin{eqnarray}\label{eq:refgal119b}
\liminf_{p\to 1} \int_{B_{g_p}(z_p,\nu_p)\cap \Omega_p}v_p^{p^*}dv_{g_p}>0.
\end{eqnarray}
By \eqref{eq:refgal118} and \eqref{eq:refgal119} given $R>0$, we get $B_{g_p}(0,R)\cap B_{g_p}(z_p,\nu_p)=\emptyset$ because $|z_p|\to\infty $ when $p\to 1$. Furthermore  
\begin{equation*}
1=\int_{\Omega_p}v^{p^*}dv_{g_p}=\int_{\Omega_p\cap B_{g_p(0,R)}}v^{p^*}dv_{g_p}+\int_{\Omega_p\setminus B_{g_p(0,R)}}v^{p^*}dv_{g_p},
\end{equation*}
on the other hand by \eqref{Eq:Comparison}, \eqref{Eq:Prop1BiLip} we have 
\begin{eqnarray*}
\left(1-C\mu_p^2\right)\int_{B_\xi\left(0,\frac R{\mu_p}\right)}v_p^{p^*}dv_\xi& \le & \int_{\Omega_p\cap B_{g_p(0,R)}}v^{p^*}dv_{g_p}\\
& \le & \left(1+C\mu_p^2\right)\int_{B_\xi\left(0,\frac R{\mu_p}\right)}v^{p^*}dv_\xi,
\end{eqnarray*} 
taking the limit when $p\to 1$ in the last two equations, using \eqref{Eq:PropLnorm}, \eqref{eq:druet2.17}, \eqref{eq:refgal119b} and remembering that $\mu_p\to0$ when $p\to1^+$ we get easily 
\begin{equation*}
0<\liminf_{p\to 1^+}\int_{B_{g_p}(z_p,\nu_p)\cap \Omega_p}v_p^{p^*}dv_{g_p}\le\liminf_{p\to 1^+}\int_{\Omega_p\setminus B_{g_p(0,R)}}v^{p^*}dv_{g_p}=0, 
\end{equation*}
which is the desired contradiction. Since this contradition comes from taking for granted \eqref{47}, we are lead to negate \eqref{47} and to have the existence of some $C > 0$ such that for any $p > 1$, and for all $x\in\Omega_p$
\begin{equation}\label{eq:druet2.24}
w_p(x)=|x|^{\frac{n}{p}-1}v_p(x) \le  C.
\end{equation}
In the same way, using \eqref{eq:druet2.24}, one proves thanks to \eqref{eq:druet2.21} that for any $R > R_0$,
\begin{equation}\label{eq:druet2.25} 
\lim _{p\to1}\sup_{\Omega_p\setminus B_{g_p}(0,R)} |x|^{\frac np-1} v_p(x)  = 0.
\end{equation}
To prove \eqref{eq:druet2.25} we argue by contradiction so we suppose that there exist $y_p\in \Omega_p$ and $\delta>0$ such that
\begin{eqnarray*}
	\lim_{p\to  1}|y_p|=+\infty, \qquad\text{and }\: w_p(y_p)\ge\delta.
\end{eqnarray*}
%Choose $\sigma>0$ small such that $B(0,\sigma)\subset \Omega_p$, and 
Define $v_p(y_p)=\nu_p^{1-\frac np}$, and $\tilde{\Omega}_p=\nu_p^{-1}\exp_{y_p}^{-1}(\Omega_p)$. Observe that $w_p(y_p)=|y_p|^{\frac{n}{p}-1}\nu_{p}^{1-\frac np}\ge\delta$. For $x\in \tilde{\Omega}_p$, let $\phi_p(x)=\nu_p^{\frac np-1}v_p(\exp_{y_p}(\nu_px))$ and $\phi_p(x)=0$ in $x\in \tilde\Omega_p^c$, and $\tilde g_p(x)=\exp_{y_p}^*g_p(\nu_px)$.

Now for any $x\in B_\xi\left(0,\frac12\delta^{\frac p{n-p}}\right)$, by the same arguments that above, we get that  $|\exp_{y_p}(\nu_p x)|\ge \frac12|y_p|$. Then using \eqref{eq:druet2.24}, we get that
\begin{eqnarray*}
	\phi_p(x)&=&\nu_p^{\frac np-1}v_p(\exp_{y_p}(\nu_px))=
	\nu_p^{\frac np-1}w_p(\exp_{y_p}(\nu_px))|\exp_{y_p}(\nu_px)|^{1-\frac np}\\
	&\le& C 2^{\frac np-1} |y_p|^{1-\frac np}\nu_p^{\frac np-1}\le C 2^{\frac np-1}\delta^{-1}.
\end{eqnarray*}
That is $\|\phi_p\|_{L^\infty\left(B_\xi\left(0,\frac12\delta^{\frac p{n-p}}\right)\right)}\le C 2^{\frac np-1}\delta^{-1}$, and by Moser's iterative scheme Lemma \ref{Lemma:DeGiorgiNashMoser} we get that 
$$
1=\phi_p(0)\le\sup_{\tilde \Omega_p\cap B_\xi\left(0,\frac14\delta^{\frac p{n-p}}\right)}\phi_p\le
C\left(\int_{\tilde \Omega_p\cap B_\xi\left(0,\frac12\delta^{\frac p{n-p}}\right)} 
\phi_p^{p^*}dv_{\tilde g_p}\right)^{\frac1{p^*}}.
$$
On the other hand, since 
$$\int_{\tilde \Omega_p\cap B_\xi\left(0,\frac12\delta^{\frac p{n-p}}\right)} 
\phi_p^{p^*}dv_{\tilde g_p}=\int_{\Omega_p\cap B_{g_p}\left(y_p,\frac12\delta^{\frac p{n-p}}\nu_p\right)} 
v_p^{p^*}dv_{g_p},
$$ 
using the same arguments as above we get that for $R>0$ for $p$ close to $1$
$$
B_{g_p}\left(y_p,\frac14\delta^{\frac p{n-p}}\nu_p\right)\cap B_{g_p}(0,R)=\emptyset.
$$
But for $R>R_0$, by \eqref{eq:druet2.21} we have $\lim_{p\to1}\sup_{\Omega_p\setminus B_{g_p}(0,R)}v_p=0$ , and $$\Omega_p\cap B_{g_p}\left(y_p,\frac12\delta^{\frac p{n-p}}\nu_p\right)\subset \Omega_p\setminus B_{g_p}(0,R),$$ thus
$$
1\le\lim_{p\to1}\sup_{\tilde \Omega_p\cap B_\xi\left(0,\frac12\delta^{\frac p{n-p}}\right)}\phi_p=0,
$$
which is a contradiction. 
%%%% END STEP 3 %%%%%%

%%%% NEW STEP 4 %%%%%%
\paragraph{Step 4} 
Unfortunately the pointwise estimates that we obtained in \eqref{eq:druet2.24} is not enough to prove our crucial \eqref{eq:druet2.34}. For this reasons we need to improve it. This is the goal to achieve in this step $4$, which culminate in the proof of \eqref{eq:druet2.26} below. Consider the following operator
$$
L_pu:=
C_p\Delta_{p,g_p}u+\alpha_p\mu_p^2 \|v_p\|_{p,g_p}^{2-p}u^{p-1}-\lambda_pv_p^{p^*-p}u^{p-1}.
$$

Choose $0 < \nu < n - 1$ and put
$$G_p(x) = \theta_p|x|^{-\frac{n-p-\nu}{p-1}},$$
where $\theta_p$ is some positive constant to be fixed later. 

%Easy computations lead to
We will use the following relation for the $p$-Laplacian for radial functions that could be found in Lemma $1.2$ of \cite{RodneyThesis} for an arbitrary Riemannian metric $h$
$$-\Delta_{p,h}u = -\Delta_{p,\xi}u + O(r)|\partial_ru|^{p-2}\partial_ru.$$
and we obtain
\begin{small}
$$
|x|^p\frac{L_pG_p(x)}{G_p(x)^{p-1}}
\ge C_p\nu\left(\dfrac{n-p-\nu}{p-1}\right)^{p-1}-C\mu_p^2|x|^2
+\alpha_p\mu_p^2\|v \|^{2-p}_p|x|^p - \lambda_p |x|^pv_p^{p^*-p},
$$
\end{small}
in $\Omega_p\setminus\{0\}$. Here $C$ denotes some constant independent of $p$.
Thanks to \eqref{eq:druet2.7}, \eqref{eq:druet2.8}, \eqref{eq:druet2.9}, \eqref{eq:druet2.25} and the fact that $r_p\to0$ as $p\to1$, one gets that for any $R>R_0$, $L_pG_p(x) \ge  0$ in $\Omega_p\setminus B_{g_p} (0, R)$
for $p$ close enough to $1$.
On the other hand,
$$
L_pv_p=0 \quad\text{in }\Omega_p.
$$
At last, it is not too hard to check with \eqref{eq:druet2.21} that
$$
v_p \le G_p \quad\text{on } \partial B_{g_p} (0, R),
$$
if we take 
$\theta_p = R^{\frac{n-p-\nu}{p-1}}$. 
Now we may apply the maximum principle as stated for
instance in Lemma $3.4$ of \cite{AubinLi} to get, 
$$
v_p(y)
\le\left(\frac R{|y|}\right)^{\frac{n-p-\nu}{p-1}}\quad\text{in } \Omega_p\setminus B_{g_p} (0, R),
$$
for $p$ close enough to $1$.
This inequality obviously holds on $B_{g_p} (0,R)$ and so we have finally obtained that for any 
$n-1>\nu > 0$ and any $R > R_0$, there exists $C(R,\nu) > 0$ such that for any $p > 1$ and any $y \in   \Omega_p$,
\begin{equation}\label{eq:druet2.26}
 \left(\frac {|y|}R\right)^{\frac{n-p-\nu}{p-1}}v_p(y)\le C(R,\nu).
\end{equation}
%%%% END STEP 4%%%%%%%

%%%%% NEW STEP 5 %%%%%
\paragraph{Step 5} We give in this Step the final arguments to conclude the proof of our Proposition \ref{Proposition:Sufficient}. 
%Since $\mu_p\to 0$ when $p\to 1$, and how  $g_p(x)=\exp^*g(\mu_px)$ we get that
%$g_p\to \xi$ in $C^2_{loc}(\R^n)$, where $\xi$ is the Euclidean metric.
We apply the $W_\xi^{1,1}(\R^n)$ Euclidean  Sobolev inequality to $\tilde v_p$:
\begin{equation}\label{eq:druet2.27}
\left(\int_{\Omega_p} \tilde v_p^{\frac n{n-1}} dv_\xi\right)^{\frac {n-1}n}\le K(n,1)\int_{\Omega_p}|\nabla_\xi  \tilde v_p|_\xi dv_\xi. 
\end{equation}
Recalling the Cartan expansion of $g_p$ around $0$, we have
\begin{equation}\label{eq:druet2.28}
dv_\xi=\left(1+\dfrac16\mu_p^2Ric_g(x_p)_{ij}x^ix^j+o(\mu_p^2|x|^2)\right)dv_{ g_p},
\end{equation}
where $Ric_g$ denotes the Ricci curvature of $g$ in the $\exp_{x_p}$ -map.
This last formula is true because $Ric_{g_p}(0)=\mu_p^2Ric_{g}(\exp_{x_p}(0))=\mu_p^2Ric_{g}(x_p)$.
Thus, by \eqref{eq:druet2.11} we obtain
$$
\int_{\Omega_p}\tilde v_p^{\frac n{n-1}}dv_\xi=1+\dfrac16\mu_p^2Ric_g(x_p)_{ij}\int_{\Omega_p}x^ix^jv_p^{p^*}dv_{g_p}+o\left(\mu_p^2   \int_{\Omega_p}|x|^2v_p^{p^*}  dv_{g_p} \right).
$$
%And since in $\Omega_p\setminus B_{g_p}(0,R)$ by \eqref{eq:druet2.2} and \eqref{eq:druet2.26}
%$$
%|x|^2v_p^{p^*}=\left( |x|^{\frac np-1}v_p\right)^{\frac{p^2}{n-p}}|x|^{2-p}v_p^p
%\le \left( \sup_{\Omega_p\setminus B_{g_p}(0,R)}  |x|^{\frac np-1}v_p\right)^{\frac{p^2}{n-p}}|x|^{2-p}v_p^p=o(1)
%$$
%we get that
%$$
%Ric_g(x_p)_{ij}\int_{\Omega_p\setminus B_{g_p}(0,R)}x^ix^jv_p^{p^*}dv_{g_p}=o(1)
%$$
To estimate the last term on the right hand side of the preceding equality we need to prove \eqref{Eq:Prop1Step5Mean} and \eqref{Eq:Prop1Step5Estimates} below 
\begin{align}\nonumber
\int_{B_\xi(0,R_0)}x^ix^jdv_{\xi}&=\frac{\delta^{ij}}{n}\int_{B_\xi(0,R_0)}|x|^2dv_{\xi}
=\frac{\delta^{ij}}{n}\int_0^{R_0}\int_{\partial B_\xi(0,r)}r^2d\sigma_{\xi}dr\\ \nonumber
&=\frac{\delta^{ij}}{n}\int_0^{R_0} r^{n+1}dr \int_{\partial B_\xi(0,1)}d\sigma_{\xi}\\ \label{Eq:Prop1Step5Mean}
 &=\dfrac{\delta^{ij}}{n(n+2)} \omega_{n-1 }R_0^{n+2}.
\end{align}
%Then
%\begin{align}
% Area (S_p(r))=&r^{n-1}\int_{S^{n-1}}\left(1-\frac {r^2}6Ric(p)_{i,j}x^ix^j+o(r^2|x|^2)\right)dx\\
% =&r^{n-1}\omega_{n-1}\left(1-\frac {r^2}{6n}S_g(p)_{i,j}+o(r^2|x|^2)\right)\\
%\end{align} 
Let $\beta_p=\frac{n-p-\nu}{p-1}$ and $R>\max\{1, R_0\}$, by \eqref{eq:druet2.26} we obtain that 
\begin{eqnarray}\nonumber
 \int_{\Omega_p\setminus B_\xi(0,R)}|x|^2v_p^{p^*}dv_{g_p}& \le & C R^{p^*\beta_p}\int_{\Omega_p\setminus B_\xi(0,R)}|x|^{2-p^*\beta_p}dv_\xi\\ \nonumber
 & \le & C (1+C\mu_p^2)R^{p^*\beta_p}\int_{\R^n}|x|^{2-p^*\beta_p}dv_\xi\\ \nonumber
 &\le & C (1+C\mu_p^2)\omega_{n-1}R^{p^*\beta_p}\int_{R}^\infty \rho^{2-p^*\beta_p}\rho^{n-1}d\rho\\ \nonumber
 &\le & C (1+C\mu_p^2)\omega_{n-1}R^{p^*\beta_p}\left( \dfrac{\rho^{n+2-p^*\beta_p}}{n+2-p^*\beta_p}\right)\Big|_{R}^\infty\\  \label{Eq:Prop1Step5Estimates}
 & = & C (1+C\mu_p^2)\omega_{n-1}R^{n+2}\tilde\gamma_{p,n}\to 0,
\end{eqnarray}
where $\tilde\gamma_{p,n}:=\dfrac{1}{p^*\beta_p-n-2}$.
Using \eqref{eq:druet2.17}, \eqref{Eq:Prop1Step5Mean}, and \eqref{Eq:Prop1Step5Estimates} we conclude that
\begin{equation}%\label{eq:druet2.29} 
\nonumber
\int_{\Omega_p}\tilde v_p^{\frac n{n-1}}dv_\xi
=1+\dfrac{Sc_g(x_p)}{6n(n+2)}\omega_{n-1}R_0^{n+2}\mu_p^2+o(\mu_p^2),
\end{equation}
and the expression on the right hand side of \eqref{eq:druet2.27} becomes 
\begin{equation}\label{eq:druet2.29}
\left(\int_{\Omega_p}\tilde v_p^{\frac n{n-1}}dv_\xi\right)^{\frac{n-1}{n}}
=1+\dfrac{(n-1)Sc_g(x_p)}{6n^2(n+2)}\omega_{n-1}R_0^{n+2}\mu_p^2+o(\mu_p^2).
\end{equation}
Denote by $l_2$ the limit of the scalar curvature function at $x_p$, i.e., 
\begin{equation}\label{Eq:PropSufficientl_2}
l_2:=\lim_{p\to1^+}Sc_g(x_p)\in\R,
\end{equation}
which exists and is finite because in strong bounded geometry $|Sc_g(x)|$ is uniformly bounded with respect to $x\in M$.
A fact that will be used often in the sequel is that thanks to the hypothesis of $C^2$ convergence of the metric to the metric at infinity we have
\begin{equation}\label{Eq:SamelimitSc}
\lim_{p\to1^+}Sc_g(x_p)=\lim_{p\to1^+}Sc_g(x_{0,p})=l_1,
\end{equation}
since $d_g(x_p, x_{0,p})\le r_p\to0$, when $p\to 1^+$. 
But observe that in the proof of Proposition \ref{Proposition:SufficientWeak1} this no longer true.
By the Cartan expansion of $g_p$ at $0$, since $r_p \to 0$ as $p\to1$, we also have
%%%%
%%%% Druet have negative signature in Rm
$$
|\nabla_\xi \tilde v_p|_{\xi}^p=|\nabla_{g_p} \tilde v|_{g_p}^p\left(1+\dfrac{\mu_p^2}6|\nabla_{g_p}\tilde  v_p|_{g_p}^{-2} Rm_g(x_p)(\nabla_{g_p}\tilde  v_p,x,x,\nabla_{g_p}\tilde  v_p)+o(\mu_p^2|x|^2)\right),
$$
 where $Rm_g$ denotes the Riemann curvature of $g$ in the $\exp_{x_p}$-map. Then, using
\eqref{eq:druet2.28}, we get
 
\begin{align}\label{eq:druet2.30}
\int_{\Omega_p} |\nabla_{\xi} \tilde v_p|_{\xi}dv_{\xi}
=&\int_{\Omega_p}|\nabla_{g_p} \tilde v_p|_{g_p}dv_{g_p}+\dfrac{\mu_p^2}6Ric_g (x_p)_{ij} \int_{\Omega_p}x^ix^j |\nabla_{\xi} \tilde v_p|_{\xi} dv_{\xi} \notag \\
&+\dfrac{\mu_p^2}6\int_{\Omega_p}|\nabla_{g_p} \tilde v_p|_{g_p}^{-1}Rm_g(x_p)(\nabla_{g_p}\tilde  v_p,x,x,\nabla_{g_p}\tilde  v_p)dv_{g_p} \notag\\
&+o\left(  \mu_p^2 \int_{\Omega_p}|x|^2|\nabla_{g_p} \tilde v_p|_{g_p}dv_{g_p}\right).
%+o\left(  \mu_p^2 \int_{\Omega_p}|x|^2|\nabla \tilde v_p|_{\xi}dv_{\xi}\right)
 \end{align}
 Let us now estimate the different terms of \eqref{eq:druet2.30}. First, by equation \eqref{eq:druet2.10} and 
relation \eqref{eq:druet2.5}, we have
\begin{eqnarray*}
\int_{\Omega_p} |\nabla_{g_p} \tilde v_p|_{g_p}dv_{g_p} & = & \tilde\gamma_{p,n}^*\int_{\Omega_p} v_p^{\frac{n(p-1)}{n-p}}|\nabla_{g_p} v_p|_{g_p}dv_{g_p}\\
& \le & \tilde\gamma_{p,n}^*\left( \int_{\Omega_p} v_p^{p^*}dv_{g_p} \right)^{\frac{p-1}{p}}\left( \int_{\Omega_p} |\nabla_{g_p}  v_p|^p_{g_p}dv_{g_p} \right)^{\frac1p}\\
& \le & \tilde\gamma_{p,n}^*\left(\lambda_p-\alpha_p\mu_p\|v\|_p^2\right)^{\frac12}\\
 & = &\tilde\gamma_{p,n}^*\lambda_p^{\frac12} \left(1-\alpha_p\mu_p\lambda_p^{-1}\|v\|_p^2\right)^{\frac12}\\
& \le & K(n,1)^{-1}\left( 1-\alpha_p\mu_p^2\lambda_p^{-1}\|v_p\|_p^2 \right)^{\frac12},
\end{eqnarray*}
where $\tilde\gamma_{p,n}^*:=\frac{p(n-1)}{n-p}$.
Since, by \eqref{eq:druet2.17} and \eqref{eq:druet2.26}, $\|v_p\|_{g_p,p} = 1 + o (1)$, we get
\begin{equation}\label{eq:druet2.31}
 \int_{\Omega_p} |\nabla_{g_p} \tilde v_p|_{g_p}dv_{g_p}\le 
 K(n,1)^{-1}-\frac{\alpha_p}2K(n,1)\mu_p^2+o(\mu_p^2).
\end{equation} 
By Holder's inequalities, we have
\begin{align*}
\int_{\Omega_p} |x|^2|\nabla_{g_p} \tilde v_p|_{g_p}dv_{g_p}&= 
\tilde\gamma_{p,n}^*\int_{\Omega_p} v_p^{\frac{n(p-1)}{n-p}}|x|^2|\nabla_{g_p} v_p|_{g_p}dv_{g_p}\\
&\le \tilde\gamma_{p,n}^*\left( \int_{\Omega_p} v_p^{p^*}dv_{g_p} \right)^{\frac{p-1}{p}}\left( \int_{\Omega_p} |x|^{2p}|\nabla_{g_p}  v_p|^p_{g_p}dv_{g_p} \right)^{\frac1p}\\
&=\tilde\gamma_{p,n}^*\left( \int_{\Omega_p} |x|^{2p}|\nabla_{g_p} v_p|^p_{g_p}dv_{g_p} \right)^{\frac1p}.
\end{align*}
Multiplying the equation \eqref{eq:druet2.10} by $|x|^{2p}v_p$ and integrating by parts, one gets
\begin{eqnarray*}%\label{eq:druet2.10}
\int_{\Omega_p} |\nabla_{g_p} v_p|^{p-2}\langle \nabla_{g_p}(|x|^{2p}v_p),\nabla_{g_p} v_p\rangle dv_{g_p} & = &
C_p^{-1}\lambda_p\int_{\Omega_p}|x|^{2p}v_p^{p^*}dv_{g_p}\\
& - & C_p^{-1}\alpha_p\mu_p^2 \|v_p\|_p^{2-p}\int_{\Omega_p}|x|^{2p}v_p^{p}dv_{g_p}.
\end{eqnarray*}
By \eqref{eq:druet2.26}, every term on the right hand side of the preceding inequality is uniformly bounded with respect to $p$, then we conclude that
\begin{equation}\label{eq:druet2567}
\int_{\Omega_p} |\nabla_{g_p} v_p|^{p-2}\langle \nabla_{g_p}(|x|^{2p}v_p),\nabla_{g_p} v_p\rangle dv_{g_p}\le C,
\end{equation}
for some $C>0$ that does not depend on $p$. Furthermore by Cauchy-Schwarz's inequality we get that
\begin{eqnarray*}
\int_{\Omega_p} |\nabla_{g_p} v_p|^{p-2}\langle \nabla_{g_p}(|x|^{2p}v_p),\nabla_{g_p} v_p\rangle dv_{g_p}
& = & \int_{\Omega_p} |x|^{2p}|\nabla_{g_p} v_p|^p\\
& + & \int_{\Omega_p} |\nabla_{g_p} v_p|^{p-2}v_p\langle \nabla_{g_p}(|x|^{2p}),\nabla_{g_p} v_p\rangle dv_{g_p}\\
& \ge & \int_{\Omega_p} |x|^{2p}|\nabla_{g_p} v_p|^pdv_{g_p}\\
& - & 2p\int_{\Omega_p} |x|^{2p-1}\nabla_{g_p}(|x|) |\nabla_{g_p} v_p|^{p-1}v_p  dv_{g_p}\\
&= & \int_{\Omega_p} |x|^{2p}|\nabla_{g_p} v_p|^pdv_{g_p}\\
& - & 2p\int_{\Omega_p} |x|^{2p-1} |\nabla_{g_p} v_p|^{p-1}v_p dv_{g_p}.
\end{eqnarray*}
Therefore we are lead to
\begin{eqnarray*}
 \int_{\Omega_p} |x|^{2p}|\nabla_{g_p} v_p|^p_{g_p}dv_{g_p} & \le & \int_{\Omega_p} |\nabla_{g_p} v_p|^{p-2}_{g_p}\langle \nabla_{g_p}(|x|^{2p}v_p), \nabla_{g_p} v_p\rangle_{g_p}dv_{g_p}\\
 & + & C\int_{\Omega_p} |x|^{2p-1}|\nabla_{g_p} v_p|^{p-1}_{g_p}v_pdv_{g_p}\\
 & \le & C\\ 
 & + & C\left(  \int_{\Omega_p} |x|^{2p}|\nabla_{g_p} v_p|^{p}_{g_p}dv_{g_p}   \right) ^{\frac{p-1}p}
  \left( \int_{\Omega_p} |x|^{p} v_p^{p}dv_{g_p}\right)^{\frac1p},
 \end{eqnarray*}
where $C$ denotes some constants independent of $p$. By  \eqref{eq:druet2.26} we see easily that 
$ \int_{\Omega_p} |x|^{p} v_p^{p}dv_{g_p}$, is uniformly bounded with respect to $p$. Then by Young's inequalities, one deduces that
\begin{align*}
 \int_{\Omega_p} |x|^{2p}|\nabla_{g_p} v_p|^p_{g_p}dv_{g_p}\le& C+C\left(  \int_{\Omega_p} |x|^{2p}|\nabla_{g_p} v_p|^{p}_{g_p}dv_{g_p}   \right) ^{\frac{p-1}p}\\
 \le& C+ \frac{C^p}{p}+\frac{p-1}{p}\left(  \int_{\Omega_p} |x|^{2p}|\nabla_{g_p} v_p|^{p}_{g_p}dv_{g_p}   \right) ^{\frac{p-1}p\frac p{p-1}},
\end{align*}
and so 
\begin{equation*}
\int_{\Omega_p} |x|^{2p}|\nabla_{g_p} v_p|^p_{g_p}dv_{g_p}\le \left(1-\frac{p-1}{p}\right)^{-1}\left(C+ \frac{C^p}{p}\right)\le\tilde C,
\end{equation*}
with $\tilde C>0$ independent of $p$.
That is
\begin{equation}\label{eq:druet2.32}
\int_{\Omega_p} |x|^{2p}|\nabla_{g_p} v_p|^{p}_{g_p}dv_{g_p}=O(1).
\end{equation}
Now for some $R>R_0$, we get readily  by \eqref{eq:druet2.18} that
\begin{eqnarray*}
\int_{\Omega_p} |\nabla_{\xi}\tilde v_p|_{\xi}x^ix^jdv_{\xi} & = & O\left( \int_{\Omega_p\setminus B_\xi(0,R)} |x|^{2}|\nabla_{\xi}\tilde v_p|_{\xi}dv_{\xi}   \right)\\
& + & \int_{\partial B_\xi(0,R_0)} x^ix^jd\sigma_{\xi}+o(1).
\end{eqnarray*}
By H\"older inequality we obtain that 
$$
\int_{\Omega_p\setminus B_{\xi}(0,R)}|x|^2|\nabla_{g_p} \tilde v_p|_{g_p}dv_{g_p} \le
\dfrac{p(n-1)}{n-p}\left( \int_{\Omega_p\setminus B_{\xi}(0,R)} |x|^{2p}|\nabla_{g_p} v_p|^p_g dv_{g_p}\right)^{\frac1p}.
$$
Multiplying the equation \eqref{eq:druet2.10} by $|x|^{2p}v_p$, integrating over $\Omega_p\setminus B_{\xi}(0,R):=\Omega_p^*$, and using Cauchy-Schwarz, H\"older inequality and later by Young inequality we obtain 
\begin{eqnarray*}
 \int_{\Omega_p^*} |x|^{2p}|\nabla_{g_p} v_p|^p_{g_p}dv_{g_p} &\le & \int_{\Omega_p^*} |\nabla_{g_p} v_p|^{p-2}_{g_p}\langle \nabla_{g_p}(|x|^{2p}v_p), \nabla_{g_p} v_p\rangle_{g_p}dv_{g_p}\\
 & + & 2p\int_{\Omega_p^*} |x|^{2p-1}|\nabla_{g_p} v_p|^{p-1}_{g_p}v_pdv_{g_p}\\
%\int_{\Omega_p^*}ertext{and now by H\"older inequality and later by Young inequality,}
& \le & \int_{\Omega_p^*} |\nabla_{g_p} v_p|^{p-2}_{g_p}\langle \nabla_{g_p}(|x|^{2p}v_p), \nabla_{g_p} v_p\rangle_{g_p}dv_{g_p}\\ 
& + & 2p\left(  \int_{\Omega_p^*} |x|^{2p}|\nabla_{g_p} v_p|^{p}_{g_p}dv_{g_p}   \right) ^{\frac{p-1}p}
  \left( \int_{\Omega_p^*} |x|^{p} v_p^{p}dv_{g_p}\right)^{\frac1p}\\
& \le & \int_{\Omega_p^*} |\nabla_{g_p} v_p|^{p-2}_{g_p}\langle \nabla_{g_p}(|x|^{2p}v_p), \nabla_{g_p} v_p\rangle_{g_p}dv_{g_p}\\
& + & 2(p-1)\int_{\Omega_p^*} |x|^{2p}|\nabla_{g_p} v_p|^{p}_{g_p}dv_{g_p}+2 \int_{\Omega_p^*} |x|^{p} v_p^{p}dv_{g_p}.
 \end{eqnarray*}
At last we obtain that
\begin{eqnarray*}
0\le(3-2p)\int_{\Omega_p^*} |x|^{2p}|\nabla_{g_p} v_p|^p_{g_p}dv_{g_p} & \le & \int_{\Omega_p^*} |\nabla_{g_p} v_p|^{p-2}_{g_p}\langle \nabla_{g_p}(|x|^{2p}v_p), \nabla_{g_p} v_p\rangle_{g_p}dv_{g_p}\\
& + & 2 \int_{\Omega_p^*} |x|^{p} v_p^{p}dv_{g_p}.
 \end{eqnarray*}
But when $p\to1$, by  \eqref{eq:druet2.26}, the terms on the right hand side go to $0$, then we can conclude that
\begin{equation}\label{Eq:Prop1QuasiFinal}
\int_{\Omega_p\setminus B_\xi(0,R)} |x|^{2}|\nabla_{\xi} v_p|_{\xi}dv_{\xi}\to 0.
\end{equation}
Thus for the second term on the right hand side of \eqref{eq:druet2.30} we see that
\begin{equation}\label{eq:druet2.33}
 \lim _{p\to1}Ric_g(x_p)_{ij}\int_{\Omega_p} |\nabla_{\xi}\tilde v_p|_{\xi}x^ix^jdv_{\xi}=\frac{\omega_{n-1}}nR_0^{n+1}l_1.%S_g(x_0).
\end{equation}
Now, we look at the third term on the right hand side of \eqref{eq:druet2.30}.
Since $\nabla V_p$, $V_p$ as in \eqref{eq:druet2.19}, and $x$ are pointwise colinear vector fields, we have
\begin{equation}\label{Eq:Prop1QuasiFinal1}
Rm_g(x_p)(\nabla_{g_p}\tilde v_p,x,x,\nabla_{g_p}\tilde v_p)\le C|x|^2|\nabla_{\xi} \tilde v_p|_\xi|\nabla_{\xi}(\tilde v_p-V_p) |_\xi.
\end{equation}

Now by %\eqref{eq:druet2.10}, 
\eqref{Eq:Prop1QuasiFinal1}, %\eqref{eq:druet2.20}-\eqref{eq:druet2.26}, 
integrating over $\Omega_p\cap B_\xi(0,R):=\hat{\Omega}_p$ we have 
\begin{small}
\begin{eqnarray*}
\int_{\hat{\Omega}_p}|\nabla_{g_p}\tilde v_p|_{g_p}^{-1}Rm_g(x_p)(\nabla_{g_p}\tilde v_p,x,x,\nabla_{g_p}\tilde v_p)dv_{g_p} & \le & \int_{\hat{\Omega}_p} C|\nabla_{g_p}\tilde v_p|_{g_p}^{-1}|x|^2|\nabla_{\xi} \tilde v_p|_\xi|\nabla_{g_p}(\tilde v_p-V_p) |_\xi dv_{g_p}\\
& \le & C_R\left(1+C\mu_p^2\right)\int_{\hat{\Omega}_p}|\nabla_{_{\xi}}(\tilde v_p-V_p) |_\xi dv_\xi.
\end{eqnarray*}
\end{small}
This last inequality combined with \eqref{eq:druet2.20} yields to
\begin{equation}\label{Eq:Prop1QuasiFinal2}
\lim _{p\to1}\int_{\Omega_p\cap B_\xi(0,R)}|\nabla_{g_p}\tilde v_p|_{g_p}^{-1}Rm_g(x_p)(\nabla_{g_p}\tilde v_p,x,x,\nabla_{g_p}\tilde v_p)dv_{g_p}=0.
\end{equation}
We want to estimate the integral of the same integrand function of \eqref{Eq:Prop1QuasiFinal2} but outside $B_\xi(0,R)$, for this we have
\begin{tiny}
\begin{eqnarray}\nonumber
\int_{\Omega_p\setminus B_\xi(0,R)}|\nabla_{g_p}\tilde v_p|_{g_p}^{-1}Rm_g(x_p)(\nabla_{g_p}\tilde v_p,x,x,\nabla_{g_p}\tilde v_p)dv_{g_p}& \le &\int_{\Omega_p\setminus B_\xi(0,R)}\Lambda_2|x|^2|\nabla_{g_p}\tilde v_p|_{g_p}^{-1}|\nabla_{g_p}\tilde v_p|_{g_p}^2dv_{g_p}\\ \nonumber
& \le & \Lambda_2(1+C\mu_p^2)\int_{\Omega_p\setminus B_\xi(0,R)}|x|^{2}|\nabla_{g_p}\tilde v_p|_{\xi}dv_{\xi}\\ \label{Eq:Prop1QuasiFinal3}
& \stackrel{\eqref{Eq:Prop1QuasiFinal}}{\to} & 0.
\end{eqnarray}
\end{tiny}
Combining \eqref{Eq:Prop1QuasiFinal2} and \eqref{Eq:Prop1QuasiFinal3} we conclude that 
\begin{equation}\label{eq:druet2.34}
 \lim _{p\to1}  \int_{\Omega_p} 
 |\nabla_{g_p}\tilde v_p|_{g_p}^{-1}Rm_g(x_p)(\nabla_{g_p}\tilde v_p,x,x,\nabla_{g_p}\tilde v_p)dv_{g_p}=0.
 \end{equation}
Finally, substituting in \eqref{eq:druet2.27}, using \eqref{Eq:SamelimitSc}, and \eqref{eq:druet2.29}-\eqref{eq:druet2.32}, we obtain,
\begin{eqnarray*}
1+\dfrac{(n-1)l_2}{6n^2(n+2)}\omega_{n-1}R_0^{n+2}\mu_p^2 & + & o(\mu_p^2)\\
& \le & K(n,1)\left[K(n,1)^{-1}-\frac{\alpha_p}2K(n,1)\mu_p^2\right]\\
& + & K(n,1)\left[\frac{\omega_{n-1}}{6n}R_0^{n+1}l_2\mu_p^2\right]+o(\mu_p^2)\\
& = & 1-\frac{\alpha_p}2K(n,1)^2\mu_p^2\\ 
& + & \frac{K(n,1)\omega_{n-1}}{6n}R_0^{n+1}l_2\mu_p^2+o(\mu_p^2).
\end{eqnarray*}
Since $\dfrac{\omega_{n-1}}n=\dfrac{1}{R_0^n}$, and $K(n,1)=\dfrac1n\left(\dfrac n{\omega_{n-1}}\right)^{\frac1n}=\dfrac{R_0}{n}$, a straightforward computation leads to
$$
\frac{K(n,1)^2}2\left(\alpha_p-\frac n{n+2}l_2
\right)\mu_p^2+o(\mu_p^2)\le0.
$$
This gives the desired contradiction by letting $p$ go to $0$, recalling here that $l_1=l_2$ by \eqref{Eq:PropSufficientl_2} it holds
\begin{equation}\label{Eq:PropositionSufficientFinal0}
\frac n{n+2}l_1- \varepsilon_0+\frac n{n+2}l_2=\lim_{p\to1^+}\alpha_p-\frac n{n+2}l_2= \varepsilon_0 > 0.
\end{equation}
This ends the proof of Proposition \ref{Proposition:Sufficient}.
\end{proof} 
We are now ready to accomplish the proof of our global comparison theorem for small diameters in $C^2$-locally asymptotically strong bounded geometry. We use the same argument used in \cite{DruetPAMS}, for completeness's sake we write the details here as pointed out to us by Olivier Druet in a private communication.  
\begin{proof}[Proof of Theorem \ref{Thm:DruetPAMS1}] The Proposition at page $2353$ of \cite{DruetPAMS} rewritten in this text as Proposition \ref{Prop:Druet} says that for any $\varepsilon > 0$, there exists $r_\varepsilon=r_\varepsilon(x_0, M, g)> 0$ such that if
$\Omega\subset B_g (x_0, r_\varepsilon)$, then
$$
V_g(\Omega)^{2\frac{n-1}n}\le K(n,1)^2A_g(\partial\Omega)^2 + K(n,1)^2\left(\dfrac n{n+2}S_g(x_0) +\varepsilon \right) V_g(\Omega)^2.
$$
By assumption we know that $S_g(x_0)<n(n-1)k_0$, so that applying the preceding inequality with $$\varepsilon=\dfrac n{2(n+2)}\left[n(n-1)k_0-S_g(x_0)\right],$$ we get that there exists $r>0$, $r_\varepsilon(x_0, M, g)\ge r>0$ such that if $\Omega\subset B_g(x_0,r) $, then 
\begin{eqnarray}
 V_g(\Omega)^{2\frac{n-1}n} & \le & K(n,1)^2A_g(\partial\Omega)^2\\ \nonumber
& + & K(n,1)^2\left(\dfrac n{n+2}n(n-1)k_0-\dfrac n{2(n+2)}\varepsilon_0 \right) V_g(\Omega)^2.
\end{eqnarray}
where $\varepsilon_0=n(n-1)k_0-S_g(x_0)>0$ fixed. Now let $B_v$ be a small ball in the model space $(\mathbb{M}_{k_0},g_{k_0})$ of constant sectional curvature $k_0$ and of volume $v$, for any $V_0>0$ (small enough in the case of the sphere, i.e., $k_0>0$) there exists $C_0=C_0(n,k,V_0)>0$ such that for 
balls of volume $0\le v\le V_0$ it holds
\begin{eqnarray}
V_{g_{k_0}}(B_v)^{2\frac{n-1}n} & \ge & K(n,1)^2A_{g_{k_0}}(\partial B_v)^2\\ 
& + & K(n,1)^2 \dfrac n{n+2}n(n-1)k_0 V_{g_{k_0}}(B_v)^2\\ 
& - & C_0v^{2\frac n{n+2}}.
\end{eqnarray}
If we assume that $V_g(\Omega)=V_{g_0}(B_v)=v\le V_0$, %the by the ball $B_V$ 
we get that
\begin{eqnarray*}
 K(n,1)^2A_{g_0}(\partial B_v)^2 & + &K(n,1)^2\dfrac n{n+2}n(n-1)k_0v^2 - C_0 v^{2\frac {n+2}n}\\ 
 & \le & v^{2\frac {n-2}n}\\
& \le & K(n,1)^2A_g(\partial\Omega)^2\\
& + & K(n,1)^2\left( 
\dfrac n{n+2}n(n-1)k_0-\dfrac n{2(n+2)}\varepsilon_0 \right) v^2
\end{eqnarray*}
that is,
\begin{align*}
A_{g_0}(\partial B_v)^2\le A_g(\partial\Omega)^2+C_0 K(n,1)^{-2} v^{2\frac {n+2}n}-\dfrac n{2(n+2)}\varepsilon_0v^2.
\end{align*}

If we choose $v <V_1 <V_0<\min\{1, V_{g_{k_0}}\left(\mathbb{M}^n_{k_0}\right)\}$ \footnote{$\min\{1,+\infty\}$ is assumed to be equal at $1$.}, with the property that
$$
C_0(n,k_0)K(n,1)^{-2} V_1^{2\frac {n+2}n}-\dfrac n{2(n+2)}\varepsilon_0V_1^2<0,
$$
which is always possible to find, then we get 
$$
A_{g_{k_0}}(\partial B_v)< A_g(\partial\Omega).
$$
Thus there exists $V_1=V_1(n,k_0, V_0, Sc_g(x_0))>0$ such that if $\Omega\subset B_g(x_0,r)$ with volume  $V_{g_0}(B_v)=v<V_1$, then the comparison inequality \eqref{Eq:Thm1DruetPAMs} of the theorem holds. Now, up to lower a little bit $r$ (depending on curvatures of $M$), we are sure that any domain $\Omega\subset B_g(x_0,r)$ has volume less than that of the ball $B_g(x_0,r)$ which can be chosen to be less than $V_1$ and the theorem is proved.
At the end of this proof we understand the subtle fact that a lower bound on $r_x$ of Theorem $1$ of \cite{DruetPAMS} does not depend explicitly on the $r_\varepsilon$ of the Proposition at page $2353$ but we need to prove the existence of $r_\varepsilon$ to prove the existence of $r_x$.
\end{proof}
We are thus led to the following version of Theorem $1$ of \cite{DruetPAMS}, namely our Theorem \ref{Res:Theorem1}, in which an uniform estimate of a lower bound on $r_x$ is obtained provided $M$ is of $C^2$-locally asymptotically strong bounded geometry at infinity.
\begin{proof}[Proof of Theorem \ref{Res:Theorem1}] We proceed as in the proof of the preceding theorem using our Proposition \ref{Proposition:Sufficient} instead of the Proposition at page 2353 of \cite{DruetPAMS}. This gives the existence of an uniform $r_\varepsilon(M,g)>0$ independent of $x_0$. With the help of Bishop-Gromov we have for every $0\le r\le 1$. 
$$V_g(B_g(x_0,r))\le V_{g_k}(B_{g_k}(x_0,r))\le C(n,k)r^n,$$ 
where $k$ is a lower bound on the Ricci curvature of $(M,g)$. So we can take a radius $r=r(n,k,k_0, V_0, S_g)>0$ such that 
\begin{equation}\label{Eq:RadiusEstimates}
r<\sqrt[n]{\frac{V_1(n,k_0, V_0, S_g)}{C(n,k)}},
\end{equation}
where $V_1$ is as in the proof of Theorem \ref{Thm:DruetPAMS1} with $\varepsilon_0$ replaced by $\tilde \varepsilon_0=n(n-1)k_0-S_g>0$. Observing that we can take for example $V_0<\min\{1, V_{g_{k_0}}\left(\mathbb{M}^n_{k_0}\right)\}$ fixed we obtain $r=r(n,k,k_0,S_g)>0$.
\end{proof}
The following proposition have a weaker conclusion with respect to the previous one but that holds for manifolds with weaker assumptions, namely for manifolds with just strong bounded geometry. 
\begin{ResProp}\label{Proposition:SufficientWeak1}  Let $(M, g)$ be a complete Riemannian manifold of dimension $n\ge2$ with strong bounded geometry. For any $\varepsilon>0$ there exists $r_{\varepsilon}=r_{\varepsilon}(M,g)>0$ such that for any point $x_0\in M$, any function $u\in C^{\infty}_c(B_g(x_0,r_{\varepsilon}))$, we have
\begin{equation}\label{Eq:RespropSufficientWeak1}
||u||^2_{\frac n{n-1},g}\le K(n,1)^2\left(||\nabla_g u||_{1,g}^2+\alpha_\varepsilon||u||_{1,g}^2\right),
\end{equation} 
where $\alpha_{\varepsilon}=\frac n{n+2}S_g+\varepsilon$ with $S_g:=\sup_{x\in M}\{Sc_g(x)\}\in\R$.
\end{ResProp}
\begin{Rem} We notice that the constant $r_{\varepsilon}=r_{\varepsilon}(M)>0$ is obtained by contradiction and that the proof does not give an explicit effective lower bound on it. 
\end{Rem} 
\begin{proof}[First Proof of Proposition \ref{Proposition:SufficientWeak1}] This proof is obtained from the proof of Proposition \ref{Proposition:Sufficient} replacing in that proof $\alpha_p$ by the constant independent of $p$ defined by $\alpha_{\varepsilon_0}:=\dfrac n{n+2} S_g+\varepsilon_0$. 
Then we find a sequence $(\tilde u_p)$ solving the same partial differential equations that the analogous $u_p$ solve in the proof of Proposition \ref{Proposition:Sufficient} but replacing $\alpha_p$ with $\alpha_{\varepsilon_0}$, after we take as $x_p$ a point of maximum of $\tilde u_p$ and observe that up to a subsequence we can assume
\begin{equation}\label{Eq:Prop1alphap}
\lim_{p\to1^+}Sc_g(x_p)=\frac n{n+2}l_2\le S_g,
\end{equation} 
for some $l_2\in [S_{inf}, S_g]$, where $S_{inf,g}:=\inf\{S_g(x):x\in M\}$. It is worth to note that $S_{inf,g}$ and $S_g$ are finite real numbers, because $S_g$ is bounded from below by $n(n-1)k$ and from above by $n(n-1)k_0$. 
Then exactly in the same way we are lead to the analog of \eqref{Eq:PropositionSufficientFinal0}, that is 
\begin{equation} 
\alpha_{\varepsilon_0}-\frac n{n+2}l_2=\frac n{n+2}S_g+\varepsilon_0-\frac n{n+2}l_2\ge\varepsilon_0> 0.
\end{equation}
This last inequality giving the desired contradiction proves our Proposition \ref{Proposition:SufficientWeak1}.
\end{proof} 
We can now proceed to the proof of Theorem \ref{Res:Theorem1.1}.
\begin{proof}[Proof of Theorem \ref{Res:Theorem1.1}] At this point it is easy to remark that Theorem \ref{Res:Theorem1.1} follows from Proposition \ref{Proposition:SufficientWeak1} exactly in the same way as Theorem \ref{Res:Theorem1} follows from Proposition \ref{Proposition:Sufficient}. We leave the details of the proof to the reader.
\end{proof}
In the remaining part of this section we give an alternative proof of Proposition \ref{Proposition:SufficientWeak1}, using smoothing via the Ricci flow and results of \cite{Shi} and \cite{Kapovitch}. Let us denote $(M,\tilde g)$ a complete Riemannian manifold with $|Sec_{\tilde g}|\le\Lambda$. Consider the Ricci flow of $\tilde g_t$ with initial data $\tilde g=\tilde g_0$  
 \begin{eqnarray}\label{eq:kapovitch0}
 \frac{\partial \tilde g_t}{\partial t} =  -2 Ric(\tilde g_t).
\end{eqnarray}
\begin{Thm}[see Theorem $1.1$ of \cite{Shi}, compare also \cite{BMOR} and \cite{Hamilton}]\label{Thm:Shi}
%It is known (see [Ham82], [Shi89]) that \eqref{eq:kapovitch1} has a solution on $[0,T]$ for some $T > 0$, and that (see [BMOR84], [Shi89])  the solution smoothes out the metric $\tilde g_t$ satisfies
 Under the above assumptions on $(M,\tilde g)$ there exists a constant $T=T(n,\Lambda)>0$ such that \eqref{eq:kapovitch0} with initial condition $\tilde g_0=\tilde g$ has a smooth solution $\tilde g_t$ for a short time $0\le t\le T$. Moreover $\tilde g_t$ satisfies the following estimates    
 %It is known  that \eqref{eq:kapovitch1} has a solution on $[0,T]$ for some $T > 0$, and that   the solution smoothes out the metric, called $\tilde g_t$ satisfies (see \cite{Kapovitch} and the references therein)
 \begin{equation}\label{Eq:Shi}
 e^{-c(n,\Lambda)t}\tilde g \le\tilde  g_t \le e^{c(n,\Lambda)t}\tilde g,
 \end{equation}
 \begin{equation}\label{Eq:Shi0}
  |\nabla-\nabla_t| \le c(n, \Lambda)t, 
\end{equation}
and for any integer $m\ge0$, there exists constants $c_m=c(n, m, \Lambda)>0$ such that  
\begin{equation}\label{Eq:Shi1}
\sup_{x\in M}  |\nabla^mR_{ijkl}(x,t)| \le c_mt^{-\frac m2},\:\:\:\: 0<t\le T(n,\Lambda).
\end{equation}

Moreover %, by [Shi89], 
the sectional curvature of $\tilde g(t)$ satisfies  
 \begin{eqnarray}\label{eq:kapovitch3}
|Sec_{\tilde g_t}| \le C(n,T), \:\:\:\: 0\le t\le T(n,\Lambda).
 \end{eqnarray}
for some constant $C(n,T)>0$. 
\end{Thm}
Then we can assume true the hypotheses of the existence of the constant $C(n,\Lambda, T)$ in the theorem below.
\begin{Thm}[Proposition at page $260$ of \cite{Kapovitch}]\label{thm:kapovitch1}
 Let $(M^n,\tilde g)$ be a complete Riemannian manifold and suppose that   $|Sec_{\tilde g}|\le \Lambda$, for some positive constant $\Lambda$. Consider the Ricci flow of $\tilde g$ given by
 \begin{eqnarray}\label{eq:kapovitch1}
 \frac{\partial \tilde g_t}{\partial t} =  -2 Ric(\tilde g_t).
\end{eqnarray}
Then sectional curvatures satisfy the following relation
\begin{small}
\begin{eqnarray}\label{eq:kapovitch2}
 \inf_{x\in M}\left\{ Sec_{\tilde {g}_t}(x)\right\}-C(n,T)t\le Sec_{\tilde g_t}(x) \le \sup_{x\in M} \left\{ Sec_{\tilde {g}_t}(x)\right\} +C(n,T)t,
\end{eqnarray}
\end{small}
where the scalar curvature satisfies $|Sec_{\tilde g_t}|\le C(n,t)$.
\end{Thm}
\begin{Rem} If $(M,g)$ have $C^2$-locally asymptotically strong bounded geometry, then using Proposition \ref{Proposition:Sufficient} it is trivial to check that $(M,g)$ satisfies
\eqref{Eq:RespropSufficientWeak1}. 
\end{Rem}
\begin{proof}[Second Proof of Proposition \ref{Proposition:SufficientWeak1}] For an arbitrary Riemannian metric $g$, over $M$ we define 
 \begin{eqnarray}
 \lambda_{1,r,x_0,g}:=\inf_{\substack{u\in C^\infty_c (B_g (x_0,r))\\u\not\equiv 0}}\left\{\dfrac{\|\nabla_g u\|_{1,g}^2 +\left(\frac n{n+2} S_g+\varepsilon_0\right) \|u\|_{1,g}^2}{\|u\|_{\frac n{n-1},g}^2}\right\}.
\end{eqnarray}
 By reduction to the absurd, we suppose that \eqref{Eq:RespropSufficientWeak1} above is false. Then there exists $\varepsilon_0>0$ such that for all $r>0$ there exist $x_{0,r}\in M$, $u_{x_{0,r}}\in C^\infty_c(B_g(x_{0,r},r))$  such that $$\|u\|_{\frac n{n-1},g}^2> K(n,1)^2\left[ \|\nabla_g u\|_{1,g}^2+\left(\frac n{n+2} S_g+\varepsilon_0\right)\|u\|_{1,g}^2\right].$$
 The last inequality is equivalent to say that  
\begin{equation}
\lambda_{1,r,x_{0,r},g}<K(n,1)^{-2}.
\end{equation}
  Take the Ricci flow $(M,\tilde g_t)$ with the initial data $\tilde g_0:=g$. Then by Theorem \ref{Thm:Shi}, i.e., Theorems $1.1$ and  $1.2$ of \cite{Shi} we get that $(M,\tilde g_t)$ is a smooth manifold with smooth metric $\tilde g_t$ that by \eqref{Eq:Shi1} satisfies $|\nabla_{\tilde g_t}^m R_{ijkl,\tilde g_t}|_{\tilde g_t},|\nabla_{g_t}^m R_{ijkl,\tilde g_t}|_{g_t}\le c_mt^{-\frac m2}$, for any integer $m\ge 0$. Moreover, by \eqref{Eq:Shi} we have $\tilde g_t\to \tilde g_0=g$ in $C^0$ topology. This combined with Klingenberg's Lemma (compare for instance Theorem III.2.4 of \cite{Chavel}) guarantees that for any $t\in[0, T]$ holds $inj_{M, \tilde g_t}\ge i_0=i_0(n,\Lambda_2, g)>0$ where $\Lambda_2$ is an upper bound on the sectional curvature, and so the manifold $(M,\tilde g_t)$ have strong bounded geometry. Furthermore for any integer $m\ge 0$ we have that $(M,\tilde g_t)$ belongs to $\mathcal{M}(n, m, \Lambda, r=r(n, \Lambda_2, i_0))$ (see Definition \ref{Def:BoundedGeometryPetersen} or page $308$ of \cite{Petersen}). 
Hence $(M,\tilde g_t)$ satisfy the hypothesis of Theorems $76$, $72$ of \cite{Petersen}, Theorems $1.2$, $1.3$ of \cite{Hebey}, and thus $(M,\tilde g_t)$ is smooth at infinity and have $C^{m,\alpha}$-locally asymptotically bounded geometry smooth at infinity, for any $m\in\N$. This implies that $(M,\tilde g_t)$ satisfies the assumptions of Proposition \ref{Proposition:Sufficient} and so trivially also the conclusion \eqref{Eq:RespropSufficientWeak1} of Proposition \ref{Proposition:SufficientWeak1}. Independently we observe that by Theorem \ref{thm:kapovitch1}, we have that there exists a sufficiently small $t_r\in ]0,T]$ such that $\lambda_{1,r,x_{0,r},\tilde g_{t_r}}<K(n,1)^{-2}$, but this in turns means that $(M,\tilde g_{t_r})$ does not satisfies the conclusion \eqref{Eq:RespropSufficientWeak1} of Proposition \ref{Proposition:SufficientWeak1}. In this way we get a contradiction, which indeed completes the proof of the proposition.
  \end{proof}   
\section{In mild bounded geometry isoperimetric regions of small volume are of small diameter}
In this section we work with just a fixed Riemannian metric $g$ defined on $M$. 
\begin{Lemme}[Lemma $3.2$ of \cite{Hebey}]\label{Lemma:Hebey3.2}
 Let (M,g) be a smooth, complete Riemannian $n$- dimensional manifold with weak bounded geometry. There exist two positive constants $C_{Heb}=C_{Heb}(n,k,v_0)>0$ and $\bar{v}:=\bar{v}(n,k,v_0)>0$, depending only on $n,k$, and $v_0$, such that for any open subset $\Omega$ of $M$ with smooth boundary and compact
closure, if   $V_g(\Omega)\le\bar{v}$, then $C_{Heb}V_g(\Omega)^{\frac{n-1}n}<A_g(\partial\Omega)$.
\end{Lemme}
\begin{Rem}\label{Rem:Hebey} After Theorem $1$ of \cite{FloresNardulli015} we know that we can extend the preceding lemma to an arbitrary finite perimeter set simply by approximating with open bounded with smooth boundary subsets having the same volume.
\end{Rem}
Let us introduce a crucial notion for the remaining part of this section.
\begin{Def}\label{Def:AppIsopSeq}
We say that a sequence $(D_j)$ of finite perimeter sets, $D_j\subseteq M$, with finite volume $V_g(D_j)\to 0$, is called an \textbf{approximate isoperimetric sequence}, if 
$$\lim_{j\to\infty}\dfrac{\P_g(D_j)}{V_g(D_j)^{\frac{n-1}{n}}}=\lambda,$$
where $\lambda:=\liminf_{v\to0^+}\dfrac{I_{M,g}(v)}{v^{(n-1)/n}}$.  
\end{Def}
\begin{Rem} %In case of $(M,g)$ of bounded geometry, Proposition $3.5$ of \cite{MJ} gives that $I_M\leq I_{\mathbb{M}^n_k}$ (see also Proposition $3.2$ of \cite{MondinoNardulli}) which implies immediately that 
Comparing with geodesic balls we have clearly that $\lambda\le c_n$, where $c_n$ is the Euclidean isoperimetric constant defined by $I_{\R^n}(v)=c_nv^{\frac{n-1}{n}},\:\forall v\in]0, V(M)[$. 
\end{Rem}
\begin{Rem} When $(M,g)$ have weak bounded geometry then $\lambda\ge C_{Heb}(n,k,v_0)>0$, because of Lemma $3.2$ of \cite{Hebey} reported here in Lemma \ref{Lemma:Hebey3.2} and the related Remark \ref{Rem:Hebey}.
\end{Rem}
\begin{Rem} When $(M,g)$ have strong bounded geometry then $\lambda=c_n$, this is an easy consequence of the Th\'eor\'eme of Appendice $C$ at page 531 of \cite{BerardMeyer}. We wrote an alternative proof of this last fact, based on Theorems \ref{Res:Theorem1.1} and \ref{Res:Theorem2}, in our Theorem \ref{Thm:BestSobolevConstant} below.
\end{Rem}
We recall here three well known lemmas (see for instance Corollary $2.1$ of \cite{NarAsian}) that we will use often in the sequel.
\begin{Lemme}\label{Cor:DoublingNonCollapsing}
Let $M^n$ be a complete Riemannian manifold with weak bounded geometry. Then for each $r>0$ there exists $c_1=c_1(n,k,r)>0$ such that $V_g(B_M(p,r))>c_1(n,k,r)v_0$, where $c_1(n,k,r)=\min\left\{\frac{r^n}{e^{\sqrt{(n-1)|k|}}}, 1\right\}$.
\end{Lemme}
\begin{Lemme} Let $M$ with weak bounded geometry. Then there exist two positive constants $C_1=C_1(n,k)>0$, $C_2=C_2(n,k)>0$ such that  for every $0<r<\bar r=\bar r(n,k):=\min\left\{1, e^{\frac{\sqrt{(n-1)|k|}}n}\right\}$ we have 
\begin{equation}
v_0C_1r^n\stackrel{doubling+noncollapsing}{\le} V_g(B_M(x,r))\stackrel{Bishop-Gromov}{\le} C_2r^n, 
\end{equation} 
where $C_1=C_1(n,k)=\frac{1}{e^{\sqrt{(n-1)|k|}}}$.
\end{Lemme}
\begin{Lemme} Let $M$ with weak bounded geometry. Then there exist two positive constants $\bar{v}_1=\bar{v}_1(n,k,v_0)>0$ and $C_3=C_3(n,k)>0$, such that  for every $0<v<\bar v_1$ we have 
\begin{equation}
\lambda\le\frac{I_M(v)}{v^{\frac{n-1}n}}\le C_3(n,k). 
\end{equation}
Here $\bar v_1:=\min\{1, \bar v\}$.  
\end{Lemme}
\begin{Lemme}\label{Lemma:Partitioning}
Let $M^n$ be a complete Riemannian manifold weak with bounded geometry. There exists a positive constant $N=N(n,k,v_0)>0$ such that, whenever $D$ is a finite perimeter set with finite volume and $0<R<\bar R=\bar R(n,k):=\min\left\{1, 2e^{\frac{\sqrt{(n-1)|k|}}n},\frac27\bar r\right\}$ there exists a partition  $(D_l)_{l}$ of $D$, i.e, $D=\mathring{\cup}_l D_l$, where every  $D_l$ is a set of finite perimeter contained in a ball of radius $R$ and such that 
\begin{equation}\label{Eq:Partitioning}
\left(\sum_l\P(D_l)\right)-\P(D)\le N(n,k,v_0)\dfrac{V(D)}{R}.
\end{equation}
\end{Lemme} 
\begin{proof}
Let $({p_l})_{l\in \N}$ be a sequence of points of $M$ such that 
$\left\{ B_M(p_l, \frac R4) \right\}$ is a maximal set of disjoint balls. It is straightforward to show that
$$
M=\bigcup_l B_M\left(p_l,\frac R2\right).
$$ Set $\mathcal{A}:=\{p_l\}_{l\in \N}$.  By coarea formula we can cut $D$ with a ball of radius $r_l$ centered at $p_l$, such that $\frac R2<r_l<R$ and
\begin{equation}
A(D\cap \partial B_M(p_l,r_l))\le \frac{2V(D)}{R}.
\end{equation}
Consider $D\setminus \left( \bigcup_l\partial B_M(p_l,r_l)\right)=\mathring{\bigcup}_l D_l$. Then there exists a constant $\tilde{N}=\tilde{N}(n,k,v_0)>0$ such that
$$
\left(\sum_l\P(D_l)\right)-\P(D)\le 4\tilde{N}\dfrac{V(D)}{R}.
$$
Note that by a simple combinatorial argument, $\tilde{N}$ could be taken as an upper bound of the greatest number of disjoint balls of radius $\frac{R}{4}$ contained in a ball of radius $\frac{7}{4}R$. This upper bound depends only on $n, k, v_0$ since for every $x\in M$ by our assumption $R<\bar R$ it holds 
\begin{eqnarray*}
\tilde{N}C_1(n,k)\left(\frac{R}{4}\right)^nv_0 & \le & \sum_{p_i\in B_M(x,\frac{7}{4}R)}V_g(B_M(p_i,\frac{R}{4}))\\
& \le & V_g\left(B_M\left(x, \frac{7}{4}R\right)\right)\\
&\stackrel{Bishop-Gromov}{\le} & V_{g_k}\left(B_{\mathbb{M}^n_k}\left(\frac{7}{4}R\right)\right)\le C_2(n,k)\left(\frac{7}{4}R\right)^n,
\end{eqnarray*}
where $C_1(n,k)=\frac{1}{e^{\sqrt{(n-1)|k|}}}$. %is easily deduced from Corollary $2.1$ of \cite{NarAsian}. 
Setting $N=4\tilde{N}$ we finish the proof of the  lemma.
\end{proof}
\begin{Lemme}\label{Lemma:Partitioning2}
Let $M^n$ be a complete Riemannian manifold with weak bounded geometry, and $D_j\subset M$ be a sequence of finite perimeter sets with finite volume. Then there exist a partition  of $D_j$ by finite perimeter sets of $D_j=\mathring{	\bigcup_l} D_{j,l}$ and a sequence of radii $R_j$,   
such that $\diam(D_{j,l})\le 2R_j$, with $\lim_{j\to\infty}\dfrac{V(D_j)^{1/n}}{R_j}=0$, $R_j\to0$, and
\begin{equation}\label{Eq:Partitioning2}
\limsup_{j\to\infty}\left[\left(\sum_l\P(D_{j,l})\right)-P(D_j) \right] \dfrac{1}{V(D_j)^{\frac{n-1}{n}}}=0.
\end{equation}
\end{Lemme}
\begin{proof} It is enough to apply \eqref{Eq:Partitioning} with $D=D_j$ and $R=R_j:=V(D_j)^{\alpha}$, with $0<\alpha<\frac{1}{n}$.
\end{proof} 
We are now ready to prove Theorem \ref{Res:Theorem2}.
\begin{proof}[Proof of Theorem \ref{Res:Theorem2}] Consider an arbitrary sequence of finite perimeter sets $\Omega_j$ such that $v_j:=V_g(\Omega_j)\to0$. By Lemma \ref{Lemma:Partitioning2} we can find a partition of $\Omega_j$ satisfying \eqref{Eq:Partitioning2}. For sufficiently large $j$ we have $R_j\le r$ where $r:=\frac{d}{2}$ and $d$ is given by Theorem \ref{Res:Theorem1.1}. Set $\eta_j:=N(n,k,v_0)\dfrac{v_j}{R_j}$, with $R_j>>v_j^{\frac{1}{n}}$, i.e., $\frac{v_j}{R_j}\to0$, when $j\to+\infty$, we obtain
\begin{eqnarray}\label{Eq:BM}
\frac{I_{\mathbb{M}^n_{k_0}}(v_j)}{v_j^{\frac{n-1}{n}}}-\frac{\eta_j}{v_j^{\frac{n-1}{n}}} & \le & \frac{\sum_lI_{\mathbb{M}^n_{k_0}}(v_{j,l})}{v_j^{\frac{n-1}{n}}}-\frac{\eta_j}{v_j^{\frac{n-1}{n}}}\\
& \le & \frac{\sum_l\P_g(\Omega_{j,l})}{v_j^{\frac{n-1}{n}}}-\frac{\eta_j}{v_j^{\frac{n-1}{n}}}\\ \label{Eq:BM3}
& \le & \frac{\P_g(\Omega_j)}{v_j^{\frac{n-1}{n}}},
\end{eqnarray}
where the first inequality is due to the strict subadditivity of $I_{\mathbb{M}^n_{k_0}}$, the second is due to Theorem \ref{Res:Theorem1} (because $\diam(\Omega_{j,l})<d$ for $j$ large enough), and the last inequality is due to Lemma \ref{Lemma:Partitioning2}. For all $j$ large enough we have that 
\begin{equation}\label{Eq:ResTheorem2Sequence}
(1-\varepsilon)\frac{I_{\mathbb{M}^n_{k_0}}(v_j)}{v_j^{\frac{n-1}{n}}}\le\frac{I_{\mathbb{M}^n_{k_0}}(v_j)}{v_j^{\frac{n-1}{n}}}-\frac{\eta_j}{v_j^{\frac{n-1}{n}}},
\end{equation}
thus $$(1-\varepsilon)\frac{I_{\mathbb{M}^n_{k_0}}(v_j)}{v_j^{\frac{n-1}{n}}}\le\frac{\P_g(\Omega_j)}{v_j^{\frac{n-1}{n}}}, $$
The last inequality combined with \eqref{Eq:BM}-\eqref{Eq:BM3} easily establish the validity of \eqref{Eq:ResTheorem2} and complete the proof of the theorem.
\end{proof}
\begin{Cor}\label{Thm:BestSobolevConstant} Let $M$ be a complete Riemannian manifold with strong bounded geometry. Then $$\liminf_{v\to0^+}\dfrac{I_M(v)}{v^{\frac{(n-1)}{n}}}=c_n,$$ where $c_n$ is the Euclidean isoperimetric constant defined by $I_{\R^n}(v)=c_nv^{\frac{n-1}{n}}$.
\end{Cor}
\begin{proof}
Take an arbitrary sequence $v_j\to 0$ and a sequence of positive real numbers $\varepsilon_j\to 0$, by the definition of $I_{M,g}$ we know that we can take a sequence of finite perimeter sets $\Omega_j$, such that $V_g(\Omega_j)=v_j$ and $I_{M,g}(v_j)\le\P_g(\Omega_j)\le I_{M,g}(v_j)+\varepsilon_j$. 
Passing to the limit in \eqref{Eq:ResTheorem2} or using \eqref{Eq:ResTheorem2Sequence} combined with the asymptotic expansion of the perimeter of geodesic balls in the model simply connected space forms in function of the volume enclosed it follows that 
\begin{equation}\label{Eq:BestSobolevConstant}
c_n\le\liminf_{j\to+\infty}\frac{\P_g(\Omega_j)}{v_j^{\frac{n-1}{n}}}.
\end{equation}
With this last inequality and Inequality \eqref{Eq:BestSobolevConstant} in mind the theorem follows without any further complications.
\end{proof}
Before to continue, let us state some results of independent interest that will be crucial in the proof of Theorem \ref{Res:Theorem3}.
\begin{Thm}[Selecting a large subdomain non effective]\label{Thm:SelectingaLarge}
Let $(M^n,g)$ be a complete Riemannian manifold with weak bounded geometry, and $(D_j)_j$ is an approximate isoperimetric sequence. Then there exists another approximate isoperimetric sequence $(D_{j}')$ such that $\lim_{j\to\infty} V_g(D_j\triangle D_j')=0$, $\lim_{j\to\infty} \frac{V_g(D_j')}{V_g(D_j)}=1$, $\lim_{j\to\infty}\dfrac{\P_g(D_j')}{\P_g(D_j)}=1$ and $diam (D_j')\to 0$, when $j\to \infty$.
\end{Thm}

\begin{proof}
We perform the same construction of a partition as in the proof of Lemma \ref{Lemma:Partitioning} applied to any $D_j$ with a suitable radius $R_j$ that we will choose later, and obtain a suitable partition $D_{j,l}$ of $D_j$ a maximal family of points $\mathcal{A}_j$ such that
\begin{equation}
\left(\sum_l\P_g(D_{j,l})\right)-\P_g(D_j)\le N(n,k,v_0)\dfrac{V_g(D_j)}{R_j}.
\end{equation}
Set $v_j:=V_g(D_j)$, by the definition of $\lambda$ and of $I_M$, we get that for large $j$ %every $0<\varepsilon<1$
 it holds
\begin{equation}
\P(D_{j,l})\ge I_{M,g}(V(D_{j,l}))\ge\lambda v_{j,l}^{(n-1)/n},
\end{equation}
where $v_{j,l}:=V_g(D_{j,l})$. Trivially for large $j$ we have $V_g(D_{j,l})\le v_j\le\bar{v}$ and the Euclidean type isoperimetric inequality for small volumes holds. This implies by Lemma \ref{Lemma:Partitioning} that
\begin{equation}
\dfrac{\sum_l\lambda V_g(D_{j,l})^{(n-1)/n}}{v_j^{(n-1)/n}}
\le\dfrac{\sum_l\P_g(D_{j,l})}{v_j^{(n-1)/n}}\le\frac{\P_g(D_j)}{v_j^{\frac{n-1}{n}}}+N\frac{v_j^{\frac{1}{n}}}{R_j}.
\end{equation}
Using the arguments of the combinatorial Lemma $2.3$ of \cite{NarCalcVar} applied to $f_{j,l}:=\dfrac{V_g(D_{j,l})}{v_j}$, we get that
$f^*_j:=\max\{f_{j,1},...,f_{j,l_j}\}$ satisfies
\begin{equation*}
\sum_lf_{j,l}{f^*_{j}}^{-\frac{1}{n}}\le\sum_lf_{j,l}f_{j,l}^{-\frac{1}{n}}=\sum_lf_{j,l}^{\frac{n-1}{n}}\leq\frac{1}{\lambda}\left[\frac{\P_g(D_j)}{v_j^{\frac{n-1}{n}}}+N\frac{v_j^{\frac{1}{n}}}{R_j}\right],
\end{equation*}
hence
\begin{equation}
f^*_j\ge\left[\frac{\lambda}{\frac{\P(D_j)}{v_j^{\frac{n-1}{n}}}+N\frac{v_j^{\frac{1}{n}}}{R_j}}\right]^n,
\end{equation}
which implies
\begin{equation}
V_g(D'_{j,1})\ge v_j\left[\frac{\lambda}{\frac{\P(D_j)}{v_j^{\frac{n-1}{n}}}+N\frac{v_j^{\frac{1}{n}}}{R_j}}\right]^n.
\end{equation}
On the other hand we recall that by construction there exists a point $p'_{j}\in\mathcal{A}_j$ depending on $D_j$ such that $D'_{j,1}\subseteq B_M(p_{D_j}, R_j)$. Fix an arbitrary sequence $\mu_j\to+\infty$ and set $R_j:=\mu_jv_j^{\frac{1}{n}}$, $D'_j:=B_M(p_{D_j}, R_j)\cap D_j$, $v'_j:=V(D'_j)$, $l_{1,j}:=\P_g(D_j, B_M(p_{D_j}, R_j))$, $l_{2,j}:=\P(D_j)-l_{1,j}$, $A_j:=\P_g(D'_j)=l_{1,j}+\frac{\Delta v_j}{R_j}$, and $\Delta v_j:=v_j-v'_j$ we have $D_{1,j}\subseteq D'_j\subseteq D_j$ thus
\begin{equation}
\lim_{j\to+\infty}\frac{v'_j}{v_j}=1,
\end{equation}
\begin{equation}
\frac{\Delta v_j}{v_j}\le1-\left[\frac{\lambda}{\frac{\P(D_j)}{v_j^{\frac{n-1}{n}}}+N\frac{v_j^{\frac{1}{n}}}{R_j}}\right]^n=1-\left[\frac{\lambda}{\frac{\P(D_j)}{v_j^{\frac{n-1}{n}}}+\frac{N}{\mu_j}}\right]^n,
\end{equation} 
\begin{equation}
\frac{\Delta v_j}{v_j}\to0,
\end{equation}
\begin{equation}
\frac{A_j}{\P_g(D_j)}\to1,
\end{equation}
\begin{equation}
\frac{l_{1,j}}{\P_g(D_j)}\to1,
\end{equation}
\begin{equation}
\frac{l_{2,j}}{v_j^{\frac{n-1}{n}}}\to0.
\end{equation}
\end{proof}
Essentially Theorem \ref{Thm:SelectingaLarge} says that for small volumes approximate isoperimetric sequences have all the mass and perimeter that stay inside a ball of small radius. What will be proved further in Lemma \ref{Lemma:smallvolumesimpliessmalldiametersmild} is that in fact the part outside this latter ball in fact does not give any contribution and actually for small volumes is empty.
\begin{Def} Let $(M^n,g)$ be a Riemannian manifold. We say that \textbf{$(M,g)$ satisfy $(H)$}, if there exists 
\begin{equation}
\lim_{v\to0^+}\frac{I_M(v)}{v^{\frac{n-1}n}}=\liminf_{v\to0^+}\frac{I_M(v)}{v^{\frac{n-1}n}}=\limsup_{v\to0^+}\frac{I_M(v)}{v^{\frac{n-1}n}}=\lambda.
\end{equation}  
\end{Def}
In the next theorem we will give a little more refined proof of Theorem \ref{Thm:SelectingaLarge} having the advantage of being effective.
\begin{Thm}[Selecting a large subdomain effective]\label{Thm:SelectingaLargeEffective}
Let $M^n$ be a complete Riemannian manifold with weak bounded geometry and $\mu>0$. Then there exists $\bar v_2=\bar v_2(n,k,v_0,\mu)>0$ such that for any finite perimeter set $D$ with volume $v\le\bar v_2$ there exists $p_D\in M$ another finite perimeter set $D':=B_M(p_D, \mu v^{\frac1n})\cap D\subseteq D$ such that 
\begin{equation}
V_g(D\triangle D')\le 1-\left[\frac{\lambda}{\frac{\P_g(D)}{v^{\frac{n-1}{n}}}+\frac{N}{\mu}}\right]^n,
\end{equation} 
\begin{equation} 
\frac{V_g(D')}{V_g(D)}\ge\left[\frac{\lambda}{\frac{\P_g(D)}{v^{\frac{n-1}{n}}}+\frac N\mu}\right]^n.
\end{equation} 
In particular $diam_g(D')\le2\mu v^{\frac1n}$.
\end{Thm}
\begin{proof}
First of all choose $\bar v_2\le\min\{\frac{\bar R}\mu, \bar v_1\}$. Then perform the same construction of a partition as in the proof of Lemma \ref{Lemma:Partitioning} applied to $D$ with radius $R:=\mu v^{\frac1n}$, and obtain a suitable partition $\{D_l\}_l$ containing a finite number $l_D$ of of components $D_j=\mathring{\cup}D_l$ joint with a maximal family of points $\mathcal{A}$ such that
\begin{equation}
\left(\sum_{l=1}^{l_D}\P_g(D_l)\right)-\P_g(D)\le N(n,k,v_0)\dfrac{V_g(D)}{R}.
\end{equation}
Set $v_l:=V_g(D_l)$, by the definition of $\lambda$ and of $I_M$,
 it holds
\begin{equation}
\P_g(D_l)\ge I_{M,g}(V_g(D_l))\ge\lambda v_l^{\frac{n-1}n}.
\end{equation}
Since $v\le\bar v_2$ we have $V_g(D_l)\le v\le\bar v_2$ and the Euclidean type isoperimetric inequality for small volumes holds. This implies by Lemma \ref{Lemma:Partitioning} that
\begin{equation}
\dfrac{\sum_l\lambda V_g(D_l)^{(n-1)/n}}{v^{(n-1)/n}}
\le\dfrac{\sum_l\P_g(D_l)}{v^{(n-1)/n}}\le\frac{\P_g(D)}{v^{\frac{n-1}{n}}}+N\frac{v^{\frac{1}{n}}}{R}.
\end{equation}
Using the arguments of the combinatorial Lemma $2.3$ of \cite{NarCalcVar} applied to $\gamma_l:=\dfrac{V_g(D_l)}{v}$, we get that
$\gamma^*:=\max\{\gamma_1,...,\gamma_{l_D}\}$ satisfies
\begin{equation*}
\sum_l\gamma_l{\gamma^*}^{-\frac{1}{n}}\le\sum_l\gamma_l\gamma_l^{-\frac{1}{n}}=\sum_l\gamma_l^{\frac{n-1}{n}}\leq\frac{1}{\lambda}\left[\frac{\P_g(D)}{v^{\frac{n-1}{n}}}+N\frac{v^{\frac{1}{n}}}{R}\right],
\end{equation*}
hence
\begin{equation}
\gamma^*\ge\left[\frac{\lambda}{\frac{\P_g(D)}{v^{\frac{n-1}{n}}}+N\frac{v^{\frac{1}{n}}}{R}}\right]^n,
\end{equation}
which implies
\begin{equation}\label{Eq:SelLargeEffective}
V(\tilde D'_1)\ge v\left[\frac{\lambda}{\frac{\P_g(D)}{v^{\frac{n-1}{n}}}+N\frac{v^{\frac{1}{n}}}{R}}\right]^n,
\end{equation}
where $\tilde D'_1$ is one of the connected components of the partition $\{D_l\}_l$ of $D$ that satisfy $\frac{V_g(\tilde D'_1)}v=\gamma^*$. On the other hand we recall that by construction there exists a point $p_D\in\mathcal{A}$ depending on $D$ such that $\tilde D'_1\subseteq B_M(p_D, R)$. Set $D':=B_M(p_D, R)\cap D$, $v':=V_g(D')$, %$l_1(D):=\P_g(D, B_M(p_D, R))$, %$l_2(D):=\P_g(D)-l_1(D)$, $A(D):=\P_g(D')=l_1(D)+\frac{\Delta v}{R}$, 
and $\Delta v:=v-v'$ we have $\tilde D'_1\subseteq D'\subseteq D$ thus by \eqref{Eq:SelLargeEffective}
\begin{equation}
\frac{v'(D)}{v}\ge\left[\frac{\lambda}{\frac{\P_g(D)}{v^{\frac{n-1}{n}}}+\frac{N}{\mu}}\right]^n,
\end{equation}
uniformly with respect to $D\in\tilde\tau_v$ ($\tilde\tau$ is defined in Definition ). Furthermore, it is also easily seen by \eqref{Eq:SelLargeEffective} that  
\begin{equation}
\frac{\Delta v}{v}\le1-\left[\frac{\lambda}{\frac{\P_g(D)}{v^{\frac{n-1}{n}}}+\frac{N}{\mu}}\right]^n.
\end{equation} 
\begin{Rem} At this stage we made the choice of not controlling the perimeter added cutting with a ball of radius $R$ by a coarea formula argument. We recall that this is always possible (by coarea formula) up to take a slightly larger radius $R+\eta_DR$ for a suitable $0<\eta_D<1$. 
\end{Rem}
\end{proof}
The following lemma have its own interest. Its proof is based on the adaptation of the arguments of the Deformation Lemma $4.5$ of \cite{RitoreGalli} and of formula (1.10) of \cite{Tamanini} also named Almgren's Lemma in some literature see for instance  the book \cite{Maggi}.
\begin{Def}\label{Def:Trace} For every $f\in BV(M)$ every $D\subset M$ with Lipschitz continuous boundary $\partial D$ we define
the \textbf{trace of $f$ on the boundary of $D$} as the function $f_{|\partial D}:\Sigma\subseteq\partial D\to\R$, where $\mathcal{H}^{n-1}(\partial D\setminus\Sigma)=0$ and defined by $$f_{|\partial B}:x\mapsto\lim_{r\to0^+}\frac1{B_M(x,r)}\int_{B_M(x,r)} f(y)dV_g(y).$$ 
\end{Def}
\begin{Lemme}[Theorem $2.10$ of \cite{Giusti}]\label{Lemma:WeakGaussGreen} 
Let $\Omega$ be a bounded open set in $M^n$ with Lipschitz continuous boundary $\partial \Omega$ and let $f\in BV(\Omega)$. Then there exists a function $\phi\in L^1(\partial \Omega)$ such that for $\mathcal H^{n-1}$-almost all $x\in\partial\Omega$ 
\begin{equation}
\lim_{\rho\to0^+} \rho^{-n}\int_{B_M(x, \rho)\cap\Omega}|f(z)-\phi(x)|dz=0.
\end{equation}  
In particular such a $\phi=f_{|\partial D}$ belongs to $L^1(\partial \Omega)$.
Moreover, for every $X\in\mathfrak{X}^1_0(M)$,
\begin{equation}
\int_\Omega f div Xdx=-\int_\Omega\langle X,Df\rangle+\int_{\partial\Omega}\phi\langle X,\nu\rangle,
\end{equation}
where $\nu$ is the unit outer normal to $\partial\Omega$.
\end{Lemme}
\begin{Lemme} [Deformation Lemma]\label{Lemma:Deformation} Let $(M^n,g)$ be a Riemannian manifold, $p\in M$, $B:=B_M(p,r)$ a geodesic ball with $0<r<inj_M$, $k\in\R$, $(n-1)k$ a lower bound on the Ricci curvature tensor inside $B$, $u:B\to[0,+\infty[$, $u:x\mapsto d_M(p,x)$, is the distance function to the point $p$, $E\subseteq B$ a set of locally finite perimeter in $B$. Then it holds 
\begin{equation}\label{Eq:DeformationStatement}
\P_g(\partial B, \left(M\setminus E\right)^{(1)})\le\P_g(E, B)+\frac crV_g(B\setminus E),
\end{equation}
where $c=c(n,k, inj_M):=1+(n-1)c_k(inj_M)>0$ is a positive constant, and $\left(M\setminus E\right)^{(1)}$ are the points of density $1$ of $\left(M\setminus E\right)$.
\end{Lemme}
\begin{proof} 
Applying Lemma \ref{Lemma:WeakGaussGreen} with $f=\chi_{E^c}$ and $X:=\varphi\frac ur\nabla u$, where $E^c:=B\setminus E$ and $\varphi\in C^1_0(B_M(p,inj_M))$ with the property that $\{x\in M :\:d_M(x, B)\le\varepsilon\}\subseteq\varphi^{-1}(1)\subseteq B_M(p,inj_M)$ for some small $\varepsilon>0$ (observe that this choice of $\varphi$ yields $X\in\mathfrak{X}^1_0(M))$ leads to  
\begin{eqnarray*}
\int_B\chi_{E^c}div_g\left(\frac ur\nabla u \right)dV_g & = & \int_{\partial B}\chi_{E^c}|_{\partial B}\frac ur\langle \nabla u,\nu_{ext}\rangle d\mathcal{H}^{n-1}\\
 & - & \int_B\left\langle \nabla \chi_{E^c},\frac {u\nabla u}r\right\rangle dV_g,
\end{eqnarray*}
then
\begin{align*}
\int_B\chi_{E^c}div_g\left(\frac ur\nabla u \right)dV_g &=\int_B\chi_{E^c}\left(\frac{\|\nabla u\|^2}{r}+\frac ur div_g(\nabla_gu)\right)dV_g\\
&\le \int_B\chi_{E^c}\left(\frac1r+\frac{u(x)(n-1)}r\cot_k(u(x))\right)dV_g(x)\\
&\le\int_B\chi_{E^c}\left(\frac1r+\frac{n-1}r c_k(r)\right)dV_g(x)\\
&\le\frac1r\left(1+(n-1)c_k(r)\right)V_g(E^c).
\end{align*}
On the other hand
\begin{align*}
-\int_B \left\langle\nabla \chi_{E^c},\frac ur \nabla u\right\rangle &=-\P(E,B).
\end{align*}
Hence
\begin{equation*}
\P(B\setminus E, (M\setminus E)^{(1)})-\P(E,B)\le \frac1r\left(1+(n-1)c_k(r)\right)V_g(E^c).
\end{equation*}
From the last inequality it is easy to deduce \eqref{Eq:DeformationStatement}, after the simple observation that $c_k(r):=r\cot_k(r)$ is a strictly increasing function in particular is bounded in $[0,inj_M]$.
\end{proof}
\begin{Rem} It is worth to recall here that by Theorem $1$ of \cite{Tamanini} (which immediately could be adapted to the Riemannian manifold because is a local theorem) an isoperimetric region have always nonempty interior as well as its complement but a lot of proofs of regularity do not give a satisfying and uniform estimates of the radius of the balls contained inside.
\end{Rem}
The following lemma have its own interest. Its proof is based on the adaptation of the arguments of the Deformation Lemma $4.5$ of \cite{RitoreGalli} which in this context are given by our Lemma \ref{Lemma:Deformation} combined with the arguments of Section $2$ of \cite{NarCalcVar} that are adapted here in Theorem $3.2$ with the use of the Heintze-Karcher comparison Theorem of \cite{HeintzeKarcher} combined with the proof of Lemma $3.8$ of \cite{NarBBMS} and Theorem $3$ of \cite{NarAsian}. 
\begin{Lemme}\label{Lemma:smallvolumesimpliessmalldiametersmild}
Let $(M^n,g)$ be a complete Riemannian manifold with mild bounded geometry satisfying $(H)$. Then there exist two positive constants $\mu^*=\mu^*(n,k,v_0,\lambda)>0$ and $v^*=v^*(n,k,v_0,\lambda)>0$ such that whenever $\Omega\subseteq M$ is an isoperimetric region of volume $0\le v\le v^*$ it holds that $$\diam_g(\Omega)\le\mu^*v^{\frac{1}{n}}.$$  
\end{Lemme}
\begin{Rem} In mild bounded geometry $v_0$ depends on $k$ and $inj_M$ so in the preceding lemma we have that $\mu^*=\mu^*(n,k,inj_M,\lambda)$ and $v^*=v^*(n,k,inj_M,\lambda)$. In strong bounded geometry condition $(H)$ is always fulfilled, moreover it is known that $\lambda=c_n$, hence in the preceding lemma when specialised to the case of strong bounded geometry we have actually $\mu^*=\mu^*(n,k,inj_M)$ and $v^*=v^*(n,k,inj_M)$. The construction made to prove the preceding lemma it is possible only because we assume positive injectivity radius. So the injectivity radius is hidden inside $\mu^*$ and $v^*$ although it is tempting to prove Lemma \ref{Lemma:smallvolumesimpliessmalldiametersmild} just assuming $M$ with weak bounded geometry instead of mild bounded geometry.  
\end{Rem}
\begin{Rem} As already observed, if $M$ have strong bounded geometry, then always exists $\lim_{v\to0^+}\frac{I_M(v)}{v^{\frac{n-1}n}}=\lambda$ and so in particular Lemma \ref{Lemma:smallvolumesimpliessmalldiametersmild} applies to manifold with strong bounded geometry. Unfortunately, we still do not know wether the existence of $\lim_{v\to0^+}\frac{I_M(v)}{v^{\frac{n-1}n}}=\lambda$ could be dropped or not in the statement of the preceding lemma. Obviously in weak bounded geometry or in mild bounded geometry, if one proves that $\lambda=c_n$, then automatically condition $(H)$ is fulfilled.
\end{Rem}
In fact the following questions are still open at the present stage of our knowledge. 
\begin{Question} If $(M^n,g)$ is with weak bounded geometry, then $M$ satisfy $(H)$? 
\end{Question}
\begin{Question} If $(M^n,g)$ is with mild bounded geometry, then $M$ satisfy $(H)$? 
\end{Question}
\begin{Question} If $(M^n,g)$ is with weak bounded geometry or with mild bounded geometry, what is the sharp value of $\lambda$? \end{Question}
%\begin{Rem} Via smoothing a lot of particular cases could be settled, for example if $(M,g)$ is a complete Riemannian manifold such that $|Ric_g|\le K$ and the conjugate radius $conj_M>0$, by a smoothing theorem as for example Theorem $1.1$ of \cite{DaiWeiYe}, we can say that $(M,g)$ satisfy Property $(H)$. In general if a property on the original manifold $(M,g)$ depend $C^0$ on the metric and the bounds on the curvatures are strong enough in general smoothing could permit to transport the same property when valid on the smoothed manifold to the original one.
%\end{Rem}
\begin{Rem} The main reason to assume positive injectivity radius in the preceding lemma is that we make a crucial use of Lemma \ref{Lemma:Deformation} which in turn uses radial deformations which are well defined only locally at a point $x\in M$ inside a ball of radius less than $inj_x$. We will see later in the proof of Lemma \ref{Lemma:smallvolumesimpliessmalldiametersmild} that we want to apply radial deformations with center at the point $p^*_{\Omega}$ defined further, but if $inj_M=0$ we have no control about the size of $inj_{p^*_{\Omega}}$ $($remember that $p^*_{\Omega}$ could go to infinity$)$ and the volume that we can put inside $B_M(p^*_{\Omega},inj_{p^*_{\Omega}})$. To avoid this problem of course it is enough to assume positive injectivity radius, but we still do not know whether this assumption could be dropped and replaced just by the noncollapsing of the volume of balls of radius $1$. 
\end{Rem}
\begin{Rem} It is well known by that in mild bounded geometry $v_0=v_0(n,k, inj_M)>0$. Thus in the statement of Lemma \ref{Lemma:smallvolumesimpliessmalldiametersmild} we can suppress the dependence of $v^*$ on $v_0$. 
\end{Rem}
The geometric idea of the proof is not too complicated but unfortunately the writing turns out to be technical, because of the effective calculations of the constants. In first by an application of Theorem \ref{Thm:SelectingaLarge} we find a point $p_\Omega\in M$ and a controlled radius $\mu v^{\frac1n}$, such that almost all the volume of $\Omega$ is recovered inside the ball $B_M(p_\Omega,\mu v^{\frac1n})$. In second, we take a ball inside $\Omega$ of controlled volume and radius and show that these two balls cannot be disjointed when the volume tends to $0$, so we take a bigger but still controlled radius. Then we proceed by contradiction and suppose that there are points of $\Omega$ very far from $p_\Omega$. Under this assumption we take a possibly bigger controlled radius such that we can pick a fraction of the volume of $\Omega$ that is far from $p_\Omega$ and compensate close to $p_\Omega$ to readjust volume. In this way a competitor $F$ of the same volume than $\Omega$ is constructed provided the volume of $\Omega$ is taken small enough. Finally it is shown that this last competitor $F$ have perimeter strictly less than the perimeter of $\Omega$ for volumes possibly smaller. This gives a contradiction and so the impossibility to find points of $\Omega$ far from $p_\Omega$.
\begin{proof}[Proof of Lemma \ref{Lemma:smallvolumesimpliessmalldiametersmild}] For simplicity of notations we consider just the case $k\le 0$. When $k>0$ the theorem of Bonnet-Myers ensures that $M$ is compact and so the lemma is already proved in Theorem $2.2$ of \cite{MJ}. Our proof works also without the restriction $k\le 0$. For now on in this proof we assume that $k\le0$.  The hypothesis $(H)$ permit to us to define the quantities
 \begin{equation}\label{Eq:Structural0}
 \tilde f^*(n, k,v_0,\lambda,\mu):=\lim_{v\to 0^+}\tilde f(v,n,k,v_0,\lambda,\mu)=\left[\frac{\lambda}{\lambda+\frac{N}{\mu}}\right]^n,
 \end{equation}
 \begin{equation}\label{Eq:Structural01}
 f^*(n,k,v_0,\lambda,\mu):=\lim_{v\to 0^+}f(v,n,k,v_0,\lambda,\mu)=1-\left[\frac{\lambda}{\lambda+\frac{N}{\mu}}\right]^n,
 \end{equation}
where $$\tilde f(v,n,k,v_0,\lambda,\mu, M):=\left[\frac{\lambda}{\frac{I_M(v)}{v^{\frac{n-1}{n}}}+\frac{N}{\mu}}\right]^n,$$ $$f(v,n,k,v_0,\lambda,\mu, M):=1-\left[\frac{\lambda}{\frac{I_M(v)}{v^{\frac{n-1}{n}}}+\frac{N}{\mu}}\right]^n.$$
In the remaining part of this proof we will use frequently the two following crucial properties  
\begin{equation}\label{Eq:Structural02}
\lim_{\mu\to+\infty}\tilde f^*(n, k,v_0,\lambda,\mu)=1,
\end{equation}
\begin{equation}\label{Eq:Structural03}
\lim_{\mu\to+\infty}\tilde f^*(n, k,v_0,\lambda,\mu)=0.
\end{equation}
Suppose until the end of the proof that $\Omega$ is an isoperimetric region of volume $v$. By Heintze-Karcher's theorem we have that in weak bounded geometry $$c_7(n,k)v^{\frac1n}\ge inrad(\Omega)\ge\frac{v}{I_M(v)}\ge \frac{v}{I_{\mathbb{M}_k^n}(v)}\ge\frac{v}{c_2v^{(n-1)/n}}=\tilde c_2 v^{\frac1n},$$ for some positives constants $c_7=c_7(n,k)>0$ and $\tilde c_2=\tilde c_2(n,k)>0$. To see this in details the reader could consult Lemma $3.1$ of \cite{NarBBMS}.
In first we observe that the hypothesis that $M$ satisfy $(H)$ permits to have \eqref{Eq:Structural0}-\eqref{Eq:Structural03} which in turn allow us to choose $\mu=\mu(n,k,v_0,\lambda)>0$ large enough to satisfy simultaneously
\begin{equation}\label{Eq:Structural1}
\mu>c_1(n,k)^{1/n},
\end{equation} 
\begin{equation}\label{Eq:Structural2}
\mu>c_7(n,k),
\end{equation} 
\begin{equation}\label{Eq:Structural3}
\dfrac{\tilde c_2}{2n}\left(2f^*(n,k,v_0,\lambda,\mu)\right)^{1/n}<\dfrac{C_{Heb}}{4n},
\end{equation}
where $C_{Heb}>0$ is the constant appearing in Lemma $3.2$ of \cite{Hebey} %and $C_4=C_4(n,k,v_0)>0$ is such that $ $
.As it is easy to see, using \eqref{Eq:Structural0} and \eqref{Eq:Structural01} we can prove the existence of $v^*=v^*(n,k,v_0, inj_M, \lambda)>0$ such that for every $v\leq v^*$ we have that the following conditions are satisfied
\begin{equation}\label{Eq:Volumic0}
v^{\frac{n-1}n}f(v,n,k,v_0,\lambda,\mu)<c_1(n,k)v_0,
\end{equation}
\begin{equation}\label{Eq:Volumic1}
f(v,n,k,v_0,\lambda,\mu)\leq 2f^*(n,k,v_0,\lambda,\mu),
\end{equation}
%\begin{equation}\label{Eq:Volumic2}
%\tilde f(v,n,k,v_0,\lambda,\mu)\geq\frac34\tilde f^*(n, k,v_0,\lambda,\mu),
%\end{equation}
\begin{equation}\label{Eq:Volumic3}
r_v=4\mu v^{\frac1n}\le\frac14inj_M,
\end{equation}
\begin{equation}\label{Eq:Volumic4-}
c_1(n, k, \mu v^{\frac1n})=\tilde C_1(n,k)\mu^nv,
\end{equation}
\begin{equation}\label{Eq:Volumic4}
v\le\min\{1,\overline v, \bar v_1,\bar v_2\},
\end{equation}
where $\overline v$ is obtained in Lemma $3.2$ of \cite{Hebey}, i.e., such that for volumes smaller than $\overline v$ it holds $I_M(v)\ge\lambda v^{\frac{n-1}{n}}$. In the remaining part of this proof we always assume that $v\le v^*$. Consider an isoperimetric region $\Omega$ of $V(\Omega)=v$, the same construction of Theorem \ref{Thm:SelectingaLarge} applied to $\Omega$ gives the existence of $p_\Omega\in M$, (notice that the point $p_{\Omega}$ could be chosen satisfying the condition $p_\Omega\in\mathring{\Omega}$, but this is not relevant for the rest of our discussion) such that
\begin{eqnarray*}
\frac{V(B_M(p_\Omega,\mu v^{1/n})\cap\Omega)}{v} & = & \frac{v_1(\Omega)}{v}\\ 
& \ge & \left[\frac{\lambda}{\frac{\P_g(\Omega)}{v^{\frac{n-1}{n}}}+\frac{N}{\mu}}\right]^n\\
& \ge & \tilde f(v,n,k,v_0,\lambda,\mu, M)=\left[\frac{\lambda}{\frac{I_M(v)}{v^{\frac{n-1}{n}}}+\frac{N}{\mu}}\right]^n.
\end{eqnarray*}
Consider $\Delta v=\Delta v(\Omega):=v-v_1(\Omega)$, observe that
\begin{equation}\label{eq:lemmasmallvolsmalldiamst}
\frac{\Delta v}{v}=\frac{v-v_1(\Omega)}{v}\le 1-\left[\frac{\lambda}{\frac{\P_g(\Omega)}{v^{\frac{n-1}{n}}}+\frac{N}{\mu}}\right]^n\le1-\left[\frac{\lambda}{\frac{I_M(v)}{v^{\frac{n-1}{n}}}+\frac{N}{\mu}}\right]^n.
\end{equation}
Observe that we can put inside $\Omega$ a geodesic ball $$B_1(\Omega):=B_M(p^*_\Omega, inrad(\Omega))\subset \Omega.$$ 
We now show that $B_0(\Omega):=B_M(p_\Omega,\mu v^{1/n})$, cannot be disjoint from $B_1(\Omega)$. We prove this last assertion by contradiction. Indeed if it was the case we would have $B_1(\Omega)\subset \Omega\setminus B_0(\Omega)$, this would implies that $V(B_1(\Omega))\le V(\Omega\setminus B_0(\Omega))=\Delta v$, and in turn by estimative \eqref{eq:lemmasmallvolsmalldiamst}
$$c_1v_0v^{\frac1n}\le v(B_1(\Omega))\le vf_M(v,\mu)\le vf(v,n,k,v_0,\lambda,\mu),$$
which manifestly contradicts \eqref{Eq:Volumic0}. Hence we necessarily have $B_1(\Omega)\cap B_0(\Omega)\neq\emptyset$. Thanks to our choice \eqref{Eq:Structural2} the ball $B_1(\Omega)\subseteq B_M(p^*_\Omega, \mu v^{\frac1n})$, since $\mu v^{\frac1n}>c_7(n,k) v^{\frac1n}\ge inrad(\Omega)$. Moreover by our choice \eqref{Eq:Structural1} the $V(B_M(p^*_\Omega, \mu v^{\frac1n}))>v$, hence there exists a radius $r^*_v<\mu v^{\frac1n}$ such that the ball $B:=B_M(p^*_{\Omega},r^*_v)\subset B_2(\Omega):=B_M(p_\Omega, 3\mu v^{\frac1n})$ have $V(B)=v$. Notice that $B$ is just contained in $B_2(\Omega)$ and cannot be chosen as a proper subset of $\Omega$. This guarantees that $V(B\setminus\Omega)>0$ and furthermore that 
\begin{equation}\label{Eq:smallvolumeimpliessmalldiameters}
V(B\setminus\Omega)=V(B)-V(B\cap\Omega)\ge v-V(B_2\cap\Omega)\ge\Delta v,
\end{equation}
 because $V(B\setminus\Omega)=v-V(B_M(p^*_{\Omega},r_v^*)\cap\Omega)$ but $V(B_M(p^*_{\Omega},r^*)\cap\Omega)\le v_1(\Omega)$ and \eqref{Eq:smallvolumeimpliessmalldiameters} follows readily. Observe that by our choice \eqref{Eq:Volumic3} we have $B_M(p_{\Omega},r_v)\subset B_M(p_\Omega,\frac14inj_M)$. The following picture illustrates well our construction
\begin{figure}[H]
\label{fig:figura1}
\centering
\begin{tikzpicture}[y=0.50pt, x=0.50pt, yscale=-1.000000, xscale=1.000000]
\def\curve1{
    (297.3723,313.6493) .. controls (238.1343,338.7156) and (124.4206,222.4380)  .. %node {$N$}
    (98.1785,269.9494) .. controls (84.0761,295.4817) and (230.7401,315.0756) .. %node {$A$}
    (238.1343,338.7156) .. controls (252.1136,383.4088) and (66.0569,393.3281) ..%node {$B$}
    (88.4358,457.0431) .. controls (99.6898,489.0844) and (270.7532,369.6434) ..%node {$C$}
    (270.7532,369.6434) .. controls (328.6205,367.3904) and (187.5224,496.2889) ..%node {$D$}
    (216.8915,519.8158) .. controls (254.0874,549.6125) and (280.3138,428.8549) ..%node {$E$}
    (321.0530,415.7712) .. controls (358.2983,397.6895) and (475.5585,557.4375) .. %node {$F$}
    (525.7058,517.8497) .. controls (536.6322,492.3513) and (379.0718,403.3563) .. %node {$G$}
    (435.0861,403.3563) .. controls (495,400) and (520,440) .. %node {$H$}
    (565.6241,405.5590) .. controls (613.4999,373.5702) and (457.6,363.1860) ..%node {$I$}
    (443.1891,332.5355) .. controls (418.4455,279.7971) and (605.5903,252.7507) ..% node {$J$}
    (650,200) .. controls (600,200) and (420,325) ..  %node {$K$}
    (377.5179,306.2222) .. controls (322.2413,286.5862) and (262.3497,152.7256) ..  %node {$L$}
    (233.7785,212.2701) .. controls (222.3731,236.0399) and (369.2765,288.1066) .. %node {$M$}
    (297.3723,313.6493) -- cycle;
}
\coordinate (A) at (340,357.5); %% p*
\coordinate (B) at (310,339); %% p
\def\inrad{57};
\def\rav{82}; %rav=r_v^*
\def\ramu{100}; %\mu v^1/n
\def\ramub{1.8*\ramu}; %3\mu v^1/n
\def\rainj{290}; %3\mu v^1/n
\def\rab{125};
\def\geodballa{ (A) circle (\rainj)} ; %B(R2,p) \tilde B
\def\balla{ (A) circle (\inrad)}:            %B(p*,inrad) B_1
\def\ballb{ (A) circle (\rav)};            %B(p*, r_v^*) 
\def\ballc{ (A) circle (\ramu)};         %B(p*, \mu v^1/n) 
\def\balld{ (B) circle (\ramub)};       %B(p, 3\mu v^1/n) B_2
\def\balle{ (B) circle (2.3*\ramu)};   %B(p, r) B_r
\def\ballf{ (A) circle (\rab)};             %B(p, rb) B
\begin{scope}
\path [pattern=north east lines, pattern color=blue] \curve1;
     \fill [white] \balle;
\end{scope}
\begin{scope}
%\clip \balle;
\path [pattern=north west lines, pattern color=red] \curve1;
\fill [white] \balld;
\end{scope}    
\begin{scope}
\path [ pattern=north east lines, pattern color=blue] \ballb;
%\path [ pattern=north west lines, pattern color=red] \ballb;
      \clip \curve1;
      \fill [white] \ballb;
\end{scope}    
\path [draw=black,line join=miter,line cap=butt,even odd rule,line width=1.5pt] \curve1;    
\path [draw] \balla;    
\path [draw] \ballb;
\path [draw , line width=0.8pt] \balld;
\path [draw, line width=0.8pt] \balle;
\path [draw] \ballf;
\path [draw, line width=1pt] \geodballa;
\draw[fill]  (B) circle (2) node [right]   {$p_\Omega$}; 
\draw[fill]  (A) circle (2) node [right]   {$p^*_\Omega$}; 
\draw[line width=0.7pt, <-]  ($(A)+(0:\inrad)$)--($(A)+(0:310)$) node [right] {$B_1$};
\draw[line width=0.7pt, <-]  ($(A)+(-15:\rav)$)--($(A)+(-15:310)$) node [right] {$B_{\rho_1}$};
\draw[line width=0.7pt, <-]  ($(B)+(-60:\ramub)$)--($(B)+(-60:245)$) node [above right] {$B_2$};
\draw[line width=0.7pt, <-]  ($(A)+(-45:\rainj)$)--($(A)+(-45:310)$) node [above] {$\tilde B$};
\draw[line width=0.7pt, <-]  ($(B)+(60:2.3*\ramu)$)--($(B)+(55:270)$) node [right] {$B_r$};
\draw[line width=0.7pt, <-]  ($(A)+(-90:\rab)$)--($(A)+(-90:255)$) node [above] {$B$};
\draw[line width=0.7pt, dashed]  ($(B)$)--($(B)+(60:2.3*\ramu)$); %%radii;
\fill [white] ($(B)+(60:2.1*\ramu)$) circle (8.0);
\node at  ($(B)+(60:2.1*\ramu)$) {$r$} ; %%radii;
\draw[line width=0.7pt, <-]  (282,300)--($(282,300)+(-150:110)$) node [above] {$F$};
\draw[line width=0.7pt, <-]  (230,510)--($(230,510)+(90:90)$) node [right] {$\Delta v$};
\draw[line width=0.7pt, <-]  (620,218) ..controls (625, 170) .. (620,120) node[above] {$\Delta^*v$};
\draw[line width=0.7pt, <-]  (330,430)--($(330,430)+(80:150)$) node[below] {$\Delta^*v$};
\node at (530,490) {$\Omega$};
\end{tikzpicture}
\caption{Construction of the competitor $F:=(B_3\cup \Omega)\setminus B_r$ used in the proof of Lemma \ref{Lemma:smallvolumesimpliessmalldiametersmild}. Here $\tilde B:=B_M(p^*_{\Omega}, inj_M)$, $B_2:=B_M(p_\Omega, 3\mu v^{\frac1n})$, $B_r:=B_M(p_{\Omega}, r)$.} 
\end{figure}
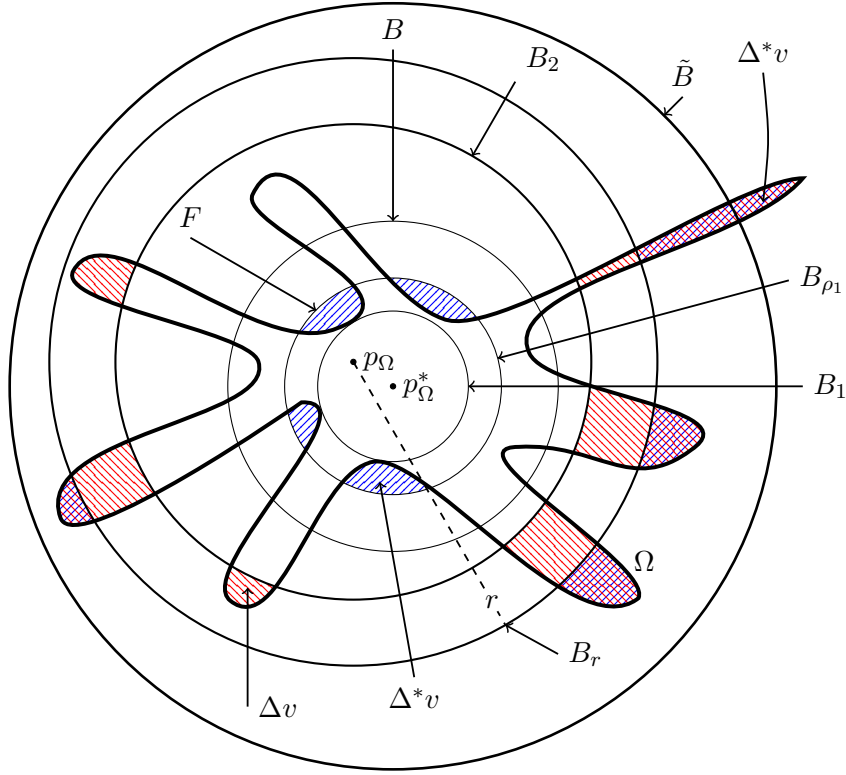
Assume the following notations
$$
\tilde d_\Omega:=\sup_{x\in\Omega}\{d(x,p_\Omega)\}, \quad d_\Omega=\tilde d_\Omega-r_v, \quad d_v:=\sup_{\Omega\in\tilde\tau, V(\Omega)=v}\{d_\Omega\}.
$$
For any $r>0$ let us define $V_\Omega(r):= V(\Omega \cap (M\setminus \bar B_r))=V(U_r)$ where $B_r := \{y \in M : d_M(p_{\Omega},y)<r\}$, $p_{\Omega}$ is given by Theorem \ref{Thm:SelectingaLarge}, and 
$U_r=\Omega \cap (M\setminus\bar{B}_r)$. The function $V_{\Omega}(r)$ is monotone decreasing and $V_{\Omega}(r)\searrow0$ as $r\to\infty$. Denote by $A_\Omega(r):=A(\partial \Omega\cap (M\setminus\bar{B}_r))$. Coarea formula gives immediately
$$
V(\Omega \cap (M\setminus\bar{B}_r))=\int_r^\infty A(\Omega\cap\partial B_r)dr,
$$
then
$$
V_\Omega'(r)=-A(\Omega\cap \partial B_r)=-A(\Omega\cap \partial (M\setminus B_r)).
$$
Consider any $r\ge3\mu v^{\frac1n}$ and put all the volume $\Delta^*v:=V_\Omega(r)$ inside $B$, by choosing a concentric ball $B_3$ with $B_1\subset B_3\subset B\subset B_2$ of radius $$Cv^{\frac{1}n}\le\rho_1=\rho_1(v,r)\le\mu v^{\frac{1}n},$$ such that $V(B_3\setminus\Omega)=\Delta^*v$, then $$F:=(B_3\cup \Omega)\setminus B_M(p_\Omega, r)=(B_M(p^*_\Omega,\rho_1(v,r))\cup\Omega)\setminus B_M(p_\Omega,r),$$ satisfy $V(F)=V(\Omega)$. Now, we make an application of Lemma \ref{Lemma:Deformation} with $E=\Omega$ inside the ball $B$. From this and the fact that $\Omega$ is an isoperimetric region
follows that $\P(\Omega)\le\P(F)$ and thus for almost all $r\ge3\mu v^{\frac1n}$ it holds   
\begin{equation}\label{Eq:smalldiam0}
l_1(\Omega)(r)+A_\Omega(r)\le l_1(\Omega)(r)-V_\Omega'(r)+(1+(n-1)c_k(\rho_1))\dfrac1{\rho_1}V_\Omega(r),
\end{equation}
where $l_1(\Omega)(r):=\P(\Omega, B_M(p_\Omega, r))$. This easily leads to
\begin{equation}\label{Eq:smalldiam}
A_\Omega(r)\le-V_\Omega'(r)+KV_\Omega(r),
\end{equation}
where $K=K(n,k,v_0,v, \mu)=\tilde c_2(n,k)\frac{1}{v^{\frac{1}{n}}}$.
Independently, by the Euclidean type isoperimetric inequality for small volumes of Lemma 3.2 of \cite{Hebey} we have that for small volumes there exists a positive constant $\bar{v}_1=\bar{v}_1(n,k,v_0)>0$, such that if $v\le\bar{v}_1$, then for every $r>0$ it holds
$$
C_{Heb}V_\Omega(r)^{(n-1)/n}\le A (\partial U_r),
$$ where $C_{Heb}=C_{Heb}(n,k,v_0)>0$ is given by Lemma $3.2$ of \cite{Hebey} too. Thus for almost every $r>0$ we have the following
\begin{align}\label{Eq:smalldiam0}
-V_\Omega'(r)+A_\Omega(r)
=&A(\Omega\cap \partial (M\setminus B_r))+A(\partial \Omega\cap (M\setminus\bar B_r))\nonumber
\\ \nonumber 
\ge& A (\partial ( \Omega\cap (M\setminus B_r)))\\ 
=& A (\partial U)\ge C_{Heb} V_\Omega(r)^{(n-1)/n}.
\end{align}
Adding the two inequalities \eqref{Eq:smalldiam} and \eqref{Eq:smalldiam0} we get that
$$
2V_\Omega'(r)\le KV_\Omega(r)-C_{Heb}V_\Omega(r)^{(n-1)/n}.
$$
Using the fact that $n(V_\Omega^{1/n})'=V_\Omega^{\frac{1}{n}-1}V_\Omega'$ we can write the preceding inequality as
$$
(V_\Omega^{1/n})'(r)\le\dfrac{\tilde c_2}{2n}\left(\frac{V_\Omega(r)}{v}\right)^{1/n}-\dfrac{C_{Heb}}{2n},
$$
for every $v\le v^*$ and $\Omega$ such that $V(\Omega)=v$, where $C_{Heb}=C_{Heb}(n,k,v_0)>0$ is the constant appearing in the isoperimetric inequality for small volumes of Lemma $3.2$ of \cite{Hebey} reported here as Lemma \ref{Lemma:Hebey3.2}. Since $r\ge3\mu v^{\frac1n}$, one have 
$$(V_\Omega^{1/n})'\le\dfrac{\tilde c_2}{2n}\left(2f^*(n,k,v_0,\lambda,\mu)\right)^{1/n}-\dfrac{C_{Heb}}{2n}.$$
By Theorem \ref{Thm:SelectingaLarge} and \eqref{Eq:Structural3} we argue that 
\begin{equation}\label{Eq:DifferentialInequality}
(V_\Omega^{1/n}(r))'\le -C'=-\dfrac{C_{Heb}}{4n}.
\end{equation}
It is worth to recall here that by Theorem $3$ of \cite{NarAsian}, in weak bounded geometry $diam(\Omega)<+\infty$, because $\Omega$ is an isoperimetric region, and hence $d_{\Omega}:=essup_{x\in\Omega}d_M(p_{\Omega},x)=||d_M(p_{\Omega},\cdot)||_{L^{\infty}(\Omega)}<+\infty$. Furthermore we have the elementary relation $diam(\Omega)\le 2d_{\Omega}$. Now, if we assume $r_v:=3\mu v^{\frac1n}<d_{\Omega}$, we can integrate \eqref{Eq:DifferentialInequality} over the interval $[r_v, d_{\Omega}]$, and noting that $V_\Omega(r_v)\le V(\Omega)=v$, $V_\Omega(d_{\Omega})=0$, we get
$$
d_{\Omega}\le\dfrac1{C'}V_{\Omega}(r_v)^{1/n}+r_v\le\frac{v^{\frac{1}{n}}}{C'}+r_v=\left(\frac1{C'}+3\mu\right)v^{\frac{1}{n}}.
$$
From this last equation we easily a constant $\mu^*$ such that for every $v\le v^*$ results $$diam_g(\Omega)\le\mu^*(n,k,v_0,\lambda)v^{\frac{1}{n}},$$ which clearly proves the lemma.
\end{proof}
\section{Isoperimetric comparison in strong bounded geometry}
Now we are in position to prove Theorem \ref{Res:Theorem3}.
\begin{proof}[Proof of Theorem \ref{Res:Theorem3}] If $|Sec_M|\le K$ and $inj_M>0$ then the assumptions of Theorem $76$ of 
\cite{Petersen} holds which implies that also the assumptions of Theorem $72$ of \cite{Petersen} are satisfied with $m=1$, see 
also Theorem $4.4$ of \cite{Peters}. The problem here is that the limit metric space have an atlas of harmonic coordinates of class 
$C^{3,\alpha}$ with just a $C^{1,\alpha}$ limit metric. Unfortunately, to apply Theorem \ref{Res:Theorem1} to a limit manifold we 
need to have in the limit a smooth Riemannian manifold with a smooth Riemannian metric, for this reason we make a stronger 
assumption on $(M^n,g)$ requiring that $M$ have strong bounded geometry smooth at infinity.  A fortiori $M$ have also $C^0$-locally 
bounded geometry. This means by Theorem $1$ of \cite{NarAsian} or Theorem $1$ of \cite{FloresNardulliGenComp} that for every 
$v\in]0, V(M)[$ there exists a generalized isoperimetric region $\tilde\Omega_v$ contained in some smooth limit manifold $
(M_{\infty}, g_{\infty})$ (that could even coincide with $(M,g)$). Now if we look at the limit manifold $(M,g_{\infty})$, by $(iv)$ in
 Theorem $4.4$ of \cite{Peters} we learn that $inj_{(M_{\infty},g_{\infty})}\ge inj_{(M,g)}>0$. Moreover, Theorem $10.7.1$ of \cite{BuragoBuragoIvanov01} permits to conclude that $(M_{\infty}, g_{\infty})$ have sectional curvature bounded below by $\Lambda_1$. On the other hand, the property of being a \textit{metric space of curvature $\le K$} (see Definition $1.2$ at page $159$ of 
 \cite{BridsonHaefligger}, i.e., being locally a $Cat(K)$ space) pass to the limit in the pointed Gromov-Hausdorff convergence 
 because distances pass to the Gromov-Hausdorff limit. This fact combined with Theorem $1A.6$ at page $173$ of 
 \cite{BridsonHaefligger} implies that for a smooth Riemannian manifold to have sectional curvature bounded above by $K$ is 
 equivalent to satisfy the condition of having curvature bounded above in the sense of Alexandrov that is in the sense of Definition 
 $1.2$ of page $159$ of \cite{BridsonHaefligger}. From this easily follows that the sectional curvature of $(M_{\infty},g_{\infty})$ 
 (that exists because $g_{\infty}$ is assumed at least $C^2$ or more regular by the assumption of strong bounded geometry 
 smooth at infinity) is bounded from above by the same constant than the sectional curvature of $M$. Hence $(M_{\infty},g_{\infty})
 $ have strong bounded geometry, in particular have also mild bounded geometry which gives the validity of Theorem 
 \ref{Res:Theorem1} and Lemma \ref{Lemma:smallvolumesimpliessmalldiametersmild} in $(M_{\infty},g_{\infty})$, with constants 
 $d(M_{\infty},g_{\infty})$ and $v^*(M_{\infty}, g_{\infty})$ that in principle depend on $(M_{\infty},g_{\infty})$. On the other hand by \eqref{Eq:Theorem1DiameterEstimates} and the explicit estimates that Lemma \ref{Lemma:smallvolumesimpliessmalldiametersmild} on $v^*$ we argue that $d(M_{\infty},g_{\infty})\ge d(M,g)>0$ and $v^*(M_{\infty}, g_{\infty})\ge v^*(M, g)>0$. 
% a careful analysis of the proof of Proposition \ref{Proposition:Sufficient} implies that $r_\varepsilon^*(M_{\infty},g_{\infty})\ge r^*_\varepsilon(M,g)>0$ from which easily follows $d(M_{\infty},g_{\infty})\ge d(M,g)>0$ as shown by the constant displayed in \eqref{Eq:Theorem1DiameterEstimates} of Theorem \ref{Res:Theorem1} shows that such a $d$ is uniformly bounded from below with respect to all $(M_{\infty},g_{\infty})$ limit manifolds, because the condition of having a lower bound depends continuously on the supremum of the scalar curvature function.
  At this stage an application of Lemma \ref{Lemma:smallvolumesimpliessmalldiametersmild} gives that for $v\le v^*$ any generalized isoperimetric region $\tilde\Omega_v\subseteq M_{\infty}$ with $V_{g_{\infty}}(\tilde\Omega_v)=v$ have $diam_{g_{\infty}}(\tilde\Omega_v)\le\mu^*v^{\frac1n}$. Thus for small values of $v\leq\tilde v_0:=\min\{v^*, \left(\frac{d}{\mu^*}\right)^n\}$ we have $\diam_{g_{\infty}}(\tilde\Omega_v)\le d$, where $d$ is given by Theorem $\ref{Res:Theorem1}$. Finally for every finite perimeter set  $\Omega\subset M$ such that $V_g(\Omega)=v$, we conclude that
\begin{equation}\label{Eq:ResTheorem3Last}
\P_g(\Omega)\ge I_{M,g}(v)=I_{M_{\infty},g_{\infty}}(v)=\P_{g_{\infty}}(\tilde\Omega_v)>\P_{g_0}(B)=I_{\mathbb{M}^n_{k_0}}(v),
\end{equation} 
where the first inequality comes from the definition of $I_M$, the first equality comes from Theorem $1$ of \cite{NarAsian} where $\Omega_v$ is a generalized isoperimetric region of $V_{g_{\infty}}(\Omega_v)=v$, the second equality from the definition of $\Omega_v$ as a generalized isoperimetric region of volume $v$, the second inequality is an application of Theorem \ref{Res:Theorem1} to $(M_{\infty},g_{\infty})$, and the last is simply the fact that the isoperimetric regions in space forms are the geodesic balls.  This finish the proof of Theorem  \ref{Res:Theorem3}.
\end{proof}
We have now all the tools needed to prove Theorem \ref{Res:Theorem3.1}.
\begin{proof}[Proof of Theorem \ref{Res:Theorem3.1}] Analogously to the second proof of Proposition \ref{Proposition:SufficientWeak1} we deform smoothly the metric $g$ along the Ricci flow $\tilde g_t$ with initial data $\tilde g_0=g$. In first we suppose that $k_0$ is also a strict upper bound for the sectional curvature of $(M,g)$, i.e., $Sec_g< k_0$. Then for any $k_0$ such that $n(n-1)k_0>S_g$ we apply Theorem \ref{Res:Theorem3} to some manifold $(M,\tilde g_t)$ since these manifolds satisfy the required hypothesis of regularity and bounded geometry at infinity. Moreover, thanks to the estimates on the supremum of the sectional curvatures of $(M,\tilde g_t)$ given by the Proposition at page $260$ of \cite{Kapovitch}  reported here in Theorem \ref{thm:kapovitch1}, there exists a sufficiently small $t$ such that $S_{\tilde g_t}<n(n-1)k_0$. Hence applying Theorem \ref{Res:Theorem3} to $(M,\tilde g_t)$ we have that there exists a $v^*=v^*(n,k,k_0, S_g)>0$ such that $\tilde{v}_0(M,\tilde g_t)\ge v^*$ for every $t\in[0, T]$ such that for every $\Omega$ with $V_{\tilde g_t}(\Omega)\le v^*(n,k,v_0, \tilde g_t)$ we have
\begin{equation}\label{Eq:ResTheorem3.10}
\P_{\tilde g_t}(\Omega)>\P_{g_{k_0}}(B_{k_0}).
\end{equation}
Consider now a finite perimeter set $\Omega$ such that $V_{g}(\Omega)<v^*$. By $C^0$ convergence of $\tilde g_t$ to $\tilde g_0=g$ as stated in \eqref{Eq:Shi}, we can find $t$ small enough such that $V_{\tilde g_t}(\Omega)<v^*$ and thus \eqref{Eq:ResTheorem3.10} hold for all sufficiently small $t$. Again by $C^0$ convergence of $\tilde g_t$ to $\tilde g_0=g$, taking the limit when $t\to 0^+$ in \eqref{Eq:ResTheorem3.10} we obtain 
\begin{equation} \label{Eq:ResTheorem3.11}
\P_g(\Omega)=\P_{\tilde g_0}(\Omega)=\lim_{t\to0^+}\P_{\tilde g_t}(\Omega)\ge\P_{g_{k_0}}(B_{k_0}).
\end{equation}
%The only thing that remains to prove is that there exists a uniform lower bound the volumes $v^*(n,k,v_0, \tilde g_t)$ w.r.t. $t$. 
To settle the general case we need just to observe that it is not too hard to obtain (from the proof of Proposition at page $260$ of \cite{Kapovitch}, i.e., Theorem \ref{thm:kapovitch1}) the following estimates on the bounds of the scalar curvature is obtained, i.e.,
\begin{eqnarray}\label{Eq:KapovitchScalar}
 \inf_{x\in M}\left\{ Sc_{\tilde {g}_t}(x)\right\}-C(n,T)t\le Sc_{\tilde g_t}(x) \le \sup_{x\in M} \left\{ Sc_{\tilde {g}_t}(x)\right\} +C(n,T)t.
\end{eqnarray}
With \eqref{Eq:KapovitchScalar} in mind the same argument used in the special case is valid again in the general case. So \eqref{Eq:ResTheorem3.11} is proved in the general case. To finish the proof of the theorem we remark that using just \eqref{Eq:ResTheorem3.11} the asymptotic expansion of the isoperimetric profile given in Corollary \ref{CorRes:PuiseuxSeries}  holds and then reasoning by reduction to the absurd we conclude that there exists $\tilde v_1=\tilde v_1(n,k,k_0, inj_M, S_g, I_M)\le\tilde v_0$ such that for every volume $v\le\tilde v_1$ the equality case in \eqref{Eq:ResTheorem3.11} cannot happen. Unfortunately we do not still have an effective estimates for $\tilde v_1$. From this the theorem follows.
\end{proof}
\section{Asymptotic expansion of the isoperimetric profile in strong bounded geometry}
We prove in this last section the asymptotic expansion in Puiseux's series up to the second nontrivial term stated in Corollary \ref{CorRes:PuiseuxSeries}.
\begin{proof}[Proof of Corollary \ref{CorRes:PuiseuxSeries}] We use Theorem \ref{Res:Theorem3} to prove the first of the following inequalities then we compare with the area of a geodesic ball centered at $x_0$ and of enclosed volume $v$, proving the second inequality of
\begin{tiny}
\begin{eqnarray}\nonumber
%c_nv^{\frac{(n-1)}n}\left(1-\gamma_nS_gv^{\frac2n}\right) & \le & \lim_{(k_0, v)\to(\frac{S_g}{n(n-1)}^+,0^+)}\nonumber
P_{g_{\frac{S_g}{n(n-1)}}}\left(B_{g_{\frac{S_g}{n(n-1)}},v}\right) & = & c_nv^{\frac{(n-1)}n}\left(1-\gamma_nS_gv^{\frac2n}\right)+O_{\frac{S_g}{n(n-1)}}\left(v^{\frac4n}\right)\\ \nonumber
& < & P_{g_{k_0}}\left(B_{g_{k_0},v}\right)\\ \nonumber
 & = & c_nv^{\frac{(n-1)}n}\left(1-\gamma_nn(n-1)k_0v^{\frac2n}\right)+O_{k_0}\left(v^{\frac4n}\right)\\ \label{Eq:PuiseuxSeries}
 & \le & I_M(v)\\ \label{Eq:PuiseuxSeries0} 
& \le & c_nv^{\frac{(n-1)}n}\left(1-\gamma_nSc_g(x_0)v^{\frac2n}\right)+O_{x_0}\left(v^{\frac4n}\right)\\ \nonumber
& = & P_g(B_g(x_0, v)),%\\
%& \le & c_nv^{\frac{(n-1)}n}\left(1-\gamma_nS_gv^{\frac2n}\right),
\end{eqnarray}
\end{tiny}
for every $x_0\in M$, $k_0>S_g$, $v\le\tilde v_0$, with $B_{g_{k_0},v}$ a ball in $\mathbb{M}^n_{k_0}$ such that $V_{g_{k_0}}(B_{k_0,v})=v$ and $B_g(x_0, v)$ a ball of $(M^n,g)$ having $V_g\left(B_g(x_0, v)\right)=v$.
%To prove uniformity with respect to $x_0$ we use all the theory developed by the first author in \cite{NarAnn}, \cite{NarCalcVar}. 
Taking the infimum with respect to $x_0$ and using \eqref{Eq:PuiseuxSeries} and \eqref{Eq:PuiseuxSeries0} we get for a fixed $\tilde x_0\in M$, that $$O_{\tilde x_0}\left(v^{\frac4n}\right)\ge\eta_1(v):=\inf_{x_0\in M}O_{x_0}\left(v^{\frac4n}\right)\ge O_{\frac{S_g}{n(n-1)}}\left(v^{\frac4n}\right).$$ Hence
\begin{tiny}
\begin{eqnarray}
c_nv^{\frac{(n-1)}n}\left(1-\gamma_nS_gv^{\frac2n}\right)+\eta_2(v) & \le & I_M(v)\\ 
& \le & c_nv^{\frac{(n-1)}n}\left(1-\gamma_nS_gv^{\frac2n}\right)+\eta_1(v),
\end{eqnarray}
\end{tiny}
with $\eta_1(v)=O\left(v^{\frac4n}\right)$, $\eta_2(v)=O\left(v^{\frac4n}\right)$. We need that the metric $g$ on the manifold $M$ is at least $C^3$ because of the Taylor expansion of the metric in normal coordinates. This ensures that $O_{\tilde x_0}\left(v^{\frac4n}\right)$ does not blows up. From this the corollary, indeed follows promptly. 
\end{proof}
\markboth{References}{References}
      \bibliographystyle{alpha}
      \bibliography{ArtigoDruetNonCompact6}

\begin{thebibliography}{BMOR84}

\bibitem[AFP00]{AmbrosioFuscoPallara}
Luigi Ambrosio, Nicola Fusco, and Diego Pallara.
\newblock {\em Functions of bounded variation and free discontinuity problems}.
\newblock Oxford mathematical monographs. Clarendon Press, Oxford, New York,
  2000.
\newblock Autres tirages : 2006.

\bibitem[AL99]{AubinLi}
Thierry Aubin and Yan~Yan Li.
\newblock On the best {S}obolev inequality.
\newblock {\em J. Math. Pures Appl. (9)}, 78(4):353--387, 1999.

\bibitem[Aub76]{AubinJDG76}
Thierry Aubin.
\newblock Probl\`emes isop\'erim\'etriques et espaces de {S}obolev.
\newblock {\em J. Differential Geometry}, 11(4):573--598, 1976.

\bibitem[BBI01]{BuragoBuragoIvanov01}
Dmitri Burago, Yuri Burago, and Sergei Ivanov.
\newblock {\em A course in metric geometry}, volume~33 of {\em Graduate Studies
  in Mathematics}.
\newblock American Mathematical Society, Providence, RI, 2001.

\bibitem[BH99]{BridsonHaefligger}
Martin~R. Bridson and Andr{\'e} Haefliger.
\newblock {\em Metric spaces of non-positive curvature}, volume 319 of {\em
  Grundlehren der Mathematischen Wissenschaften [Fundamental Principles of
  Mathematical Sciences]}.
\newblock Springer-Verlag, Berlin, 1999.

\bibitem[Bie03]{RodneyThesis}
Rodney~Josue Biezuner.
\newblock Best constants, optimal {S}obolev inequalities on {R}iemannian
  manifolds and its applications.
\newblock {\em Rutgers Ph.D. Thesis}, pages 1--131, 2003.

\bibitem[BK13]{KloeckerKuperberg}
Kloeckner Beno\^{i}t and Greg Kuperberg.
\newblock The cartan-hadamard conjecture and the little prince.
\newblock {\em arXiv:1303.3115}, pages 1--31, 2013.
\newblock
  \href{https://arxiv.org/pdf/1303.3115.pdf}{https://arxiv.org/pdf/1303.3115.pdf}.

\bibitem[BM82]{BerardMeyer}
Pierre B{\'e}rard and Daniel Meyer.
\newblock In\'egalit\'es isop\'erim\'etriques et applications.
\newblock {\em Ann. Sci. \'Ecole Norm. Sup. (4)}, 15(3):513--541, 1982.

\bibitem[BMOR84]{BMOR}
Josef Bemelmans, Min-Oo, and Ernst~A. Ruh.
\newblock Smoothing {R}iemannian metrics.
\newblock {\em Math. Z.}, 188(1):69--74, 1984.

\bibitem[Cha06]{Chavel}
Isaac Chavel.
\newblock {\em Riemannian geometry}, volume~98 of {\em Cambridge Studies in
  Advanced Mathematics}.
\newblock Cambridge University Press, Cambridge, second edition, 2006.
\newblock A modern introduction.

\bibitem[Cro84]{Croke}
Christopher~B. Croke.
\newblock A sharp four-dimensional isoperimetric inequality.
\newblock {\em Comment. Math. Helv.}, 59(2):187--192, 1984.

\bibitem[Dru98]{DruetJFA}
Olivier Druet.
\newblock Optimal {S}obolev inequalities of arbitrary order on compact
  {R}iemannian manifolds.
\newblock {\em J. Funct. Anal.}, 159(1):217--242, 1998.

\bibitem[Dru99]{DruetMathAnn1999}
Olivier Druet.
\newblock The best constants problem in {S}obolev inequalities.
\newblock {\em Math. Ann.}, 314(2):327--346, 1999.

\bibitem[Dru00]{DruetPRSEdinburgh}
Olivier Druet.
\newblock Generalized scalar curvature type equations on compact {R}iemannian
  manifolds.
\newblock {\em Proc. Roy. Soc. Edinburgh Sect. A}, 130(4):767--788, 2000.

\bibitem[Dru02a]{DruetGeoDed2002}
Olivier Druet.
\newblock Isoperimetric inequalities on compact manifolds.
\newblock {\em Geom. Dedicata}, 90:217--236, 2002.

\bibitem[Dru02b]{DruetPAMS}
Olivier Druet.
\newblock Sharp local isoperimetric inequalities involving the scalar
  curvature.
\newblock {\em Proc. Amer. Math. Soc.}, 130(8):2351--2361 (electronic), 2002.

\bibitem[Dru10]{DruetSurvey}
Olivier Druet.
\newblock Isoperimetric inequalities on nonpositively curved spaces.
\newblock 2010.
\newblock
  \href{http://math.arizona.edu/~dido/presentations/Druet-Carthage.pdf}{http://math.arizona.edu/~dido/presentations/Druet-Carthage.pdf}.

\bibitem[Giu84]{Giusti}
Enrico Giusti.
\newblock {\em Minimal surfaces and functions of bounded variation}, volume~80
  of {\em Monographs in Mathematics}.
\newblock Birkh\"auser Verlag, Basel, 1984.

\bibitem[GMT83]{Tamanini}
E.~Gonzalez, U.~Massari, and I.~Tamanini.
\newblock On the regularity of boundaries of sets minimizing perimeter with a
  volume constraint.
\newblock {\em Indiana Univ. Math. J.}, 32(1):25--37, 1983.

\bibitem[GR13]{RitoreGalli}
Matteo Galli and Manuel Ritor{\'e}.
\newblock Existence of isoperimetric regions in contact sub-{R}iemannian
  manifolds.
\newblock {\em J. Math. Anal. Appl.}, 397(2):697--714, 2013.

\bibitem[Ham82]{Hamilton}
Richard~S. Hamilton.
\newblock Three-manifolds with positive {R}icci curvature.
\newblock {\em J. Differential Geom.}, 17(2):255--306, 1982.

\bibitem[Heb99]{Hebey}
Emmanuel Hebey.
\newblock {\em Nonlinear analysis on manifolds: {S}obolev spaces and
  inequalities}, volume~5 of {\em Courant Lecture Notes in Mathematics}.
\newblock New York University, Courant Institute of Mathematical Sciences, New
  York; American Mathematical Society, Providence, RI, 1999.

\bibitem[Heb02]{Hebey02}
Emmanuel Hebey.
\newblock Sharp {S}obolev-{P}oincar\'e inequalities on compact {R}iemannian
  manifolds.
\newblock {\em Trans. Amer. Math. Soc.}, 354(3):1193--1213 (electronic), 2002.

\bibitem[HK78]{HeintzeKarcher}
Ernst Heintze and Hermann Karcher.
\newblock A general comparison theorem with applications to volume estimates
  for submanifolds.
\newblock {\em Ann. Sci. \'Ecole Norm. Sup. (4)}, 11(4):451--470, 1978.

\bibitem[Kap05]{Kapovitch}
Vitali Kapovitch.
\newblock Curvature bounds via {R}icci smoothing.
\newblock {\em Illinois J. Math.}, 49(1):259--263 (electronic), 2005.

\bibitem[Kle92]{Kleiner92}
Bruce Kleiner.
\newblock An isoperimetric comparison theorem.
\newblock {\em Invent. Math.}, 108(1):37--47, 1992.

\bibitem[Lio84]{LionsCCIAIHP}
P.-L. Lions.
\newblock The concentration-compactness principle in the calculus of
  variations. {T}he locally compact case. {I}.
\newblock {\em Ann. Inst. H. Poincar\'e Anal. Non Lin\'eaire}, 1(2):109--145,
  1984.

\bibitem[Lio85]{LionsCCIRMIA}
P.-L. Lions.
\newblock The concentration-compactness principle in the calculus of
  variations. {T}he limit case. {I}.
\newblock {\em Rev. Mat. Iberoamericana}, 1(1):145--201, 1985.

\bibitem[Mag12]{Maggi}
Francesco Maggi.
\newblock {\em Sets of finite perimeter and geometric variational problems},
  volume 135 of {\em Cambridge Studies in Advanced Mathematics}.
\newblock Cambridge University Press, Cambridge, 2012.
\newblock An introduction to geometric measure theory.

\bibitem[MFN15]{FloresNardulli015}
Abraham~Enrique Munoz~Flores and Stefano Nardulli.
\newblock Continuity and differentiability properties of the isoperimetric
  profile in complete noncompact riemannian manifolds with bounded geometry.
\newblock {\em arXiv:1404.3245}, V3:1--31, 2015.
\newblock
  \href{https://arxiv.org/pdf/1404.3245.pdf}{https://arxiv.org/pdf/1404.3245.pdf}.

\bibitem[MJ00]{MJ}
Frank Morgan and David~L. Johnson.
\newblock Some sharp isoperimetric theorems for {R}iemannian manifolds.
\newblock {\em Indiana Univ. Math. J.}, 49(2):1017--1041, 2000.

\bibitem[MN15]{FloresNardulliGenComp}
A.~E. {Mu{\~n}oz Flores} and S.~{Nardulli}.
\newblock {Generalized compactness for finite perimeter sets and applications
  to the isoperimetric problem}.
\newblock {\em ArXiv e-prints}, April 2015.

\bibitem[MN16]{MondinoNardulli}
A.~Mondino and S.~Nardulli.
\newblock Existence of isoperimetric regions in non-compact {R}iemannian
  manifolds under {R}icci or scalar curvature conditions.
\newblock {\em Comm. Anal. Geom.}, 24(1):115--1138, 2016.

\bibitem[MPPP07]{MPPP}
M.~Miranda, Jr., D.~Pallara, F.~Paronetto, and M.~Preunkert.
\newblock Heat semigroup and functions of bounded variation on {R}iemannian
  manifolds.
\newblock {\em J. Reine Angew. Math.}, 613:99--119, 2007.

\bibitem[Nar14a]{NarAsian}
Stefano Nardulli.
\newblock Generalized existence of isoperimetric regions in non-compact
  {R}iemannian manifolds and applications to the isoperimetric profile.
\newblock {\em Asian J. Math.}, 18(1):1--28, 2014.

\bibitem[Nar14b]{NarCalcVar}
Stefano Nardulli.
\newblock The isoperimetric profile of a noncompact {R}iemannian manifold for
  small volumes.
\newblock {\em Calc. Var. Partial Differential Equations}, 49(1-2):173--195,
  2014.

\bibitem[Nar15]{NarBBMS}
Stefano Nardulli.
\newblock Regularity of solutions of the isoperimetric problem that are close
  to a smooth manifold.
\newblock {\em arXiv:0710.1849v3}, 2015.

\bibitem[Nik91]{Nikolaev}
I.~G. Nikolaev.
\newblock Bounded curvature closure of the set of compact {R}iemannian
  manifolds.
\newblock {\em Bull. Amer. Math. Soc. (N.S.)}, 24(1):171--177, 1991.

\bibitem[Pet87]{Peters}
Stefan Peters.
\newblock Convergence of {R}iemannian manifolds.
\newblock {\em Compositio Math.}, 62(1):3--16, 1987.

\bibitem[Pet06]{Petersen}
Peter Petersen.
\newblock {\em Riemannian geometry}, volume 171 of {\em Graduate Texts in
  Mathematics}.
\newblock Springer, New York, second edition, 2006.

\bibitem[Rit05]{RitoreCartanHadamard05}
Manuel Ritor{\'e}.
\newblock Optimal isoperimetric inequalities for three-dimensional
  {C}artan-{H}adamard manifolds.
\newblock In {\em Global theory of minimal surfaces}, volume~2 of {\em Clay
  Math. Proc.}, pages 395--404. Amer. Math. Soc., Providence, RI, 2005.

\bibitem[RS10]{RitoreSinestrari}
Manuel Ritor{\'e} and Carlo Sinestrari.
\newblock {\em Mean curvature flow and isoperimetric inequalities}.
\newblock Advanced Courses in Mathematics. CRM Barcelona. Birkh\"auser Verlag,
  Basel, 2010.
\newblock Edited by Vicente Miquel and Joan Porti.

\bibitem[Sch08]{Schulze08}
Felix Schulze.
\newblock Nonlinear evolution by mean curvature and isoperimetric inequalities.
\newblock {\em J. Differential Geom.}, 79(2):197--241, 2008.

\bibitem[Ser64]{Serrin}
James Serrin.
\newblock Local behavior of solutions of quasi-linear equations.
\newblock {\em Acta Math.}, 111:247--302, 1964.

\bibitem[Shi89]{Shi}
Wan-Xiong Shi.
\newblock Deforming the metric on complete {R}iemannian manifolds.
\newblock {\em J. Differential Geom.}, 30(1):223--301, 1989.

\bibitem[Str08]{Struwe}
Michael Struwe.
\newblock {\em Variational methods}, volume~34 of {\em Ergebnisse der
  Mathematik und ihrer Grenzgebiete. 3. Folge. A Series of Modern Surveys in
  Mathematics [Results in Mathematics and Related Areas. 3rd Series. A Series
  of Modern Surveys in Mathematics]}.
\newblock Springer-Verlag, Berlin, fourth edition, 2008.
\newblock Applications to nonlinear partial differential equations and
  Hamiltonian systems.

\bibitem[Tal76]{Talenti}
Giorgio Talenti.
\newblock Best constant in {S}obolev inequality.
\newblock {\em Ann. Mat. Pura Appl. (4)}, 110:353--372, 1976.

\bibitem[Wei26]{Weil26}
Andr{\'e} Weil.
\newblock Sur les surfaces {\`a} courbure n\'egative.
\newblock {\em C.R. Acad. Sci. Paris}, 182(2):1069--1071, 1926.

\end{thebibliography}
      \addcontentsline{toc}{section}{\numberline{}References}
      \emph{Stefano Nardulli\\ Departamento de Matem\'atica\\ Instituto de Matem\'atica\\ UFRJ-Universidade Federal do Rio de Janeiro, Brazil\\ email: \href{mailto:nardulli@im.ufrj.br}{nardulli@im.ufrj.br}\\ \\ Luis Eduardo Osorio Acevedo\\ Ph.D. student\\Instituto de Matem\'atica\\UFRJ-Universidade Federal do Rio de Janeiro, Brazil\\ email: \href{mailto:mat.eduardo@gmail.com}{mat.eduardo@gmail.com}} 
\end{document}